\documentclass[10pt,a4paper]{article} 
\usepackage[utf8]{inputenc}
\usepackage{amsmath}
\usepackage{amsfonts}
\usepackage{amssymb}
\usepackage{amsthm}
\usepackage{dsfont}
\usepackage{float}
\usepackage{url}
\PassOptionsToPackage{hyphens}{url}
\usepackage{mathtools}
\usepackage{pdflscape}
\usepackage{multirow}
\usepackage{enumitem}

\usepackage{graphicx}
\usepackage{grffile}

\usepackage[noend, ruled, linesnumbered]{algorithm2e}
\makeatletter
\newcommand{\nosemic}{\renewcommand{\@endalgocfline}{\relax}}
\newcommand{\dosemic}{\renewcommand{\@endalgocfline}{\algocf@endline}}
\let\oldnl\nl
\newcommand{\nonl}{\renewcommand{\nl}{\let\nl\oldnl}}
\makeatother

\usepackage[left=3cm,right=3cm,top=2cm,bottom=2cm]{geometry}

\usepackage{url}

\usepackage[colorlinks,
            citecolor=blue,
            linkcolor=red,
            urlcolor=blue]{hyperref}

\usepackage{bm}

\usepackage[authoryear, round]{natbib}

\usepackage{booktabs}
\usepackage{colortbl}
\usepackage{graphicx}

\usepackage{tikz}
\usetikzlibrary{shapes,snakes}
\usetikzlibrary{positioning, decorations.pathreplacing}
\usepackage{varwidth}
\tikzset{>=latex}

\newcommand{\co}{c^{\boldsymbol o}}
\newcommand{\ci}{c^{\boldsymbol i}}
\newcommand{\cw}{c^{\boldsymbol w}}
\newcommand{\cm}{c^{\boldsymbol m}}
\newcommand{\E}{\mathbb{E}}
\renewcommand{\vec}[1]{\boldsymbol{#1}}

\newcommand{\changed}[1]{#1}

\setitemize{labelindent=0.2em,labelsep=0.5em,leftmargin=1em,itemsep=0mm,topsep=0.5em}
\setenumerate{labelsep=0.5em,leftmargin=1.5em,itemsep=0pt,topsep=0.5em}

\usepackage{xspace}
\newcommand{\EESPlong}{\textsc{Elective and Emergency Surgery Planning Problem}\xspace}
\newcommand{\EESP}{EESPP\xspace}
\newcommand{\APPlong}{\textsc{Advanced Planning and Appointment Scheduling Problem}\xspace}
\newcommand{\APP}{APASP\xspace}
\newcommand{\OSPlong}{\textsc{Online Scheduling Problem}\xspace}
\newcommand{\OSP}{OSP\xspace}

\usepackage{authblk}
\title{Surgery Scheduling in Flexible Operating Rooms by using a Convex Surrogate Model of Second-Stage Costs}
\author{Mohammed Majthoub Almoghrabi\footnote{corresponding author}\ \ 
and Guillaume Sagnol\footnote{supported by the Deutsche Forschungsgemeinschaft (DFG, German Research Foundation) under Germany's Excellence Strategy --- The Berlin Mathematics Research Center MATH+ (EXC-2046/1, project ID: 390685689).}\\
\texttt{\{majthoub,sagnol\}@math.tu-berlin.de}}
\affil{Technische Universität Berlin, Insitut für Mathematik\\
Stra\ss e des 17 Juni 136, 10623 Berlin, Germany}
\date{}

\usepackage{setspace}

\begin{document}
\maketitle
\setcounter{page}{1}

\begin{abstract}
We study the elective surgery planning problem in a hospital with operating rooms shared by elective and emergency patients. 
This problem is split in two distinct phases. 
First, a subset of patients to be operated in the next planning
period is selected and the selected patients are assigned to a block and a tentative starting time. 
Then, in the online phase of the problem, a policy decides how to insert the emergency patients in the schedule and may cancel
planned surgeries. The overall goal is to minimize the expectation of a cost function representing the assignment of patient to blocks, case cancellations, overtime, waiting time and idle time.
We model the offline problem by a two-stage stochastic program,
and show that the optimal second-stage costs can be approximated by a
convex piecewise linear surrogate model that can be computed
in a preprocessing step. This results in a mixed integer program which can be solved very fast, even for very
large instances of the problem. We also describe a greedy 
policy for the online phase of the problem, and
analyze the performance of our approach by comparing it to
both heuristic methods
or approaches relying on sampling average approximation (SAA)
on a large set of benchmarking instances. Our simulations indicate that our approach can reduce
the expected costs by as much as \changed{30}\% compared to heuristic
methods and it can solve problems with 1000 patients in about one minute, while SAA-approaches fail to obtain good solutions within 30 minutes on small instances.
\end{abstract}

\medskip
\noindent\textbf{Keywords:} OR in health services, Surgery Scheduling, Two-stage Stochastic Programming

\onehalfspacing

\section{Introduction}

The healthcare industry is currently facing a crisis as hospital
surgical services are being stretched to their limits.
Operating rooms (ORs) are one of the most critical resources,
concentrating soaring costs and an increasing demand.
Aggregate surgical expenditures amount to 40\% of hospital costs~\citep{pham2008surgical}, which represents, e.g., 
5 to 7\% of the gross national product in the United States, and  this share is constantly growing~\citep{munoz2010national}.
The demand is so high that many patients have to wait months or even years for an elective surgery, with large variations across countries~\citep{oecd2020waiting}. 
The situation was further aggravated by the COVID-19 pandemic, which led hospitals to cancel planned operations,
thus increasing the elective surgery backlog to record levels in many countries~\citep{mehta2022elective}.
This is not only dramatic from an organizational point of view, but also from a clinical one.
Indeed, long waiting times are not only a discomfort for the patient, but also a risk factor that can jeopardize the success of an operation~\citep{garbuz2006delays}.

\begin{figure}[t]
 \begin{center}
  \includegraphics[width=\textwidth]{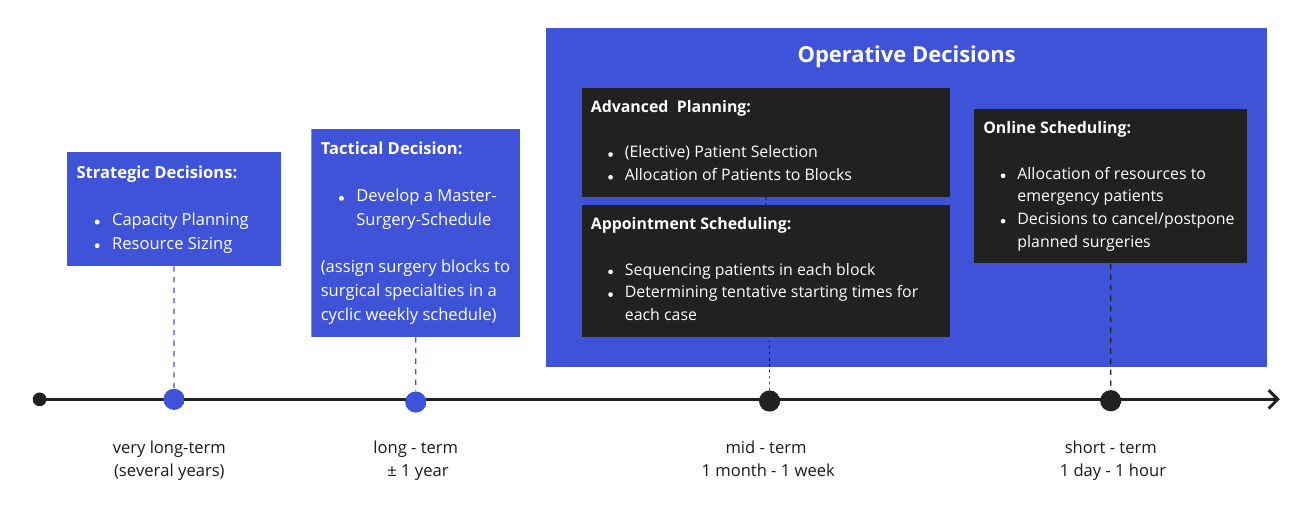}
 \end{center}
\caption{\small Standard breakdown of OR management tasks into successive decision stages. This paper focuses on operative decisions (advanced planning, appointment scheduling and online scheduling) \label{fig:decisions}}
\end{figure}

In this context, optimizing the utilization of operating rooms is a priority of hospital management~\citep{dexter2002schedule}.  Appropriate
decision-making must occur at all stages of the management process~\citep{guerriero2011operational},
as depicted in Figure~\ref{fig:decisions}. Strategical decisions are made on the very long-term, and are typically concerned with the sizing and staffing of the facilities; In the most common OR-management method called \emph{block-scheduling},
the tactical decisions involve establishing a \emph{master surgery schedule} (MSS), which assigns blocks of OR-time to the different surgical specialties of the hospital in a cyclic weekly schedule. An MSS is used for a longer period ranging from several months to a year, and has to take various constraints into account, such as the required volume of OR-time for each specialty, the availability of staff during the week and the impact of operations on downstream recovery units~\citep{fugener2014master}. 

\changed{This paper focuses on the mid-term decision stage, dealing with uncertainties related to both the surgical duration of an operation and the arrival of emergency patients. 
This includes the 
\emph{Advanced Planning}, which
requires to select patients to be operated during the next planning period (typically one week), 
and to assign them to blocks of the MSS, and the 
\emph{Appointment Scheduling}, in which each
patient must be assigned a tentative starting
time for the operation.
Although the primary focus of this paper lies in 
mid-term decisions, it is essential to bear in mind that the quality of a planning strategy can only be assessed if we 
consider the operational decisions as a whole.
Decisions that can be taken to
manage process uncertainty during execution,
i.e., insertion of 
add-on cases (emergencies) in the schedule or rescheduling/cancelling planned operations, 
can have a major impact on costs.
}

\paragraph{Our contribution.}
The problem studied in this article was the subject of the 
Optimization Modeling Competition
AIMMS-MOPTA 2022~\citep{MOPTA22}.
This paper aims at presenting the approach of the authors to solve this problem\footnote{The authors were part of a team who ranked third in this competition, together with Przemyslav Bartman.}.
Our proposed formalization of this problem will be called the \EESPlong (\EESP) throughout.
The \EESP is an integrated approach combining the two distinct phases of the operative stage, cf.\ Figure~\ref{fig:decisions}:
\begin{itemize}
\item An offline phase in which a problem called \APPlong (\APP) is solved.
\item An online phase in which the uncertainty is to be handled by means of a policy: this is the
\OSPlong (\OSP).
\end{itemize}
\changed{Our contributions can be summarized as follows:}
\begin{enumerate}
  \item To the best of our knowledge, the advanced planning problem (selection of patients and assignment to blocks) and the appointment scheduling problem (sequencing and decision of tentative starting times for each patient in each block) have always been treated separately in the literature. In contrast, our approach solves these two problems at once, and attempts to reserve some capacity for future ---but yet unknown--- emergency cases.
  \item \changed{We model the problem to be solved in the offline phase as} a two-stage stochastic program,
  \changed{which we approximate with a compact mixed integer program (MIP)} 
  relying on a piecewise linear surrogate model of the optimal second-stage costs, that can be computed in a preprocessing step. It runs very fast and yields better results than both simple heuristic approaches as well as state-of-the-art sample average approximation (SAA) approaches.
  \item  Our computational study based on benchmarking instances
shows that our approach can improve the costs by up to
\changed{30}\% compared to heuristic approaches such as
First-Fit, or a deterministic approach in which
the random duration of each surgery is replaced by 
its mean value. More advanced methods relying on SAA
to handle uncertainty (integrated MIP formulation obtained
by sampling scenarios or Benders decomposition)
only beat our approach on very small instances,
even when {emergencies} are ignored.
\item Unlike previous work on the advanced planning problem, our approach explicitly takes the cost related to idle time and waiting time between surgeries into account. Other approaches not considering the scheduled starting times could not do this, and thus ignored the
costs related to turnovers between surgeries, even though minimizing such costs is one of the priorities of hospital management~\citep{dexter2000statistical, dexter2002schedule}.
\item We propose a greedy heuristic policy which can be used to handle emergency patients in the online phase of this problem
\changed{and to measure the quality of a solution to the \APP}.
\item An implementation of our solution 
is available online,
along with a simulator that can be used to compare the performance of different policies under different scenarios. The online tool to compute and visualize surgery schedules for random instances generated using the competition data~\citep{MOPTA22} is presented at \url{https://wsgi.math.tu-berlin.de/esp_demonstrator/};
see Appendix~\ref{sec:demonstrator}.
\end{enumerate}

\vspace{-1em}
\changed{
\paragraph{Restrictions.}
We assume that the set of emergency patients arriving on a particular day is revealed at the beginning of the day, rather than in a continuous process. This allows us to handle one of the rules of the
AIMMS-MOPTA competition without ambiguity: all emergency patients must be handled on the day of their arrival, without having to consider ORs dedicated to emergencies that are open at night.
We point out that this assumption remains fairly realistic,
since in practice the majority of emergencies are in fact \emph{soft emergencies}, which must be operated within 12 to 24 hours after their arrival~\citep{Charite-personal}.
Another restriction of our paper is that we do not consider
restriction on the availability of downstream resources such
as recovery beds.
}

\paragraph{Outline.} This paper is organized as follows.  Section~\ref{sub:relwork} reviews the relevant literature. In Section~\ref{sec:model}, we define the problem formally and discuss
the assumptions of this study. 
Then, we present our approach to solve
the offline problem (\APP) in Section~\ref{sec:offlinephase}
and the greedy policy we used for the online phase (\OSP)
in {Subsection}~\ref{sec:onlinephase}. In addition, we introduce four alternative approaches to solve
the \APP in Section~\ref{sec:other_approaches}: these methods
are used for the sake of comparison in our computational study in Section~\ref{sec:study}.
Finally, we conclude and discuss
possible future work in
Section~\ref{sec:conclusion}.

\section{Related work}
\label{sub:relwork}

There has been a long stream of research on the elective surgery planning problem, 
so in this section we focus only on the literature most relevant to our problem.
For comprehensive reviews about the problems arising in the area of surgery scheduling,
as well as optimization approaches to compute planning strategies
and interesting challenges for future directions, we refer
to~\citet{may2011surgical,cardoen2010operating,shehadeh2022stochastic}.

\paragraph{Advanced planning problem.}
Recall that the goal of the advanced planning is to select a subset of elective patients to be operated in the next planning period, and to assign them to blocks of the MSS. 
This problem is usually rewritten as a two-stage optimization problem,
in which the first-stage variables represent the assignment of patients to blocks (and sometimes, the decision to open a block), while the second stage variables are used to
represent the costs associated with overtime (exceeding the regular working time of a block)
and undertime (not using a block to its full capacity). However, the approaches differ widely depending on the nature of the considered constraints and the approach used to handle uncertainty (stochastic programming or robust optimization approaches).

\medskip
\emph{Two-Stage stochastic programming.} One of the most natural {ways} to take
the uncertainty of the problem into account is to minimize the expected costs, which 
leads to a two-stage stochastic program.
This line of research was started
by~\citet{Lamiri2008Xie}; here, only overtime costs are considered, the elective surgery durations
are deterministic, and there is a single block for each day,
in which a random amount of additional emergency cases has to be operated.
This problem can be approximated using a sampled average approximation (SAA) approach,
which yields a Mixed-Integer Programming (MIP) formulation.
This work was continued in~\citet{lamiri2009optimization}, who analyzed the convergence
of the SAA  theoretically and proposed fast heuristics for this problem.

An approach closer to ours is that of~\citet{min2010scheduling},
who consider several surgical specialties and a given MSS. 
This paper takes the uncertain duration of the surgeries into account, as well the availability of downstream resources, with patients occupying recovery beds for a random length of stay.
In this paper as well, an SAA approach is used to solve the model. 

The article of~\citet{jebali2015stochastic} proposes a two-stage stochastic mixed-integer programming model with the goal of minimizing the costs incurred by overtime and under-utilization time, but no emergency operations are considered. The model takes into account the uncertainty related to the surgery duration and availability of downstream resources.
This approach was generalized in~\citet{jebali2017chance}, who consider emergency surgeries and penalty costs for exceeding the capacity of recovery areas. Similar to~\citet{Lamiri2008Xie},
the model assumes that there is a single block per time period (i.e., per day), so no decision has to
be made concerning the allocation of emergencies. The model is solved using an SAA-based approach.

All the aforementioned stochastic programming approaches rely on SAA. One of the main drawbacks of the sample average approximation is that it tends to be computationally intensive, as it yields a deterministic problem in which there is a copy of all
second-stage variables for each sampled scenario.
As a result, the considered MIPs grow quickly and can become intractable, even for rather small instances. 
To get around this problem, \citet{zhang2020column} proposed a column-generation-based approach to solve the large deterministic problem obtained by SAA. They use heuristics
to solve the pricing subproblems.

\medskip
\emph{Robust Optimization and chance-constrained problems.}
Unlike stochastic programming approaches, in which the expected value of the costs are minimized,
robust optimization aims at finding solutions in which the costs are small for all values of
the uncertain parameters lying in some uncertainty set. Obviously, such a worse-case approach
can be well adapted to the health care sector, where unexpected events can have dramatic
consequences. \citet{neyshabouri2017two} proposed a robust optimization approach 
for the advanced planning problem, in which the availability of downstream constraints is taken into
account. They formulate a two-stage robust optimization
problem with mixed-integer recourse, which is solved using a column-and-constraint
generation method.

A few works have also considered a setting in which the elective patients to be operated in the
next planning period have already been selected, and they must \emph{all} be allocated to one
of the available blocks of this period. In this case, the challenge is to find an allocation
minimizing the risk of overtime. 
In~\citet{Wang2014Tang}, a stochastic programming formulation is proposed for this purpose, with a chance constraint ensuring that the probability that the overtime exceeds some threshold remains small in order to avoid case cancellation. This problem is solved
by using the SAA approach and a column-generation based heuristic. In~\citet{sagnol2018robust},
a robust optimization approach is proposed to minimize the worse possible overtime cost, assuming
all processing times lie in some uncertainty set related to a lognormal distribution; 
a cutting-plane algorithm is proposed to solve this problem.

\medskip
\emph{Distributionally Robust Optimization.} A recent line of research offering 
a tradeoff between stochastic programming and robust optimization is called
\emph{distributionally} robust optimization; here, the expected cost function is taken with respect to an unknown distribution of the uncertain parameters, which is assumed to lie
in a so-called \emph{ambiguity set}.
This approach was used in the field of elective surgery planning by the authors of~\citet{shehadeh2021distributionally}, who
assume the surgery durations and lengths of stay in downstream units
have an unknown distribution, but the mean and support of these distributions is known.

When the ambiguity set is
a ball centered at the empirical distribution given by some historical data,
and the distance between distributions is a Wasserstein distance,
compact reformulations of the problem are sometimes possible; see e.g.~\citet{mohajerin2018data}.
This data-driven approach is used by~\citet{shehadeh2022data}, who
considers uncertain surgery durations and capacity requirement for emergency cases,
and gives a compact MIP reformulation.

\paragraph{Appointment scheduling within a single block.}
Given a list of patients to be operated within a block, a crucial step is to
define tentative starting times for their operations.
An operation is not allowed to start earlier than its tentative starting time,
as the patient might not be ready yet, and some members of the medical team
might not have arrived in the operation theater yet.
These tentative times incur two different types of cost. On the one hand, if
the operating room is available earlier, then idle time occurs. On the other hand, 
if the operation starts later than its scheduled starting time, waiting time occurs.

\citet{denton2003sequential} consider the problem of finding the 
optimal tentative starting times for a fixed sequence of the patients
by formulating the problem as a two-stage stochastic linear problem
and using the L-shaped method.

\citet{denton2007optimization} provide a two-stage stochastic mixed-integer program model to optimize simultaneously the sequence of operations and their tentative starting times,
with the goal of minimizing the total expected cost associated with waiting time, idle time, and overtime. The SAA is used to deal with the uncertainty in surgery duration, and it is shown
empirically that sequencing the operations in the order of nondecreasing variances is a very good
heuristic for this problem. This model was further enhanced in subsequent 
work~\citep{mancilla2012sample, Berg2014Denton, shehadeh2019analysis}, who consider
slightly different cost structures.

Recently, the near-optimality of the Shortest-Variance first rule (SVF)
was investigated theoretically by~\citet{de2021performance},
for an appointment scheduling problem with idle time and waiting time.
They give worse case bounds on the performance of the SVF rule compared to
an optimal solution, under different assumptions on the processing time distributions.
On the other hand, \citet{mansourifard2018heuristic} claim that
ordering the jobs according to their \emph{newsvendor-index} is better than SVF for
instances with unbalanced waiting and idle time. We also point out that while the SVF-rule is 
good to handle problems with idle, waiting and overtime, it could perform poorly
for other objectives. One such example has been considered by~\citet{van2012minimizing}, who suggest to order the surgeries in different blocks so as to minimize the length of the so-called maximal break-in interval (BII), i.e., the maximal time between two successive completion times on different blocks, which corresponds to the maximum amount of time an emergency patient in critical situation would have to wait to be operated, regardless of the time he arrives at the hospital. 

\paragraph{Online scheduling}

The computation of optimal online scheduling policies 
that make very short-term decisions to handle emergencies
has been very little investigated in the literature. The only reference we know of
in this area is~\citet{silva2020surgical}, who use an approximate dynamic programming
approach embedding an integer programming formulation to make decisions
about how to schedule arriving emergency patients and whether to 
cancel planed elective cases.

\section{Model and assumptions} \label{sec:model}

In this section, we formally define the \EESPlong (\EESP). We essentially follow the guidelines given in the description of the AIMMS-MOPTA competition~\citep{MOPTA22}, but explain some
aspects of the problem with more details, as
some modelling decisions were left to the
participants. A summary of symbols used in our model is given in Table~\ref{tab:symbols}.

\paragraph{Problem definition and notation.}
We focus on the \EESP in a hospital with flexible operating rooms shared with two classes of patients: elective patients and emergency patients. 
Following the standard block-scheduling paradigm  mentioned in the introduction,
each elective operation 
must take place within a block of OR-time of the MSS \emph{belonging} to its respective surgical specialty. However, the ORs are called \emph{flexible} because emergency patients can be inserted in the schedule of any block, independently {of} its specialty.

An instance consists of a pool of $n$ elective surgeries $I =\{1,2,\dotsc, n \}$, where  
each surgery belongs to one of several
surgical specialties $s\in S$ such as
cardiology, gastroenterology or orthopedics.
We denote the specialty of the surgery $i\in I$ by $s(i)\in S$, and we write $I_s$ for the elective surgeries of specialty $s$, i.e.,  $I_ s:=\{ i \in  I \mid s(i) = s \} .$

We denote by $P_i$ the duration of the elective surgery $i\in I$, which is a nonnegative random variable.
In our study we assume that $P_i$ has a lognormal distribution
$P_i \sim \mathcal{LN}(\mu_i,\,\sigma_i)$,
which is a reasonable assumption used in
different works, e.g.~\citet{strum1998surgical}
(note that lognormality of durations is not essential for our approach however).
In practice, the parameters $\mu_i$ and $\sigma_i$ of 
$P_i$ can be \emph{learned} from individual features of patient~$i\in I$ such as the type of surgical procedure, the patient's age and diagnosis, or an estimation of the duration done by the surgeon,
by training a model on the available historical data, see e.g.~\citet{sagnol2018robust}.

The MSS is given in the input and consists of a set of blocks~$B$ available during the planning period $D$,
which in this study corresponds to a week $D=\{\texttt{Mo},\texttt{Tu},\texttt{We},\texttt{Th},\texttt{Fr}\}$. 
Each block\linebreak $b\in B$  corresponds to some operating room of the hospital 
available for a regular working time of $T$~minutes 
at a given day $d\in D$.
Each block $b\in B$ is associated to a surgical specialty $s(b)\in S$, which means that 
an elective surgery $i\in I$ can only be assigned
to block $b$ if $i\in I_{s(b)}$, or, equivalently, $s(i)=s(b)$. We denote by $d(b)\in D$ the day of the block $b$.
In addition, we use a dummy block $b'$ which is used for the patients not operated during the current planning horizon.
We write $B'=B\cup\{b'\}$ and 
$B_s'=B_s \cup\{b'\}$, where $B_s$ stands for the set of blocks of specialty~$s$, i.e.,  $B_s := \{ b \in B \mid s(b)=s \}$.

A random set $E$ of emergency patients arrives during the planning period, and we denote by $E_d$ the subset of emergency patients arriving on day $d\in D$.
We assume that the set $E_d$ becomes known 
\emph{at the beginning of day d}, and $|E_d|$ follows a Poisson distribution with parameter $\lambda$.
Patients $e\in E_d$ must be operated on day $d$ and cannot be postponed, i.e., they must be allocated to a block $b\in B_d:=\{b\in B \mid d(b)=d\}$. As for the case of elective patients, 
we assume that a random distribution $P_e\sim \mathcal{L}\mathcal{N}(\mu_e,\sigma_e)$ is available for the duration of each patient $e\in E_d$, but this distribution becomes known only on day $d$ (because an individual prediction for patient~$e$ can only be inferred after observing the features of patient~$e$). Before day $d$, a scheduling policy is allowed to use some \emph{prior information} on the duration of emergency patients, though, which we assume is available in the form of the marginal mean $m_e$ and marginal variance $v_e$ across all emergency surgery durations. 

Finally, an instance of \EESP also depends
on the choice of different cost parameters. We denote by $\co$, $\ci$ and $\cw$ the cost of overtime, idle time and waiting time, respectively, and by $\cm$ the cost per migration. There is also a cost $c_{ib}$
for assigning patient $i$ to block $b$; The precise meaning of all these costs is described in the paragraph \emph{cost function} at the end of this section.

\paragraph{Output.}
As mentioned in the introduction, the \EESP can be decomposed in two phases. A solution to this problem is given by as solution to the \APP, which can be computed offline, as well as a policy (i.e., an online strategy) for the \OSP.

\medskip
\noindent\emph{\APPlong (\APP).} In this phase, we compute an assignment of each elective {patient} $ i \in I_s$ to a block $ b(i) \in B_s'$.
Patient $i\in I$ is called \emph{scheduled}
if $b(i)\in B$ and  \emph{postponed} if $b(i)=b'$. In addition, we must also decide a \emph{tentative starting time} $t(i)$ for each scheduled patient. 
Our approach for this phase is described in Section~\ref{sec:offlinephase}.

\begin{figure}[t]
 \begin{center}
\includegraphics[width=13cm]{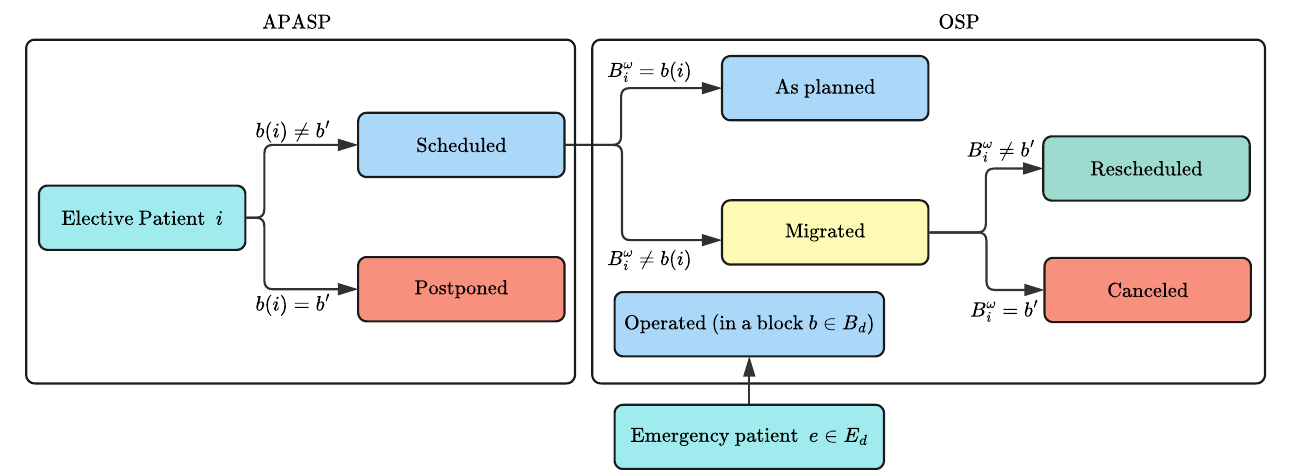}
 \end{center}
\caption{\small Classification of patients after solving the \APP (leftmost frame) and after executing the online policy in scenario $\omega$ (rightmost frame). \label{fig:flowchart}}
\end{figure}

\medskip
\noindent\emph{\OSPlong (\OSP).}
\changed{In this phase, the goal is to select a \emph{policy}}
which is able to take online decisions, such as starting the surgery of an emergency patient or reassigning the surgery of an elective patient
to another block. Reassigned surgeries are called
\emph{canceled} if their new block is $b'$ (which means that the surgery will take place in a later planning period), and \emph{rescheduled} otherwise; 
In the latter case, the online policy must also set a new tentative time.
We refer to~Figure~\ref{fig:flowchart} for a
schematic view summarizing the possible status of each patient
at the end of the online phase. 
\changed{The \OSP can in fact be modelled as} a
Markov Decision Process (MDP), as was done by~\citet{silva2020surgical}
for a similar problem.
\changed{While a formal, comprehensive definition of this MDP is beyond the scope of this article, we aim to provide a succinct overview of its key components
in the next paragraph.}
Our approach for the \OSP is
presented in Section~\ref{sec:onlinephase}.

\paragraph{MDP model of the the online phase.}

\changed{The main idea is to represent the available information at the $k$th decision time by a \emph{state vector} $\vec{s}_k$.
A state $\vec{s}$ is defined
by the current time $\tau_{\vec{s}}$ and day $d_{\vec{s}}$,
the current assignment $b_{\vec{s}}(i)$ and tentative starting time $t_{\vec{s}}(i)$
of all surgeries $i\in I$, the number $m_{\vec{s}}$ of migrations already performed,
the set of remainaing emergencies $E_{\vec{s}}$ for the current day,
the sets $F_{\vec{s}}$ and $O_{\vec{s}}$ of finished and ongoing operations,
and
the starting time of each $i\in O_{\vec{s}} \cup F_{\vec{s}}$.
The initial state $\vec{s}_0$ is defined by the output of the \APP (Phase 1),
with $\tau_{\vec{s}_0}=0$, 
$d_{\vec{s}_0}=\texttt{Mo}$,
$m_{\vec{s}_0}=0$, $F_{\vec{s}_0}=O_{\vec{s}_0}=\emptyset$,
$b_{\vec{s}_0}(i)=b(i)$ and $t_{\vec{s}_0}(i)=t(i)$, for all $i\in I$.
The set of remaining emergencies is intialized with
$E_{\vec{s}_0}=E_{\texttt{Mo}}$, the emergencies revealed on Monday,
and thus depends on} 
the scenario $\omega \in \Omega $ that represents a particular realization of the uncertainty (number of emergencies, posterior lognormal parameters $\mu_e,\sigma_e$ for each $e\in E$, and
surgery duration of all cases $i\in I\cup E$).

\changed{
A policy for the \OSP is a function $\pi$ that associates
each state $\vec{s}$ with a valid \emph{action} $a$.}
Feasible actions \changed{in state $\vec{s}$ at time $\tau$} of day $d$ are either of the following:
(1) starting a scheduled surgery $i$ in its current block \changed{$b_{\vec{s}}(i)$} if the block is available and \changed{$\tau \geq t_{\vec{s}}(i)$};
(2) starting an emergency surgery \changed{$e\in E_{\vec{s}}$} in any available block \changed{$b\in B_{d_{\vec{s}}}$};
(3)~deciding to migrate a planned surgery to a later day, by updating its current block and tentative starting time;
(4) or waiting for some fixed amount of time or until some surgery completes.
\changed{
The dynamics of our MDP are described by a  critical component known as the \emph{transition function},
which gives the next state $\vec{s}_{k+1}$, given
the action $a_k$ taken in state $\vec{s}_k$
and the scenario $\omega \in \Omega $.}
\changed{
 We avoid a lengthy, formal description of
$\vec{s}_{k+1}$ for each state-action pair $(\vec{s}_k,a_k)$ and
scenario $\omega\in\Omega$,
and opt for a more intuitive approach based on two illustrative examples,
the other cases being constructed in a similar fashion:
}
\vspace{-1em}

\changed{
\begin{itemize}
\item Assume that the action $a_k$ is to migrate patient $i$ to a new block $\tilde{b}$ with
tentative time $\tilde{t}$. Then, we simply set $b_{\vec{s}_{k+1}}(i)=\tilde{b}$,
$t_{\vec{s}_{k+1}}(i)=\tilde{t}$ and increment the migration counter: 
$m_{\vec{s}_{k+1}}=m_{\vec{s}_{k}}+1$.
\item As a second example, assume that the action $a_k$ is to wait until a surgery completes, so a block
becomes available. We denote by $p_i^{\omega}$ the realization of
the surgery duration $P_i$ in scenario $\omega$, and by
$S_i^\omega$ the realized starting time of a surgery $i\in O_{\vec{s}_k}\cup F_{\vec{s}_k}$.
Then, the next state $\vec{s}_{k+1}$ occurs at
time $\tau_{\vec{s}_{k+1}}:=\min_{i\in O_{\vec{s}_k}} \{S_i^\omega+p_i^\omega\}$.
The set $I^*$ of operations terminating at time
$\tau_{\vec{s}_{k+1}}$ are moved from \emph{ongoing} to \emph{finished},
i.e., we set $O_{\vec{s}_{k+1}}=O_{\vec{s}_{k}}\setminus I^*$ and $F_{\vec{s}_{k+1}}=F_{\vec{s}_{k}}\cup I^*$.
\end{itemize}
}

\changed{
The execution of the policy $\pi$ in scenario $\omega$ continues until
we reach a final state $\vec{s}_f=\vec{s}_f(\vec{s}_0,\pi,\omega)$ in which all
patients have been operated, postponed or canceled.
For a fixed solution of the \EESP, that is, an initial state $\vec{s}_0$
calculated after solving the \APP and a policy $\pi$  for the \OSP,
we denote by
$B_i^{\omega}:=b_{\vec{s}_f}(i)$ the block to which patient $i$ was finally assigned,
$T_i^\omega:=t_{\vec{s}_f}(i)$ its final tentative starting time
and by
$M^{\omega}=m_{\vec{s}_f}$ the total number of migrations.
Note that if a patient has been operated \emph{as planned},
i.e., not migrated during the execution of the schedule,
then we have $b_{\vec{s}_f}(i)=b_{\vec{s}_0}(i)$, that is, $B_i^{\omega} = b(i)$,
and simlilarly $T_i^{\omega} = t(i)$.}
The actual starting time of surgery~$i$ is denoted by $S_i^\omega$, and $L_b^{\omega}$ represents the load of block $b$, which is the latest completion time of a surgery scheduled in this block under scenario $\omega$, i.e., 
$L_b^{\omega}=\max_{\{i\in I\cup E_{d(b)} \mid \ B_i^\omega=b\}}\ S_i^\omega + p_i^\omega$.

\paragraph{Cost function.}

\changed{Our goal is to minimize the expected value
$\operatorname{Cost}(\vec{s}_0,\pi) := \mathbb{E}_\omega [\operatorname{cost}(\vec{s}_0,\pi,\omega)]$
of the total costs associated with the final state $\vec{s}_f$},
which are defined as a combination of
scheduling, waiting, idling, overtime and migration costs.
\changed{We use the notation $\operatorname{cost}(\vec{s}_0,\pi,\omega)$
to emphasize the dependence of the final state $s_f$ on $\vec{s}_0$, $\pi$ and $\omega$}:

\begin{equation} \label{costfunction}
\operatorname{cost}(\vec{s}_0,\pi,\omega)=
\sum_{i \in I} c_{i B^{\omega}_i} + 
\sum_{\substack{i \in I\\ B_i^\omega\neq b'}} \cw ( S^{\omega}_i - T^{\omega}_i)  
+ \sum_{b \in B} \ci( L^{\omega}_b  -\, \sum_{\mathclap{\substack{i\in I\cup E_{d(b)}\\ B^{\omega}_i = b}}} p^{\omega}_i) + \sum_{b \in B} \co  (L^{\omega}_b - T) +  \cm\, M^{\omega}
.
\end{equation}

The first term of the cost function adds
up the scheduling costs 
$c_{i B_i^{\omega}}$ of all elective patients $i\in I$. 
Recall that $c_{ib}$ represents the cost for allocating patient $i$ in block $b\in B'$.
In practice, theses costs can be tuned to achieve different goals, such as prioritizing
patients depending on the severity of their condition and their total waiting time,
minimizing hospitalization costs or 
preventing the usage of certain blocks for certain patients; see e.g.~\citet{lamiri2009optimization}.

\begin{table}[b!] 
	\centering
	 \caption { {\small List of symbols used to describe an instance of \EESP and its output}
	 \label{tab:symbols} }
	 \medskip
	 
{\small
\begin{tabular}{m{0.7cm}p{5.5cm}@{\qquad}m{1.2cm}p{5.9cm}}
\toprule
\multicolumn{4}{l}{\textbf{Input}}\\
\midrule
\vspace{0.5em} 
$B$   & set of available surgery blocks $b$ & $s(b)$ & specialty of block $b$ \\[0.5em]
$I$   & set of elective surgeries  & $s(i)$   & specialty of surgery $i\in I$\\[0.5em]
$T$ & regular working time of block $b$
& $\mu_i, \sigma_i$ &  lognormal parameters of elective surgery~$i\in I$\\[0.5em]
$d(b)$ & day of block $b$ & $m_e, v_e$ & marginal mean and variance of an emergency surgery \\[0.5em]
$\lambda$ & emergency rate\\[0.5em]
\midrule
 \multicolumn{4}{l}{\textbf{Costs}}\\
\midrule
\vspace{0.5em} 
$c_{ib}$ & cost to assign patient $i$ to block $b\in B'$   &  $\co$   &  cost parameter for overtime  \\[0.5em]
 $\ci$   &  cost parameter for the idling time & $\cw$   &   cost parameter for the waiting time    \\[0.5em] 
  $\cm$   & cost parameter for migrations & &  \\[0.5em]
\midrule
\multicolumn{4}{l}{\textbf{Output (\APP)}}\\
\midrule
 $b(i)$ & block where patient $i\in I$ is assigned at the end of the Advanced Planning & $t(i)$  &tentative start time of a scheduled surgery $i\in I$\\[0.5em]
 \midrule
\multicolumn{4}{l}{\textbf{Additional Input revealed online}} \\
\midrule
$E_d$ & set of emergency surgeries to be operated on day $d$, revealed at the beginning of day~$d\in D$, with $|E_d|\sim \mathcal{P}(\lambda)$  & $\mu_e, \sigma_e$ &  lognormal parameters of emergency surgery $e$. Revealed at the beginning of day~$d$ for $e\in E_d$\\
\midrule
\multicolumn{4}{l}{\textbf{Output of the \OSP for a given scenario $\omega$ and a policy $\pi$}}\\
\midrule
$B^{\omega}_i$ & block where $i$ was (finally) assigned to, with $B^{\omega}_i=b(i)$ if i was not migrated & $S_i^\omega$& starting time of surgery $i$ \\[0.5em]
$L_b^\omega$& load of block b & $T^{\omega}_i$ & tentative start time of surgery $i\in I$, with $T^{\omega}_i=t(i)$ if $i$ has not been rescheduled
  \\[0.5em]
 $M^{\omega}_i$& number of times that $i$ was migrated\\
\midrule
\multicolumn{4}{l}{\textbf{Other symbols used for modelling}}\\
\midrule
\vspace{0.5em}
 $b'$ & Dummy block &  $I_s$   & =$\{ i \in I \mid s(i)=s \}$; Elective patients of specialty $s$  \\[0.5em]
$B'$   & =$  B \cup \{ b' \}$
& $n$ & =$|I|$; number of elective surgeries \\[0.5em]
$B_s$   & =$ \{ b \in B \mid s(b)=s \}$; blocks of specialty~$s$ & $\mathcal{LN}(\mu,\,\sigma)$ & lognormal distribution with parameters $\mu, \sigma$ \\[0.5em]
$B'_s$   & =$B_s \cup \{ b' \}$  & $\mathcal{P(\lambda)}$ & Poisson distribution with parameter $\lambda$  \\[0.5em]
$B_d $   &=$\{ b \in B \mid d(b)=d \}$; blocks of day $d$ & $P_i$ & random duration of surgery \mbox{$i\in I\cup E$}; $P_i\sim\mathcal{LN}(\mu_i,\sigma_i)$ \\[0.5em] 
\vspace{0.5em}
$D$ & $=\{d(b) \mid b\in B\}$; Set of days in planning horizon & $p_i^\omega$ & Realization of $P_i$ in scenario $\omega$ \\ \vspace{0.5em}
$E$ & $=\cup_{d \in D} E_d$; Set of emergency surgeries & $S$   &  $=\{s(i) \mid i\in I\}$; Set of all specialties\\
\bottomrule
\end{tabular}
}
\end{table} 

\changed{
The next part of the objective function counts the waiting time, weighted with a factor $\cw$. The waiting time for an elective surgery $i\in I$ is the difference between its actual starting time and its tentative starting time. Note that this difference is always non-negative {as a surgery cannot start before its
scheduled (tentative) starting time.} Then, we consider a cost penalizing the idle time, i.e., the amount of time between $0$ and $L_b^\omega$ where no surgery takes
place in block $b$, weighted by a factor $\ci$. There is also a cost of $\co$ per unit of overtime, that is, the time
between the end $T$ of the regular working time
and the completion of the last operation in a block.
Finally, the last term in~\eqref{costfunction} accounts for the total number of migrations, weighted with a factor $c^{\mathbf{m}}$.
}

\begin{figure}[t]
\begin{center}
\includegraphics[width=0.8\textwidth]{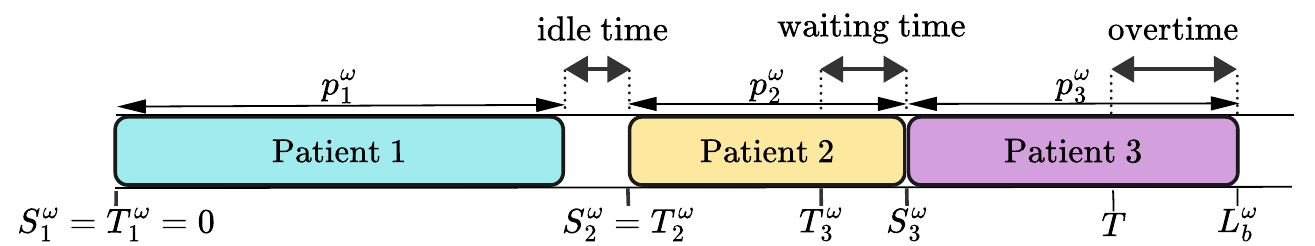}  
\end{center}
\caption{\small Schedule in a block for a realization $\omega$ of the uncertainty. In this realization, the tentative starting time $T_2^\omega$ is too large, which causes idle time in the block, while $T_3^\omega$ is too small, which results in waiting time. The last operation in the block completes after time $T$, which generates some overtime. Note that the tentative starting time $T_i^\omega$ usually does not depend on the scenario $\omega$ and is set to the value of $t(i)$ computed in the \APP, unless a migration of patient $i$ occurs in the scenario $\omega$ during the \OSP.
\label{fig:surgery_block}}
\end{figure}

An illustration of the concepts of waiting time, idle time and overtime in a block is shown in Figure~\ref{fig:surgery_block}.
Note that scheduling costs, waiting costs, and migration costs are counted for elective patients only. There are no explicit costs for the emergency patients, who arrive randomly and
have to undergo urgent surgery on the same day. However, the insertion of these patients into
one of the available blocks incur some indirect costs, as they may delay other
surgeries, create overtime or impose some case cancellations, for example.

\noindent \emph{Remark:} Unlike other works, our formulation does not take into account costs for opening an operating room or underutilization costs, but it would be straightforward to modify our approach to handle these costs, \changed{see~\cite{sagnol2018robust} and~\cite{denton2003sequential} for
articles handling opening and underutilization costs, respectively}.

\begin{figure}[t]
	\begin{center}
	\includegraphics[width=\textwidth]{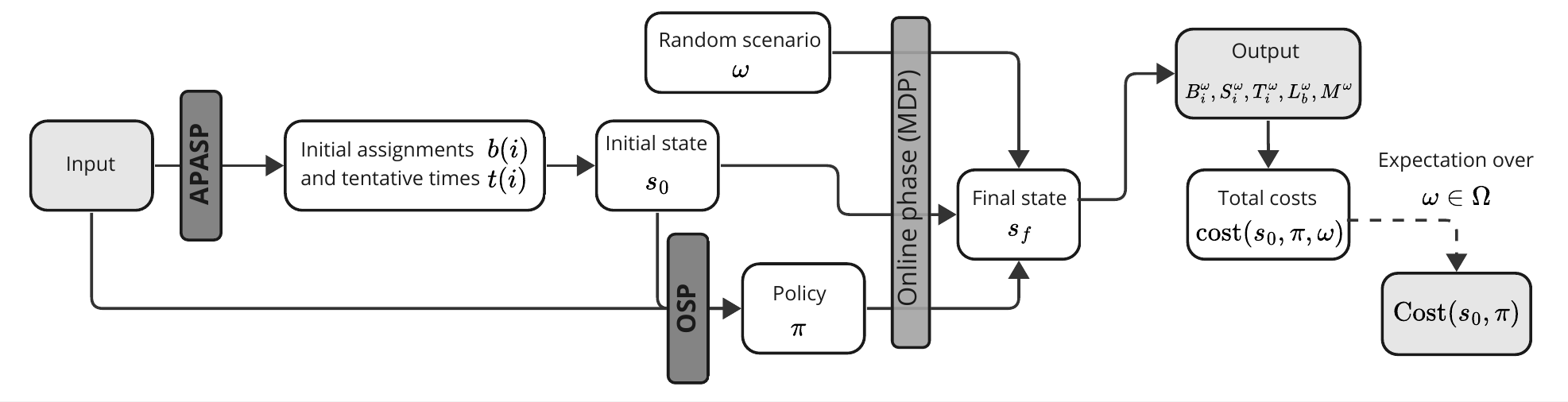}  
	\end{center}
	\caption{\small Summary of interactions between the
	different components of the \EESP
	\label{fig:APASP-OSP}}
\end{figure}

\paragraph{Relationship between different models.}
\changed{
In summary, the \EESP asks to minimize
$\operatorname{Cost}(\vec{s}_0,\pi)$
jointly over the initial state $\vec{s}_0$ and the online policy $\pi$.
The online phase of the problem corresponds to the \OSP,
which aims at finding a policy $\pi^*$ minimizing $\operatorname{Cost}(\vec{s}_0,\pi)$ for a given initial state $\vec{s}_0$.
The primary focus of this article in on the offline phase of
the problem (\APP), in which the goal is to compute a state $\vec{s}_0$
corresponding to initial patient-to-block assignments and
tentative times. The interactions between the different models
are summarized in Figure~\ref{fig:APASP-OSP}.
}

\section{Surrogate model-based two-stage stochastic programming for the \APP} \label{sec:offlinephase}

\changed{
Solving the \OSP to provable optimality is a very challenging task,
because the solution methods for MDPs suffer from
what is called the \emph{curse of dimensionality}; see, e.g., \citet{bertsekas2012dynamic}.
As a result, the \APP is computationally intractable, too, since we would need
to solve the \OSP to evaluate the quality of an assignment.
Recognizing this computational infeasibility, we adopt an optimization
strategy for the \APP that selects the initial state $\vec{s}_0$ under the assumption
that a simplified policy is employed within the MDP framework.
Specifically, we compute the patient-to-block assignment and tentative starting times
by assuming that no migration decision is ever made during the online phase,
and assuming that a simple rule to insert emergencies at the end of each day is used.
The greedy online policy which is presented in~\ref{sec:onlinephase}
and which we use in our computational study
is designed in a way that decisions differing
from these simplified rules
only occur if they are likely to improve the costs, so that our approach for the
\APP can be interpreted as the minimization of an upper bound for
the effectively achieved costs, providing a pragmatic and computationally tractable
solution for this problem. 
}

\changed{
To solve the \APP, we rely on a simplified policy
in which no patient migration is
ever carried out, and emergencies are assigned to the end of predefined blocks
depending solely on their order of arrival.
This means that we reserve capacity in some block
for the first emergency of the day,
for the second one, for the third one, and so on.}
Since in our model the number of
emergency arrivals can be arbitrarily large ($|E_d|$ is Poisson-distributed),
in practice we use a parameter $n_e$ for the maximum number of emergencies
we consider, 
\changed{
and we truncate the Poisson distribution at $n_e$. That is,
we estimate the probability distribution of $|E_d|$ on any day $d\in D$ by
$\mathbb{P}(|E_d|= j)\simeq \pi_j:= \frac{\lambda^j/j!}{\sum_{k=0}^{n_e} \lambda^k/k!}$,
for all $j=0,\ldots,n_e$.
}

\changed{
Under these assumptions, 
we can model the problem to be solved as the following two-stage stochastic program:}
\begin{align}
\min_{x,\changed{z}, t\geq 0} &\quad \sum_{i\in I} \sum_{b\in B'_{s(i)}} x_{ib}\, c_{ib}  + \sum_{b\in B} \mathbb{E}_\omega [\; \varphi_{b}(x, \changed{z}, t, \omega)], \label{P0}\\
\text{s.t.} &\quad
\sum_{b\in B_{s(i)}'} x_{ib} =1  && \forall \; i \in I,\nonumber \\
&\quad \sum\limits_{b \in B_d}   z_{jb}  = 1, && \forall d \in D , \, \forall j \in \{1,\ldots,n_e\},\nonumber \\
&\quad x_{ib} \in \{ 0, 1\}, && 
\forall \; i \in I, \; \forall b\in B_{s(i)}' \nonumber\\
&\quad z_{jb} \in \{ 0, 1\} && 
\forall \; j \in \{1,\ldots,n_e\}, \; \forall d\in D,\; \forall b\in B_{d}, \nonumber
\end{align}
where
$x_{ib}=1$ indicates that the elective patient $i\in I$ is allocated to block $b\in B_{s(i)}\cup \{b'\}$,
\changed{$z_{jb}=1$ indicates that the $j$th emergency
of day $d(b)$ is assigned to block $b\in B_d$}
and
$\varphi_b(x,\changed{z},t,\omega)$ 
represents the idle, waiting and overtime costs in block $b$
for the assignment $(x,z)$, tentative starting times $t$ and scenario $\omega$.
\changed{
The constraints state that each elective patient
must be allocated to a block of his specialty
or to the dummy block $b'$ (i.e., the patient is postponed),
and that for each day, a block has been planned for
each of the first $n_e$ arrivals.}

\changed{
We call $\varphi_b(x,z,t,\omega)$ the second-stage costs,
because they depend on the realization $\omega$,
and we will see that these costs can be computed by
solving an optimization problem
involving \emph{second-stage variables}, such as the
actual starting times of each operation in scenario $\omega$.
In contrast, $x$ and $t$ are called \emph{first-stage variables} because they must be set before observing
$\omega$.
}

The above problem can be rewritten as
\begin{equation} \label{P1}
\min_{\changed{(x,z)}\in X}\quad \sum_{i\in I} \sum_{b\in B'_{s(i)}} x_{ib}\, c_{ib}  + \sum_{b \in B} \phi_{b}(x,\changed{z})
,
\end{equation}
where $X$ denotes the feasible domain of the $x$ \changed{and $z$}-variables in Problem~\eqref{P0} and 
\begin{equation}\label{P2}
\phi_b(x\changed{,z}):=
\min_{ t} \; \mathbb{E} [\varphi_{b}(x, \changed{z}, t, \omega)]
\end{equation}
is the solution of the inner minimization problem. We call Problem~\eqref{P2} the 
\emph{appointment problem}.
\begin{table}[t!] 
	\centering
	 \caption {\small List of symbols used for the \APP \label{tab:symbols_2s} }
	 \medskip
{\small
\begin{tabular}{m{0.5cm}p{5.7cm}@{\qquad}m{0.9cm}p{5.9cm}}
\toprule
\multicolumn{4}{l}{\textbf{\changed{Surrogate model (MIP~\eqref{2MIPs}) for
the outer minimization problem~\eqref{P1}}
}}\\
\midrule
\vspace{0.5em} 
$X$ & Domain of all feasible assignments &  $\phi_{b}(x,z)$ &Optimal second-stage costs
for the assignment $(x,z)\in X$ \\[0.5em]
$x_{ib}$ & Binary variable to assign patient $i$ to block~$b$ & $z_{jb}$ & Binary variable for the allocation of the $j$-th emergency of day $d(b)$ in block $b$ \\[0.5em]
$y_{jb}$ & Approximation of optimal second-stage
costs in block $b$ if there are $j$ emergency patients
on day $d(b)$ & \(\bar{L}_{jb}\) & Expected surgery duration assigned to block $b$ if there are $j$ emergency patients on day $d(b)$ \\
\midrule
\multicolumn{4}{l}{\textbf{\changed{LP model~\eqref{LPstage2} of the appointment problem~\eqref{P2}}}}\\
\midrule
$T$ & Regular working time of block $b$ &$K$& Number of sample scenarios ($\omega_1,\ldots,\omega_K$)\\[0.5em]
$s_{i,k}$   & Starting time of surgery $i$ in scenario $\omega_k$   & $t_i$   & Tentative start time of surgery $i$\\[0.5em]
$L_k$ & Load of the block in scenario $\omega_k$& $O_k$ & Overtime in the block in scenario $\omega_k$ \\[0.5em]
\bottomrule
\end{tabular}
}
\end{table} 

In what follows, we explain our approach to solve the above two-stage stochastic program. We first show how to solve the appointment problem
by linear programming (LP), and then we will see that the optimal second-stage costs~\eqref{P2}
can be well approximated by a convex piecewise linear function. Finally, we use this approximation to reformulate~\eqref{P1} as a mixed integer program (MIP).

\paragraph{Solving the \changed{appointment} problem.} \label{sec:ssp}
We solve the appointment problem~\eqref{P2}
by using a sample average approximation, i.e., we sample $K$ scenarios $\omega_1,\ldots,\omega_K$ from the universe $\Omega$, and
replace the expectation $\mathbb{E}_\omega[\varphi_b(x,\changed{z},t,\omega)]$
with the approximation $\frac{1}{K}\sum_{k=1}^K
\varphi_b(x,\changed{z},t,\omega_k)$.
We refer to Section~\ref{sec:study} for a discussion
\changed{on the choice of an appropriate value of $K$}.

For a fixed block $b$ and an assignment $x\in X$, denote by $I_b(x):=\{i\in I \mid x_{ib}=1\}$ the set of patients assigned to block~$b$.
To ease the notation we reindex the patients
of $I_b(x)$ as $1,2,\ldots,n_b$.
\changed{
In theory, it would be necessary to optimize the order in which patients
are operated in the block in order to solve~\eqref{P2}.
Given that the number $n_b$ of patients in any block is small
(typically between 2 and 7), and that the approach presented below takes only 
a fraction of a second to solve the problem for a fixed ordering,
we note that brute-force enumeration of the $n_b!$ possible permutations of $I_b(x)$
would be a viable solution.
For the computational study in Section~\ref{sec:results}, 
we nevertheless preferred to use the 
shortest variance first (SVF) sequencing rule,
i.e., we only solve the appointment problem for the ordering of $I_b(x)$
such that $\operatorname{Var} (P_1) \leq \cdots \leq \operatorname{Var}  (P_{n_b})$.
Indeed, it is known that SVF is}
a very good heuristic for appointment scheduling
---see the discussion on this subject in Section~\ref{sub:relwork}---
\changed{and we have verified numerically that the SVF rule is near-optimal for the
appointment problem in a dedicated paragraph of Section~\ref{sec:results}; see Table~\ref{tab:svf-opt}.
}

We will now reformulate~\eqref{P2} as a 
linear program.
A similar approach was proposed in~\citet{denton2003sequential}, following the work
of~\citet{strum1999surgical}. Let us first introduce some variables.
For a given scenario $\omega_k$, $k=1,\ldots,K$, and 
\changed{an elective} patient $i$, we denote by $s_{i,k} \geq 0$ the starting time of the operation of patient $i$.
The variable $t_i \geq 0$ is the tentative starting time of patient $i$,
\changed{and does not depend on $\omega$ because it
is a first-stage variable.}
\changed{
We denote by $e_k$ the total duration
of emergency patients assigned to the end of block $b$
in scenario
$\omega_k$, that is, 
$e^{\omega_k} = \sum_{\{e\in E_{d(b)}: z_{j_e b}=1\}} p_e^{\omega_k}$,
where $j_e\in\{1,\ldots,n_e\}$ is the index of $e$ in $E_{d(b)}$.
}
The load $L_k$ of the block in scenario $\omega_k$ is defined
as the completion time of the last patient and is thus equal to 
$L_k = s_{n_b,k} + p^{\omega_k}_{n_b} \changed{+e^{\omega_k}} $.
We also introduce a variable $O_k \geq 0 $ for the overtime in the block in scenario $\omega_k$, and
which must satisfy the constraint $O_k \geq L_k - T $.
The overall cost incurred by waiting time, idle time and overtime in scenario $\omega_k$  
is thus
$
\sum_{i=1}^{n_b} \cw ( s_{i,k} - t_i)  + \ci ( L_k -\sum_{i=1}^{n_b} p^{\omega_k}_i - \changed{e^{\omega_k}}) + \co O_k 
$.
Every surgery has to start after its tentative starting time, i.e., for all $k\in\{1,\ldots,K\}$ and for all $i$ we enforce $s_{i,k} \geq t_i$. Further, a surgery can only start if the previous  one is completed, i.e.,  $s_{i,k} \geq s_{i-1,k} + p^{\omega_k}_{i-1}$, $i=2, \ldots, n_b$.
We thus arrive at the following LP-formulation for the appointment problem in block~$b$:

\begin{align}
\phi_b(x, \changed{z}) = \min&\quad \frac{1}{ K} \label{LPstage2}
\rlap{$\displaystyle{\sum_{k=1}^K} \Big( \sum_{i=1}^{n_b} {\cw} ( s_{i,k} - t_i)  + {\ci} ( L_k -\sum_{i=1}^{n_b} p^{\omega_k}_i \changed{- e^{\omega_k}}) + {\co} O_k \Big)$,}\\
&\quad s_{i,k} \geq t_i,  && i=1, \dotsc,  n_b,\quad k=1,\ldots,K \nonumber\\
&\quad s_{i,k} \geq s_{i-1,k} + p^{\omega_{k}}_{i-1}, &&i=2, \dotsc, n_b,\quad k=1,\ldots,K \nonumber\\
& \quad L_k = s_{n_b,k} + p^{\omega_{k}}_{n_b} + \changed{e^{\omega_k}}, &&k=1,\ldots,K\nonumber\\
&\quad O_k \geq L_k - T, &&k=1,\ldots,K\nonumber\\
&\quad s_{i,k} \geq 0, &&i=1, \ldots, n_b,\quad k=1,\ldots,K \nonumber\\
&\quad t_i \geq 0, &&i=1, \ldots, n_b \nonumber\\
&\quad O_k \geq 0, &&k=1,\ldots,K.\nonumber
\end{align}

\paragraph{Surrogate model of second-stage costs.} \label{sec:approx}

\changed{
Assume (temporarily) that no emergency is ever allocated
to block $b$, i.e., $z_{jb}=0$, for all $j=1,\ldots,n_e$.
The optimal second-stage costs of block $b$ are thus
equal to $\phi_b(x,0)$.}

We observed numerically that $\phi_b(x\changed{,0})$ can be  very well approximated by a convex function depending solely on the specialty $s=s(b)$ and 
the total expected processing time assigned to block $b$, i.e., there exist convex functions $f_s:\mathbb{R}\mapsto\mathbb{R}$ such that
$\phi_b(x,\changed{0}) \approx f_{s(b)}(\sum_{i\in I_b} x_{ib} \mathbb{E}[P_i])$. 
To compute these approximations, in a preprocessing step we solve the LP~\eqref{LPstage2} $N=1000$ times  for each surgical specialty $s\in S$, for a different choice of $n_b$ and different random variables $P_1,\ldots,P_{n_b}$ generated for the specialty $s$,
following the procedure described in Section~\ref{sec:assumptions}.
From this process, we get a set of $N$ points 
with coordinates 
$x_k:= \sum_{i=1}^{n_b} \mathbb{E}[P_i]$
and $y_k := \phi_b^*$, the optimal value of the LP.
Note that the computing time of this pre-processing step
\changed{does not really matter since it is carried out a single time},
so $K$ can be set to a large value; in our experiments we used $K=1000$.

Then, in order to compute the approximation $f_s$,
we split the points into $R  = 3 $ parts\footnote{We also tried splitting the data into more parts ($R=5$ or $R=7$), but the gain in accuracy was very small.} based on their $x$-values,
and we use a linear regression in every part, yielding a linear function $f_{sj} (x) = \alpha_{sj} x + \beta_{sj}$ for the part $j$.
In the end, we define the function $f_s: \mathbb{R}_{\geq 0} \mapsto \mathbb{R}$, with
\[
f_s(x) = \max_{j  = 1, \dotsc, R}\ f_{sj} (x)  = \max_{j  = 1, \dotsc, R }\  \alpha_{sj}\, x + \beta_{sj},
\]
as the elementwise maximum of the linear functions $f_{sj}, j = 1, \dotsc R$. Therefore $f$ is convex by construction. 

\begin{figure}[t]
\centering
\begin{tabular}{cc}
\includegraphics[width=6cm]{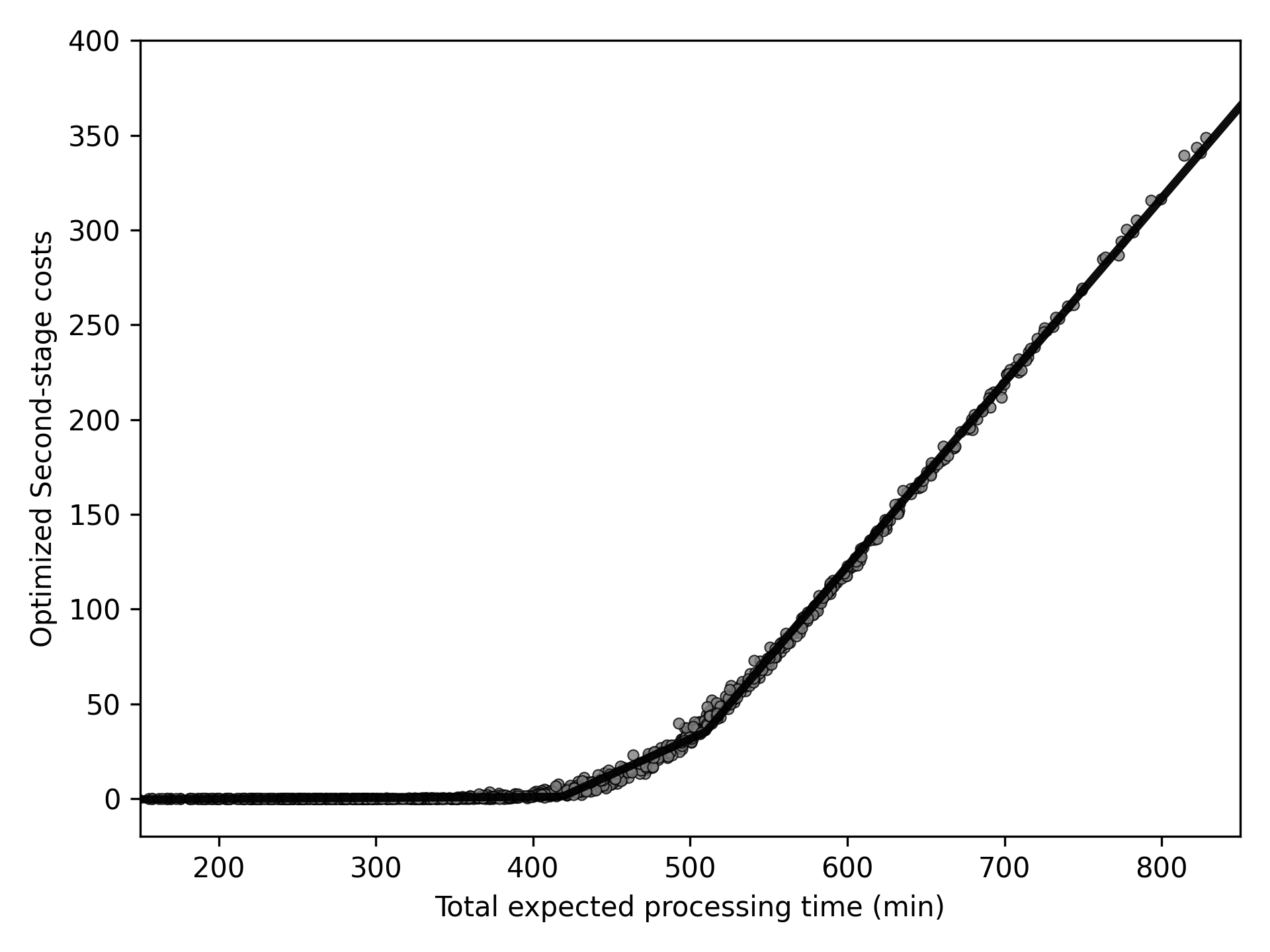}&
\includegraphics[width=6cm]{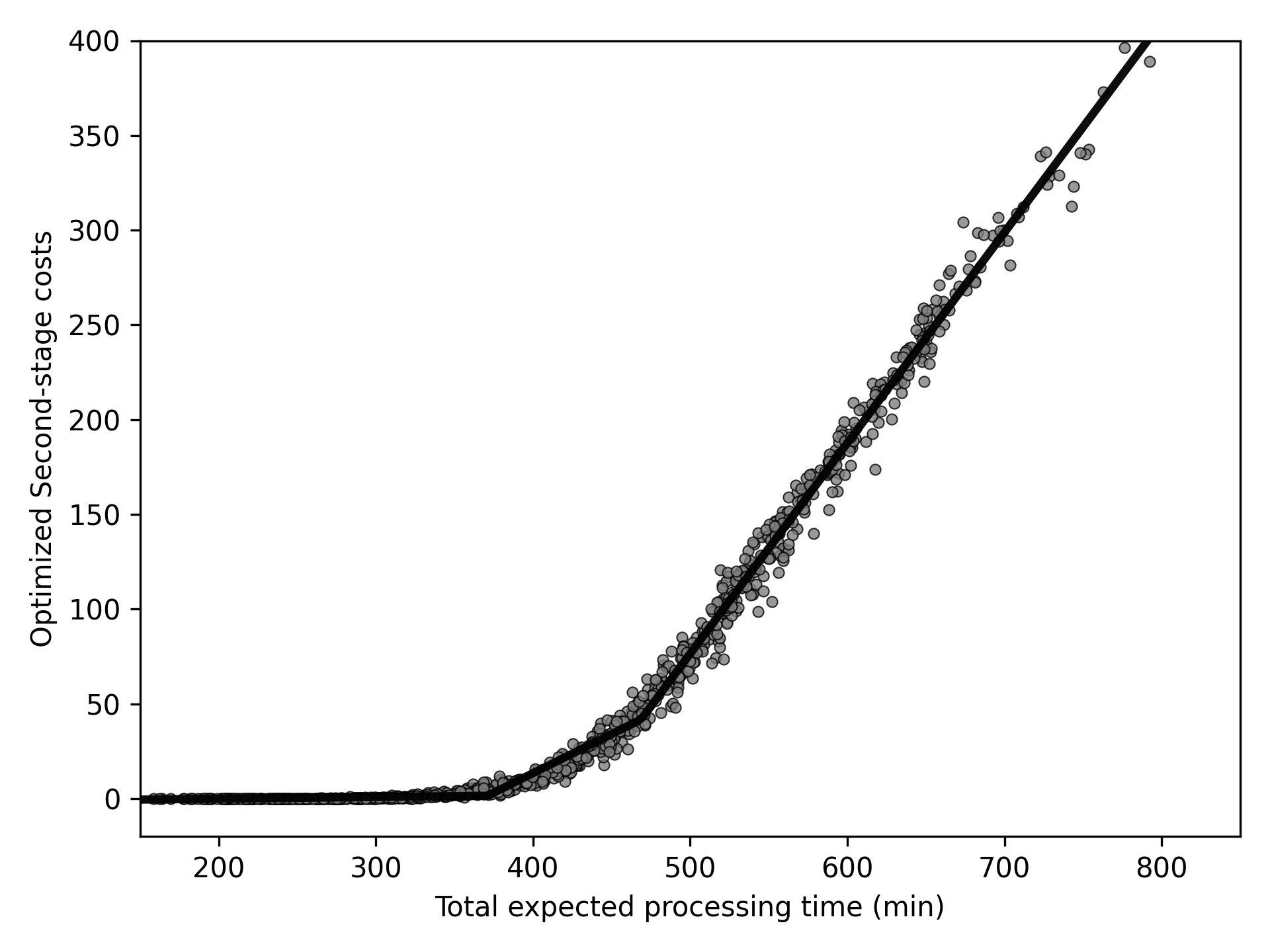}\\
{\small Cost structure $cs_1$ ($\co=1,\cw=0,\ci=0$)}
&
{\small Cost structure $cs_2$ ($\co=1,\cw=1,\ci=0$)}\\
\includegraphics[width=6cm]{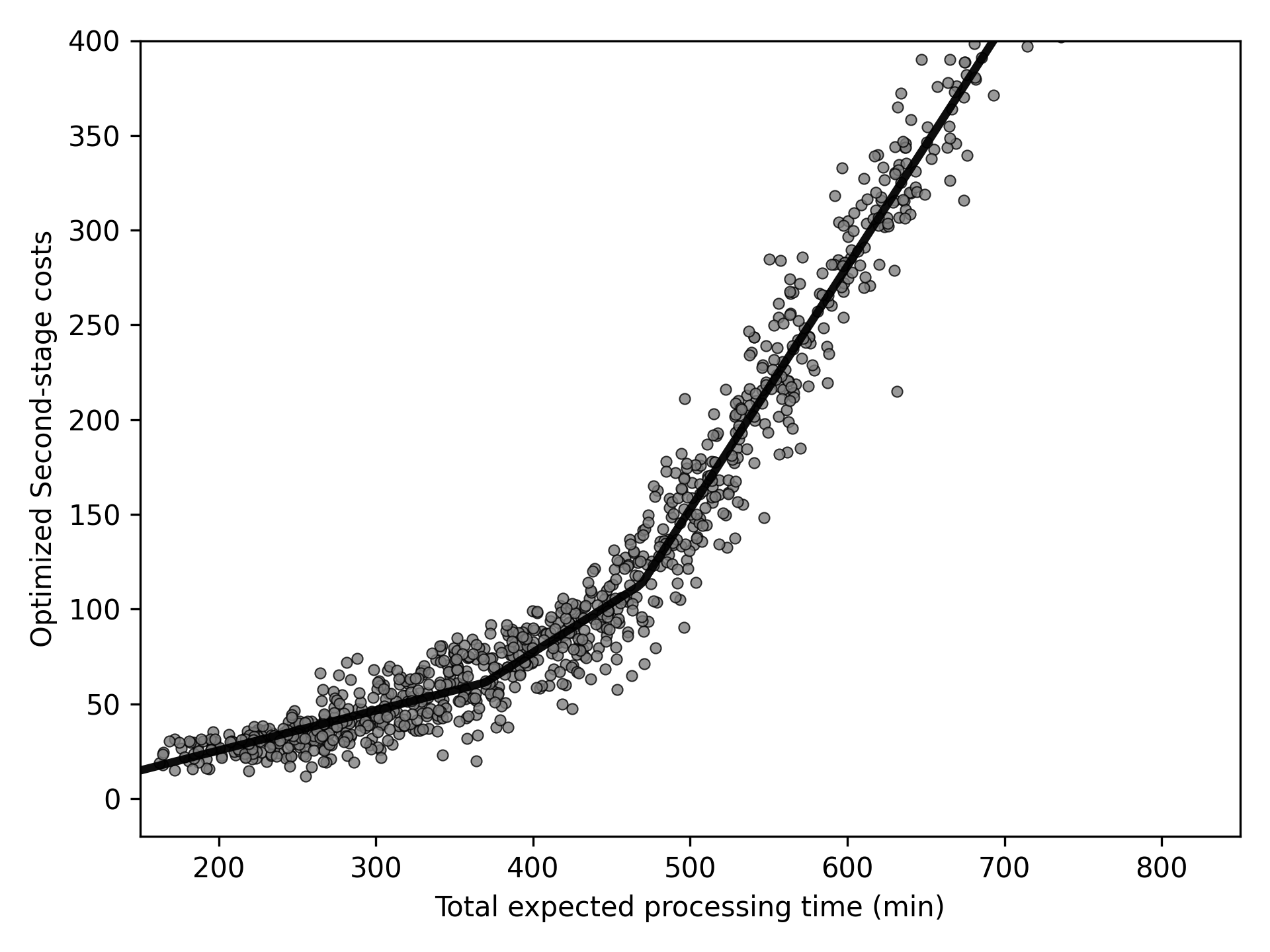}&
\includegraphics[width=6cm]{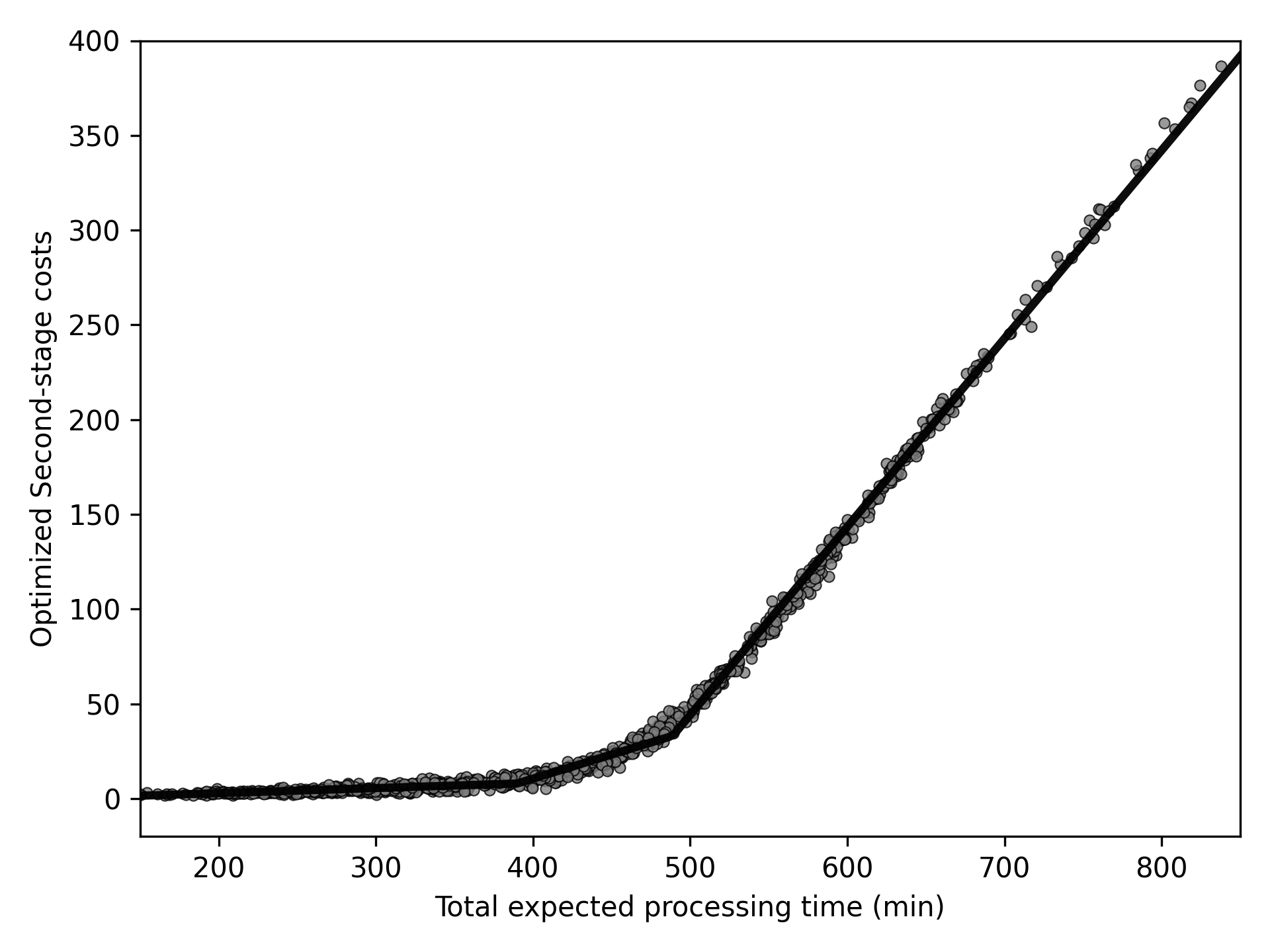}\\
{\small Cost structure $cs_3$ ($\co=1,\cw=2,\ci=2$)}
&
{\small Cost structure $cs_4$ ($\co=1,\cw=\frac{2}{15},\ci=\frac{2}{3}$)}\\

\includegraphics[width=6cm]{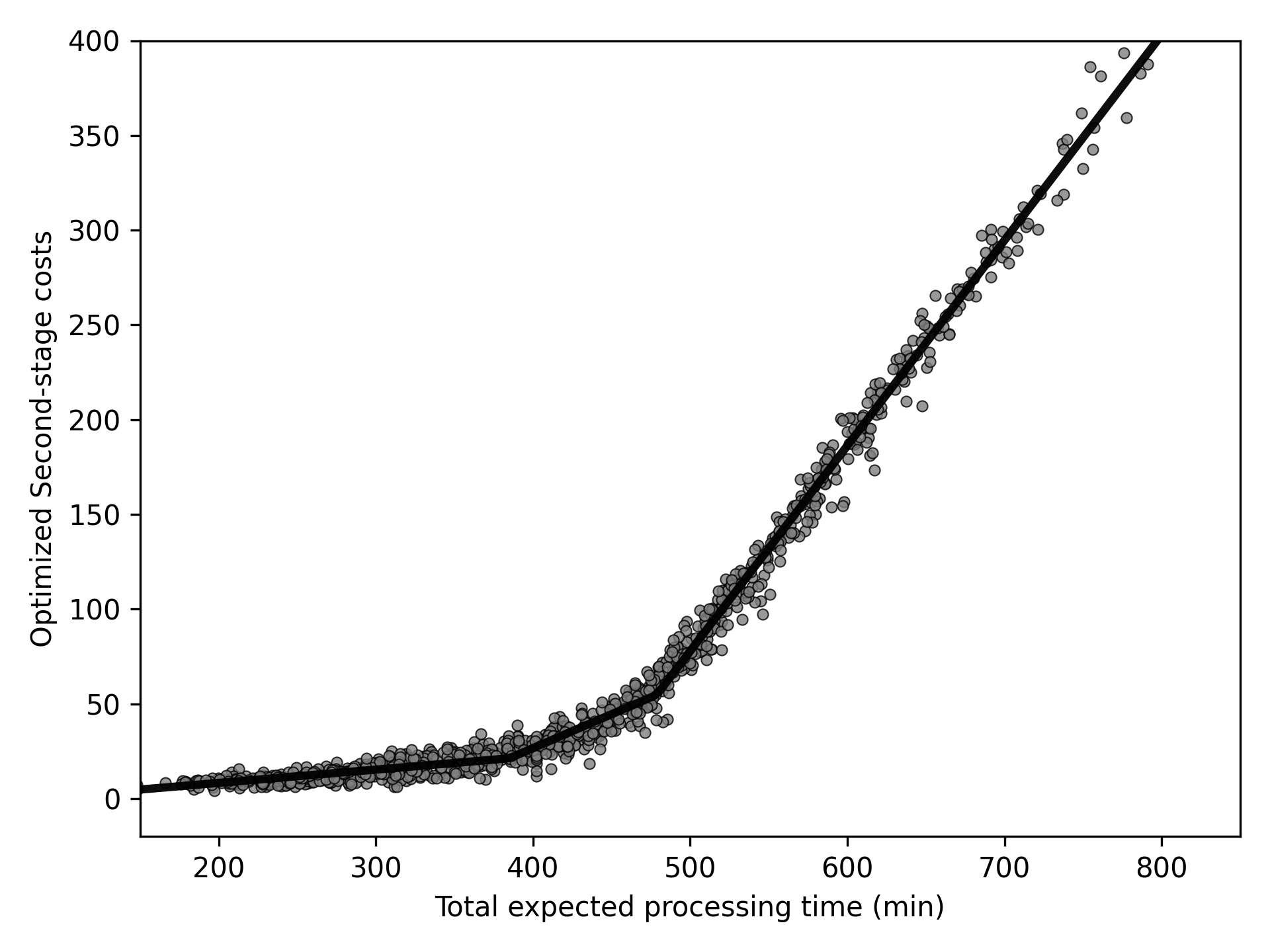}&
\includegraphics[width=6cm]{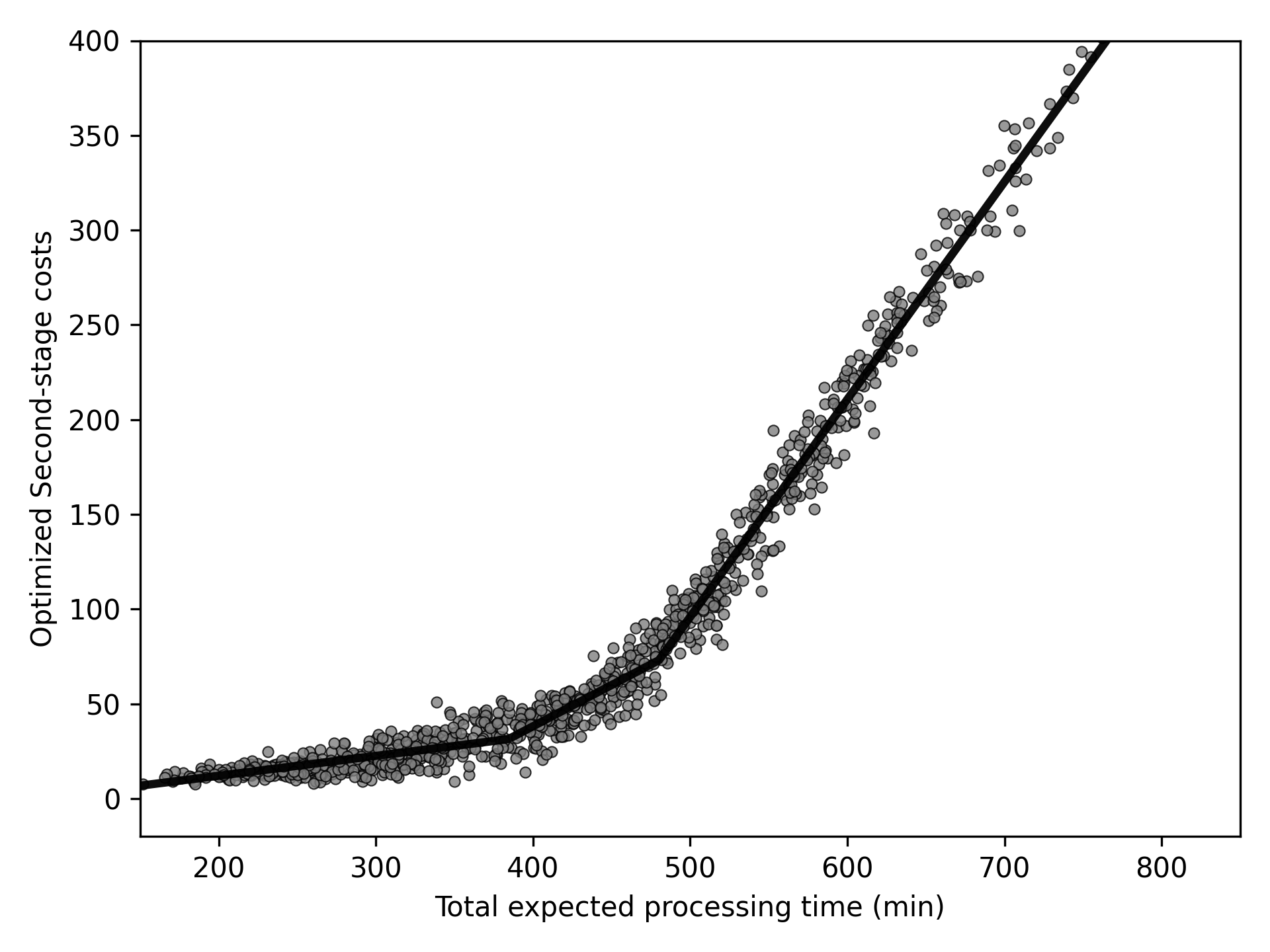}\\
{\small Cost structure $cs_5$ ($\co=1,\cw= \frac{2}{3},\ci=\frac{2}{3}$)}
&
{\small Cost structure $cs_6$ ($\co=1,\cw=1,\ci=1$)}
\end{tabular}
	
\caption{\small Cost of $N=1000$ instances of the second-stage problem
for a block $b$ of specialty $s=\texttt{ORTH}$,
as a function of the total expected duration.
The solid line shows the piecewise-linear approximation $f_{\texttt{ORTH}}$ with $R=3$ pieces,
 for $6$ different cost structures.
\label{fig:approx}}
\end{figure}

The cloud of $N=1000$ points $(x_k,y_k)$ obtained for the specialty $s=\texttt{ORTH}$ is shown on Figure~\ref{fig:approx}, together with its piecewise linear approximation, for the six
different cost structures used in our computational study; cf.\ Section~\ref{sec:assumptions}. 
\changed{For all cost structures and all specialties $s$, the mean deviation
$\frac{1}{N} \sum_{k=1}^{N} |f_s(x_k)-y_k|$
represents between 6\% and 12\% of
the second-stage costs for a regular load of 480 min.
We note that the individual relative deviation $|f_s(x_k)-y_k|/y_k$
is much larger for some data points, 
especially when 
$y_k$ is close to $0$, 
but this does not detract from the validity
of our approach: the key point is only to capture 
the general trend of increasing second-stage costs with the load of a block.
}

A remark on the convexity of $f_s$ is in order: it is well-known that the
optimal value of a minimization-LP like~\eqref{LPstage2} is a convex function of its vector
of right-hand sides; cf.\ e.g.~\citet{boyd2004convex}.
Therefore, we know that $\phi_b(x,0)$ is a convex function of the durations $p_i^\omega$, $i=1,\ldots,n_b$, $\omega\in\Omega$. The novelty of our approach lies on the fact that this 
convex function can be well approximated by a \emph{one-dimensional} function applied to the total expected processing time allocated to the block. A parallel can thus be drawn with the Benders decomposition approach, in which a convex underestimator of $x\mapsto \phi_b(x,0)$ is constructed and sequentially improved; see Section~\ref{sec:other_approaches}. In our case, we only \emph{estimate} the second-stage costs, but the restriction to a one-dimensional estimator $f_s$ significantly reduces the computational burden; see Section~\ref{sec:results} for a comparison.

\changed{
Now, let us return to the general case \emph{with} emergency arrivals.
For an emergency assignment~$z$,
we denote by $k(j):=\sum_{j'=1}^j z_{j'b}$ the number
of emergencies operated in block $b$
when there is a total of $j$ emergencies
on day $d(b)$, i.e., 
$|E_{d(b)}|=j$.
}
\changed{Following the law of total expectation,
we propose to estimate the optimal second stage costs
by
\[
\phi_b(x,z) \approx \sum_{j=0}^{n_e} \pi_j \cdot f_{s(b)} 
  (\sum_{i\in I_b} x_{ib} \mathbb{E}[P_i] + k(j)\cdot m_e),
\]
where we recall that
$m_e$ represents the mean duration of an emergency.
Note that the above expression is in fact
likely to slighly overestimate the optimal
second-stage costs in block $b$, because the emergencies
added at the end of the block contribute to the
overtime, but not to the idle time or waiting time of the block.
}

\paragraph{Reformulation of the two-stage stochastic program as a MIP.} \label{sec:mip}
We can now give a surrogate model for Problem~\eqref{P1} which can be formulated as a mixed integer program (MIP).
To this end, we replace $\phi_b$ with the above convex piecewise linear approximation:

{\small
\begin{subequations}\label{2MIPs}
\begin{align}
\min &\quad  \sum\limits_{i \in I} \sum\limits_{b \in B_{s(i)}'} x_{ib} \cdot c_{ib}+  \changed{\sum\limits_{b \in B }\, \displaystyle\sum_{j=0}^{n_e} \pi_j\cdot  y_{jb}}  \label{2MIPs-objective}\\
\textrm{s.t.} &\quad
\sum\limits_{b \in B'_{s(i)}}   x_{ib}  = 1, && \forall i \in I \label{2MIPs-assignment} \\
&\quad \sum\limits_{b \in B_d}   z_{jb}  = 1, && \forall d \in D , \, \forall j \in \{1,\ldots,n_e\},\ \label{2MIPs-dummyassignment} \\
&\quad \changed{\bar{L}_{jb}} = \displaystyle \sum_{i \in I_{s(b)}} x_{ib} \cdot  \mathbb{E}[P_i] 
\changed{ + \displaystyle \sum_{j'=1}^j z_{j'b}} \cdot m_e & & \forall b \in B, \, \changed{\forall j \in\{0,\ldots,n_e\}}, \label{2MIPs-load}\\
&\quad \alpha_{s(b), k} \cdot \changed{\bar{L}_{jb}} + \beta_{s(b), k}  \leq  \changed{y_{jb}} && \forall k \in \{1,\ldots,R\}, \, \forall b \in B, \,
\changed{\forall j\in \{0,\ldots,n_e\}},
\label{2MIPs-secondstage} \\
&\quad x_{ib}   \in  \{0,1\} && \forall i \in I,\, \forall b \in B_{s(i)}', \label{2MIPs-xdomain}\\
&\quad z_{jb}   \in  \{0,1\} && \changed{\forall j \in \{1,\ldots,n_e\}},\, \forall b \in B, \label{2MIPs-zdomain}\\
&\quad \changed{y_{jb}} \geq 0  && \changed{\forall j \in \{0,\ldots,n_e\}}, \, \forall b \in B. \label{2MIPs-ydomain}
\end{align}
\end{subequations}
}

\changed{
In addition to the $x$ and $z$ variables which define
the assignment of elective and possible emergency patients,
respectively,
we use a variable $y_{jb}$ that represents
our approximation of the
optimal second-stage costs in the event that 
$j\in\{0,\ldots,n_e\}$ emergencies
arrive on day $d=d(b)$. These variables are weighted
by the corresponding probabilities in the objective function~\eqref{2MIPs-objective}.
}
The first constraint~\eqref{2MIPs-assignment} makes sure that each elective surgery \(i \in I \) is allocated to a block \(b \in B_{s(i)} \) of its specialty, or to the dummy block $b'$.
The constraint~\eqref{2MIPs-dummyassignment} ensures that each day,
\changed{
a block is planned for up to $n_e$ emergencies}.
Next, we define the idle-free mean load $\bar{L}_{jb}$
\eqref{2MIPs-load} of a non-dummy-block \(b \in B\)
if $|E_{d(b)}|=j$, and
hence $k(j)=\sum_{j'=1}^j z_{j'b}$ emergencies are assigned to block $b$.
With~\eqref{2MIPs-secondstage}
we enforce $y_{jb}\geq f_{s(b)}(\bar{L}_{jb})$,
so that $y_{jb} = f_{s(b)}(\bar{L}_{jb})$
holds in an optimal solution, as desired.

\paragraph{\changed{Integration of Previous Models.}}
\changed{We wrap up this section by explaining
the workflow that integrates the previously introduced models for the \APP}.
In the preprocessing step, we solve the LP~\eqref{LPstage2}
multiple times for randomly generated instances corresponding
to each specialty, to obtain piecewise linear
approximations of the optimal second-stage costs.
\changed{This preprocessing step has to be done only once, since the resulting
piecewise linear approximations can be used to solve multiple instances.}
\changed{For a given instance, we first have to} solve the MIP~\eqref{2MIPs}, which 
searches for an assignment 
minimizing the first-stage costs and an approximation of the
second-stage costs. 
Finally, for each block we again solve the LP~\eqref{LPstage2} to determine the
optimal tentative starting time of each elective patient that has been assigned to this block
after solving the MIP~\eqref{2MIPs}.
\changed{Note that after this step, the assignment $z$
of emergency cases can be discarded, as our greedy
policy for the online phase uses a different
strategy to insert emergencies in the schedule,
see Section~\ref{sec:onlinephase}.}

\section{Greedy policy for the \OSP}
\label{sec:onlinephase}

In this section we present the policy which we use in our simulations of the \OSP. 
Recall that this policy has to make two types of decisions. First, it must insert
the emergencies arriving on day $d$ in one of the blocks of $B_d$, and second
it should decide to migrate elective surgeries (i.e., rescheduling them to a later
day of the planing horizon or cancelling them) when appropriate.

Our policy is executed in a 
discrete event simulation
which depends on two parameters
$\Delta>0$ and $\alpha>0$
controlling the cancellation
of surgeries and the insertion of
emergencies at non-final positions, respectively.
Its pseudo-code is presented in Algorithm~\ref{alg:alg1}.
In scenario $\omega$, the simulation raises an \emph{event} each time 
a tentative time $T_i^\omega$ is reached, or when a surgery is completed, i.e., 
at a time of the form $S_i^\omega + p_i^\omega$.
Each event is processed by a routine called \texttt{process()} which is sketched in Algorithm~\ref{alg:alg2},
and takes the decisions of starting or migrating a surgery.
In our greedy policy, \texttt{process} is implemented as follows:
\begin{itemize}
 \item If the event that has just been triggered
 is the tentative starting time of a patient $i$, then we simply mark the surgery $i$ as \emph{released}.
 \item Otherwise, the triggering event is that some block $b$ has become available,
 and we have to decide whether to start a surgery on it. We always first try to start
 an elective surgery on $b$ if there {is} a \emph{released} case in $R_b$;  if
 this is not the case, then it means that block $b$ is available until
 the next tentative time and we consider inserting an emergency  surgery in this idle time.
 \item The decision to start an elective surgery depends on the parameter $\Delta$.
 We first calculate an estimation $\tilde{L}_b$ of the load $L_b$, for all $b\in B_d$, conditioned to the knowledge about current state, 
 with a procedure called \texttt{estimate\_expected\_load()}. Then, we check whether
 $\tilde{L}_b\leq T+\Delta$, in which case we start the first released surgery in $R_b$; otherwise, we migrate the last surgery $i\in R_b$ and decide its new assignment by using a procedure 
called \texttt{reassign()}. We always trigger another event right after a migration, so the opportunity to start a surgery on~$b$ is reevaluated using an updated estimate $\tilde{L}_b$.
 \item The decision to insert an emergency in the idle-time until the next elective operation on~$b$
 depends on a parameter $\alpha$.
 We always start an emergency (the one with the longest expectation in $R_e$) if $R_b=\emptyset$.
 Otherwise, we check whether there is an emergency $e\in R_e$
 such that $\alpha\cdot\mathbb{E}[P_e]$
 \emph{fits in the idle time}
 until the next tentative starting time planned on  block~$b$. In that case, we start the longest such emergency.
 \item Whenever no surgery is started upon calling \texttt{process()}, we let the block $b$ idle until the next event.
\end{itemize}

{\small
\IncMargin{1.5em}
\begin{algorithm}[b!]
\nonl\textbf{Input:} assignment $b(i)$ and tentative starting time $t(i)$, for all $i\in B$\\
\smallskip
Initialize $R_b \gets \{ i \in I, b(i) \},\ \forall b\in B$\\
Initialize $T_i^\omega \gets t(i)$ for all scheduled surgeries (with $b(i)\neq b'$)\\
 \For {$d\in D$}
 {
 $R_e \gets E_d^\omega$\\
 Schedule a \texttt{TENTATIVE} event for patient $i$ at time $T_i^\omega$, for all $i\in \bigcup_{b\in B_d} R_b$\\
 Schedule an \texttt{AVAILABLE} event on block $b$ at time $0$ for all $b\in B_d$\\
 \While{$\bigcup_{b\in B_d} R_b \cup R_e \neq \emptyset$}
 {
 Wait for the next \texttt{event}\\
 \texttt{process(event)}
 }
 }
 \caption{\small Greedy Policy in Scenario $\omega$ \label{alg:alg1}}
\end{algorithm}
\DecMargin{1.5em}
}

It remains to explain how
\texttt{estimate\_expected\_load()} and \texttt{reassign()}  are implemented.
When \texttt{estimate\_expected\_load()} is called at time \changed{$\tau$} of day~$d$, we first replace the 
random duration $P_i$ of each surgery $i\in \bigcup_{b\in B_d} R_b \cup R_e$ by a deterministic
value $d_i$, set to the expected remaining surgery time of $i$. More precisely,
for unstarted surgeries we set $d_i=\E[P_i]$ and for ongoing surgeries we use
the conditional expectation $d_i=\E[P_i- \changed{\tau}|P_i > \changed{\tau}-S_i^\omega]$.
Then, for all $b\in B_d$ we start each remaining (deterministic) case in $R_b$ as early as possible,
in the order of their tentative starting time. After this step, it remains to handle the remaining emergency cases in $R_e$. For this  we use a Longest-Processing-Time First (LPT) heuristic: The emergencies in $R_e$ 
are considered by decreasing order of their expected surgery time and added to a block of $B_d$
sequentially, by starting them on the least loaded block. At the end of this procedure, the latest completion of a surgery on block $b$ serves as our estimation $\tilde{L}_b$.

{\small
\IncMargin{1.5em}
\begin{algorithm}[t!]
\nonl\textbf{Input:} \emph{\texttt{event}} that was triggered at current time \changed{$\tau$}\\
\smallskip
 \If{\texttt{event} is of type \texttt{TENTATIVE} for patient $i$}
 {
 	mark patient $i$ as ``released''
 }
 \If{\texttt{event} is of type \texttt{AVAILABLE} in block $b$}
 {
 	\uIf{there is a ``released'' surgery in $R_b$}
 	{
 		$\tilde{L}_b \gets$ \texttt{estimate\_expected\_load()}\\
 		\uIf{$\tilde{L}_b \leq T + \Delta$}
 		{
 		$i \gets$ released surgery in $R_b$ with smallest tentative starting time\\
 		Start $i$ on block $b$\\
 		Remove $i$ from $R_b$ and schedule an \texttt{AVAILABLE} event on $b$ at $\changed{\tau}+p_i^\omega$ 		
 		}
 		\Else{
 		$i \gets$ surgery in $R_b$ with largest tentative starting time\\
		$\tilde{b},\tilde{t}\gets$ \texttt{reassign($i$)}\\
		Move $i$ from $R_b$ to $R_{\tilde{b}}$ and set $T_i^\omega \gets \tilde{t}$\\
		Trigger immediately a new \texttt{AVAILABLE} event on $b$
 		}
 	}
 	\uElse{
 		$\changed{\tau'} \gets$ time of next \texttt{TENTATIVE} event in block $b$ (or $\changed{\tau'}\gets\infty$ if $R_b=\emptyset$)\\
 		\uIf{There is an emergency $e\in R_e$ such that $\alpha \mathbb{E} [P_e] \leq \changed{\tau'}-\changed{\tau}$}
 		{
 		Start $e$ on $b$\\
 		Remove $e$ from $R_e$ and schedule an \texttt{AVAILABLE} event on $b$ at $\changed{\tau}+p_e^\omega$
 		}	
 	}
 }
 \caption{\small \texttt{process}(\emph{\texttt{event}}) \label{alg:alg2}}
\end{algorithm}
\DecMargin{1.5em}
}

The procedure \texttt{reassign(i)} simply works by running a first-fit heuristic: We replace all elective cases on future days by their (deterministic) expected value, and search for the first block on a future day in which the surgery $i$ \emph{fits} without exceeding the nominal capacity $T$ of the block. If such a block $\tilde{b}$ is found,
the case $i$ is rescheduled in~$\tilde{b}$, and
the expected starting time $\tilde{t}$ of $i$ in block $\tilde{b}$ is used as its new tentative starting time.
Otherwise, the surgery $i$ is canceled, i.e., migrated to the dummy block $b'$.

\section{Comparative Framework: Alternative Approaches}\label{sec:other_approaches}
We are not aware of any existing method in the literature for solving the \APP in the setting that
we have presented (simultaneous computation of a patient-to-block assignment and tentative starting times with an objective function including waiting and idle time), even if emergency cases are ignored. We thus
define four other approaches against which we shall compare 
in the computational study of Section~\ref{sec:study}.

\paragraph{Deterministic approach} \label{para:Det}

In this approach, we compute an optimal schedule that results from a deterministic instance in which the random durations are replaced by deterministic ones. This deterministic duration $d_i$ can be set either to the mean processing times $\mathbb{E}[P_i]$,
or for the sake of robustness, to a certain percentile
of the distributions. 
\changed{We have run preliminary experiments with different
percentiles (50th, 60th, 70th, 80-th, 90-th) and found that
using the 70-th percentile produced the best resuts for
the majority of the considered cost structures.
}

With deterministic durations, we can set perfect tentative starting times for each patients, so in this case the idle and waiting times are equal to $0$. To compute an assignment of each elective patient, we thus
solve the following optimization problem:
\begin{subequations}
\begin{align}
\min_{x} &\quad \sum_{i\in I}\sum_{b\in B_{s(i)}'}\ x_{ib}\,  c_{ib}  + \co \cdot \sum_{b\in B} y_b \label{det-obj}\\
&\quad \sum_{b \in B_s'}   x_{ib}  = 1, && \forall s\in S, \forall i \in I_s \label{det-assignment}\\
&\quad y_b \geq  \sum_{i\in I_{s(b)}} x_{ib}\ d_i - T  & & \forall b\in B \label{det-overtime}\\
&\quad y_b \geq 0 && \forall b\in B \label{det-y-domain}\\
&\quad x_{ib} \in \{ 0, 1\} && \forall i \in I, \forall b \in B_{s(i)}.' \label{det-x-domain}
\end{align}
\end{subequations}
The objective function~\eqref{det-obj} is to minimize the sum of the scheduling costs and
the overtime costs. The constraint~\eqref{det-assignment} ensures that
every patient is either scheduled in a block or postponed (allocated to $b'$), and 
constraints~\eqref{det-overtime}-\eqref{det-y-domain} force
$y_b$ to be an overestimator of the overtime in block $b$
(for the problem with deterministic durations $d_i$).

After solving this program, we calculate the tentative starting times.
If the patients allocated to block $b$ are $\{i_1,i_2,\ldots,i_{n_b}\}$, we set
\[
t(i_k) = \sum_{j=0}^{k-1} {d_{i_j}},\quad \forall k\in\{1,\ldots,n_b\}.
\]

\paragraph{First-Fit approach}
In this approach, we also compute an assignment using a deterministic instance where the random processing times are no longer stochastic but replaced by a deterministic duration $d_i$. As for the deterministic MIP above, we found that using the $70$th percentile of $P_i$ for $d_i$ gives \changed{the best results}.

The first-fit heuristic is a well-known approach for
packing problems.
In particular, the \emph{min-weighted sum bin-packing} problem (MWSBP) shares many similarities with the elective surgery planning problem. There, a list of items are given, each with a weight $w_i$ and a size $s_i$, and the goal is to assign items to bins so as to minimize the total cost of the assignment, when assigning item $i$ to the $j$th bin results in a cost $c_{ij}=w_i \cdot j$. The MWSBP problem could hence be used to represent the problem of assigning patients to blocks if overtime was forbidden, there is one block per day in a given specialty, and the goal is to minimize the weighted sum of the dates where the operations take place. For the case of MWSBP, a popular heuristic is to use first-fit, with items sorted 
\emph{in nonincreasing order of the ratio $w_i/s_i$}.
This heuristic is known to be a $2$-approximation algorithm~\citep{epstein2008minimum,sagnol2023}.

In the \APP, the objective function is slightly different, so we have to adapt the order in which the patients are inserted by first-fit in each specialty. Instead of $w_i/s_i$, we use
the ratio $\Delta c_i/d_i$, where $\Delta c_i=c_{ib_2}-c_{ib_1}>0$ 
is the difference between the two smallest values in the set $\{{c_{ib}}\mid b\in B_{s(i)}'\}$.
For each specialty, patients are considered in this order and
sequentially inserted in the first block with enough
remaining capacity (ignoring possibility of overtime) to accomodate
a surgery of duration $d_i$.

\paragraph{Integrated SAA approach.}
Another natural idea is to formulate the \APP as a two-stage stochastic program, 
which
we reformulate into a deterministic equivalent using the SAA approach.
In what follows, 
\changed{we ignore the emergencies so that the problem
decomposes into independent subproblems for each specialty;
this is possible because the only coupling between blocks of different specialties is through emergency cases.
Hence,} we only present a large MIP for the case of a single specialty $s\in S$ only. We note that
it could easily be extended to handle emergencies by introducing variables $z_{e,b}^{\omega_k}$ for the insertion of surgery $e$ at the end of block $b$ in scenario $\omega_k$. Nevertheless, we will see that the proposed MIP is already very hard to solve, so we restrict ourselves to the formulation without emergency cases and for a single specialty.

We assume that the patients of specialty $s$ have been relabeled and ordered according to the SVF-rule, such that $I_s=\{1,\ldots,n_s\}$ and $\operatorname{var}[P_1]\leq \ldots \leq \operatorname{var}[P_{n_s}]$.
\bigskip
\begin{subequations}\label{BMIPs}
\begin{align}
\min_{x}\quad \sum\limits_{i \in I_{s}}  \sum\limits_{b \in B_s'} x_{ib} \, c_{ib} + \frac{1}{ K} \sum_{k=1}^K \Big( \sum_{i\in I_s}\sum_{b\in B_s} \cw \delta_{ib}^{\omega_k}  + \sum_{b \in B_s} \ci ( L^{\omega_{k}}_b -\sum_{i\in I_s} p^{\omega_{k}}_i x_{ib}) + \sum_{b \in B_s} \co \Delta^{\omega_{k}}_b  \Big), \hspace{-13cm} \nonumber\\
\label{BMIPs-objective}\\[1.5em]
s.t. &\quad \sum_{b \in B'_s}x_{ib} = 1, && i \in I_s \label{BMIPs-everypatientblock}\\
&\quad s^{\omega_{k}}_{{i,b}} \geq t_{ib},  && i \in I_s,\quad b \in B_s,\quad k=1,\ldots,K \label{BMIPs-startingtime}\\
&\quad s^{\omega_{k}}_{{i,b}} \geq s^{\omega_{k}}_{i-1, b} + x_{{(i-1)}b}\cdot p^{\omega_{k}}_{i-1} , &&i= 2,\ldots n_s,\quad b \in B_s,\quad  k=1,\ldots,K \label{BMIPs-order} \\
& \quad L^{\omega_{k}}_b \geq s^{\omega_{k}}_{n_s, b} + x_{n_s\, b} \cdot p^{\omega_{k}}_{n_s}, && b \in B_s ,\quad k=1,\ldots,K \label{BMIPs-Load}\\
&\quad \Delta^{\omega_{k}}_b \geq L^{\omega_{k}}_b - T, && b \in B_s , \quad k=1,\ldots,K  \label{BMIPs-overtimeduration}\\
&\quad \Delta^{\omega_{k}}_b \geq 0, && b \in B_s, \quad k=1,\ldots,K \label{BMIPs-overtimegeq0}\\
&\quad \delta_{ib}^{\omega_k}\geq s_{{i,b}}^{\omega_k}-t_{ib} - M\cdot(1-x_{ib}) && i\in I_s,\quad b\in B_s,\quad k=1,\ldots,K \label{BMIPs-deltaWaitingTime}\\
&\quad \delta_{ib}^{\omega_k}\geq 0 &&  i\in I_s,\quad b\in B_s,\quad k=1,\ldots,K \label{BMIPs-delta}\\
&\quad s^{\omega_{k}}_{{i,b}} \geq 0, &&i\in I_s,\quad b \in B_s,\quad k=1,\ldots,K  \label{BMIPs-startover0}\\
&\quad t_{ib} \geq 0, &&i\in I_s, \quad b \in B_s \label{BMIPs-tentativeover0}\\
&\quad x_{ib} \in \{0,1\} , && i\in I_s,\quad b\in B_s \label{BMIPs-xdomain}
\end{align}
\end{subequations}

\bigskip
This large MIP combines the sample average approximation from ~\eqref{LPstage2} for the second stage costs and the assignment-MIP from ~\eqref{2MIPs}, \changed{but restricted to the blocks of
specialty $s$ and without emergencies (so $n_e=0$ and there are no $z$-variables)}. 
The main difference is that we pretend that each patient $i\in I_s$ must be scheduled in all blocks $b\in B_s$, but with a duration of $P_i \cdot x_{ib}$.
Also, we need to introduce some variables $\delta_{ib}^{\omega}$ and $\Delta_b^\omega$ for the delay of surgery $i$ in block $b$ and the overtime of block $b$ in scenario $\omega$, respectively.
The constraints~\eqref{BMIPs-deltaWaitingTime} and~\eqref{BMIPs-delta} ensure that
$\delta_{ib}^{\omega_k} = x_{ib} \cdot \max(0,s_{{i,b}}^{\omega_k}-t_{ib})$ holds
in an optimal solution, so no waiting time is counted for patient $i$ in block $b$ whenever
$x_{ib}=0$. This works because of the large value of the ``big-M'' constant: if 
$x_{ib}=0$, then the right-hand-side of~\eqref{BMIPs-deltaWaitingTime} becomes negative
so $\delta_{ib}^{\omega_k}\geq 0$ has to hold.
In our experiments, we have set $M=1000$. The rest of the constraints are in all points similar to those of~\eqref{LPstage2}--\eqref{2MIPs}.

\paragraph{Benders decomposition approach.}
The main issue with the SAA approach from previous section is the computing time
required to solve Problem~\eqref{BMIPs}, especially for large sample sizes~$K$. To overcome
this issue, we can use a Benders decomposition approach, so {that} the number of samples $K$ only
impacts the solving time of the subproblems for the second stage, which are linear programs (without integer variables). Moreover, with this approach we get rid off the ``big M-constraints''~\eqref{BMIPs-deltaWaitingTime}, which are known to yield weak formulations 
(the fractional solution does not give much information about the optimal integer solution). 
The drawback of this approach is that it can require a large
number of iterations to converge.

\changed{
The main idea is to sequentially construct a convex piecewise linear underestimator of
the second-stage costs. In every iteration,
we solve an approximation of Problem~\eqref{P1} with the current
approximation of the second-stage costs; then, we solve the
dual of Problem~\eqref{P2} to obtain a \emph{cutting plane}
that refines our approximation. 
For the sake of clarity and conciseness, extensive technical details of this approach are relocated to the appendix (Appendix~\ref{sec:benders}).
}

\section{Computational study}
\label{sec:study}

In our computational study, we assume that the availability of blocks is given by the Master Surgery Schedule (MSS) that was given in the description of the AIMMS-MOPTA 2022 competition~\citep{MOPTA22}; see Table~\ref{tab:mss}. 
This MSS is actually taken from~\citet{shehadeh2022data} and appears to be a modification of an MSS first used in~\citet{min2010scheduling}.
There is a total of 32 surgery blocks over the week, each available for a regular working time $T$ of 8 hours. The blocks are allocated to one of six surgery specialties, $s\in S:=\{\texttt{CARD}, \texttt{GASTRO}, \texttt{GYN}, \texttt{MED}, \texttt{ORTH}, \texttt{URO}\}$,
namely cardiology, gastroenterology, gynaecology, general surgery,
orthopedics and urology. 

\subsection{Assumptions for data generation}\label{sec:assumptions}
\paragraph{Cost structures.}
In our study, we consider six different types of structures 
for the waiting, idling and overtime costs, denoted by $cs_1,\ldots,cs_6$.
In all settings we set $\co=1$, which means that one can interpret cost values as
\emph{equivalent minutes of overtime}.

In the first two cost structures, we assume that the idling costs are set to $\ci=0$, because idling already results in either more overtime or less scheduled patients, so there is no need to penalize it even more.  Such a setting was considered by~\citet{Berg2014Denton}, who
further vary the waiting costs $\cw$ between $\frac{1}{33}$ (a constant coming from
an estimation by specialists) and 1. In this paper, we consider the two extreme cases of 
setting the waiting cost to either $\cw=0$ in $cs_1$ (only overtime is considered) or 
$\cw=1$ in $cs_2$ (waiting time costs as much as overtime). 
Note that when the appointment scheduling aspects of the problem disappear whenever 
$\cw=0$, because in that case it is optimal to set all tentative times to $t_i=0$ in the
appointment problem (so all operations are done one after the other without idling time).

In the cost structure $cs_3$ we take the parameter values as $\co=1, \cw=2, \ci=2$,
that is, we assume that the cost of idling and waiting are the same, and equal to twice
the cost overtime. This is similar to the setting used in~\citet{denton2007optimization}.

The cost structures $cs_4$ and $cs_5$ are taken from~\citet{deceuninck2018outpatient}
and assume that $\co=1$ and $\ci=\frac{2}{3}$, i.e., 
overtime costs 50\% more than idling time. The
waiting costs are set to $\cw=\frac{2}{15}$
and $\cw=\frac{2}{3}$ in $cs_4$ and $cs_5$, 
which corresponds to idling costs that are 5 times or as
large as waiting costs, respectively.

Finally, we also consider the cost structure $cs_6$, in which
the importance of overtime, waiting time and idling time are equal,
i.e., $\co=\cw=\ci=1$.

\medskip
As mentioned in the introduction, the scheduling costs $c_{ib}$ can be used to model
a variety of scheduling objectives or preferences. In this work, we use these costs
to incentivize a fair treatment of all patients, 
that is, we want to favor schedules in which 
the total waiting time of a patient is 
inversely proportional to the severity of his or her condition.
Here, the patient waiting time is the total a patient 
\emph{remains in the system}, i.e., the delay between the day a specialist assessed the necessity of doing a surgery and the day where it finally takes place.
In scheduling terminology, this delay is called \emph{flowtime} $F_i$, and it is equal to $t(b)+e(i)$ for a patient $i$ operated in block $b$, 
where $t(b)$ is time associated with block $b$ and $e(i)>0$ denotes the entry-time of patient $i$ in the system. 
Many authors proposed to set the cost $c_{ib}$ of the form $w_i \cdot F_i$, where the weight $w_i$ 
of a patient is related to its priority~\citep{Lamiri2008Xie, jebali2015stochastic}. However, 
because of the linearity with respect to the flowtime, such a cost structure does not distinguish between a schedule in which two patients $i$ and $j$ with the same weight
have flowtimes $F_i=F_j=3$, and a schedule where $F_i=1$ and $F_j=5$.
In order to favor the former, fairer solution, we consider the weighted quadratic flowtime, i.e., we set
\begin{align*}
c_{ib} &= w_i \; F_i^ 2 = w_i  \; \Big( t(b) + e(i) \Big)^2, \forall i\in I,\ \forall b  \in B.
\end{align*}
Note that this objective (weighted squared flowtime) tries to balance the flowtime of each patient: Given a constant total flowtime for all patients, simple calculations show that an ideal schedule would yield a flowtime proportional to $\frac{1}{w_i}$ for patients with weight $w_i$.

To set the postponing costs $c_{ib'}$ associated with the dummy block $b'$, we
use the following two natural conditions: On the one hand, $c_{ib'}$ should be larger
than the cost of assigning $i$ to any block $b\in B$ of the planning horizon.
On the other hand, it should be less expensive to postpone a patient than scheduling him or her \emph{completely in overtime}, in any block of the planning period. According to our study in Section~\ref{sec:offlinephase},
the total second-stage costs in block $b$ of specialty $s$ can be approximated by a piecewise linear function $f_s$
of the total expected duration assigned to this block. Denote by $\bar{\alpha}_s = \max_{k=1,\ldots,R} \alpha_{sk}$
the coefficient of the rightmost ``piece'' of the function $f_s$, i.e., the slope of $x\mapsto f_s(x)$ for large values of $x$. In the worst case, adding a patient in block~$b$ increases
the total cost by $c_{ib} + \bar{\alpha}_{s(b)} \cdot \mathbb{E}[P_i]$. We thus require that
\begin{equation} \label{bounds_dummy_block_cost}
\max_{b\in B} c_{ib} \leq c_{ib'} \leq  \min_{b\in B} c_{ib} + \bar{\alpha}_{s(b)} \cdot \mathbb{E}[P_i]
\end{equation}
holds for all patients $i\in I$. In practice we set $c_{ib'}$ equal to the average of
the above bounds, that is,\linebreak \mbox{$c_{ib'}:= \frac{1}{2} (\max_{b\in B} c_{ib} + \min_{b\in B} c_{ib} + \bar{\alpha}_{s(b)} \cdot \mathbb{E}[P_i])$}. (The lower bound $\max_{b\in B} c_{ib}$
was always smaller than the upper bound $\min_{b\in B} c_{ib} + \bar{\alpha}_{s(b)} \cdot \mathbb{E}[P_i])$ in our experiments.)

In our study, we considered two settings $u\in\{\texttt{day}, \texttt{week}\}$, differing in the time unit used to measure the flowtime. In the first setting, the flowtime is measured in days. In that case, the entry-time of each patient is drawn at random in $\{1,\ldots,7\}$,
that is we assume that patients are waiting for an operation for less than a week 
and $t(b)$ is set equal to the index of the day of block~$b$ (0 for Monday, 1 for Tuesday, etc.) for all $b\in B$.
In the second setting, the time is measured in weeks, so $t(b)=0$ for all blocks $b\in B$
and there is no advantage to operate a patient on Monday rather than Friday. Under this setting, $e(i)$ is drawn an random in $\{1,2\}$, i.e., we assume that patients have already been waiting for 1 or two weeks at the beginning of the time period.

We assume that the patient weights vary within a ratio of $1$ to $4$, hence for some constant $w_0$ we draw $w_i$ at random in the interval $[w_0, 4w_0]$. The parameter $w_0$ serves to define the tradeoff between the scheduling costs and the other costs, it can thus be seen
as a parameter controlling our aversion to overtime.
Setting it to a large value would violate the ordering of the lower and upper bounds
in~\eqref{bounds_dummy_block_cost}, while setting $w_0$ to a small value would put the
focus only on the second-stage costs, thus favoring \emph{empty schedules} where too many patients
are postponed. In our experiments, setting $w_0=0.05$ if $u=\texttt{day}$ and $w_0=1$ if $u=\texttt{week}$ yields 
\changed{a reasonable tradeoff between the overtime and the number of cancellations}.

Finally, we also assume that $\cm = 120$, 
meaning that a patient migration costs as much as 2 hours of overtime.

\paragraph{Generation of probability distributions
for the surgery durations.}
It remains to explain how we generate random parameters $\mu_i$ and $\sigma_i$ for the lognormal probability distribution of
the duration $P_i$ for each patient.
We use the 
the marginal mean $m_s$ and marginal variance $v_s$ of surgeries
given in the description of the AIMMS-MOPTA competition~\citep{MOPTA22}; see Table~\ref{tab:sv}.
In addition, we set $m_e=90$ and $v_e=70^2$ for emergency patients $e\in E$, so $P_e$ has a large probability to fall in the interval between 1 and 3 hours.

\begin{table}[t!]
	\centering	
 \caption {{\small Marginal mean and standard deviation of surgery duration (in minutes), for each specialty
 and for emergency patients.
 \label{tab:sv}}}
 \smallskip
 {\small
\begin{tabular}{crr}
\toprule
 Specialty $s$ & Mean $m_s$ &  Standard dev.\ $\sqrt{v_s}$\\
 \midrule
\texttt{CARD}  & $99$ & $ 53$ \\
\texttt{GASTRO}  & $132$ & $76$ \\
\texttt{GYN}  & $78$ & $52$ \\
\texttt{MED}  & $75$ & $72$ \\
\texttt{ORTH}  & $142$ & $ 58$  \\
\texttt{URO}  & $72$ & $38$  \\
\texttt{EMER}  & $90$ & $70$  \\
\bottomrule
\end{tabular}
}
\end{table} 

\begin{table}[t]
\caption{\small Example of MSS used in our computational study~\label{tab:mss}}
 \begin{center}
  \begin{tabular}{rccccc}
  \hline
   \textbf{OR} & \textbf{Monday} & \textbf{Tuesday} & \textbf{Wednesday} & \textbf{Thursday} & \textbf{Friday}\\ \hline 
   1 & \texttt{GASTRO} & \texttt{GASTRO} & \texttt{GASTRO} & & \\[-1mm]
   2 &        &        & \texttt{GASTRO} & \texttt{GASTRO} & \texttt{GASTRO}\\[-1mm]
   3 & \texttt{CARD}   &        & \texttt{CARD}   &        & \texttt{CARD} \\[-1mm]
   4 & \texttt{ORTH}   & \texttt{ORTH}   &        & \texttt{ORTH}   & \texttt{ORTH}\\[-1mm]
   5 &        & \texttt{ORTH}   & \texttt{MED}    &        & \\[-1mm]
   6 & \texttt{GYN}    & \texttt{GYN}    & \texttt{GYN} & \texttt{GYN} \\[-1mm]
   7 &        & \texttt{GYN}    & \texttt{GYN}    & \texttt{GYN}    & \texttt{GYN}\\[-1mm]
   8 & \texttt{URO}    & \texttt{URO}    &        & \texttt{URO}    & \texttt{URO}\\[-1mm]
   9 & \texttt{CARD}   &        & \texttt{URO}    &        & \texttt{CARD}\\[-1mm]
   10& \texttt{URO}    &        & \texttt{ORTH}\\\hline
  \end{tabular}
 \end{center}
\end{table}

\changed{
In practice, it is possible to obtain a more precise 
estimation of the surgery duration of an individual 
patient than simply saying
that a patient's duration follows the marginal probability
law of all patients in his or her specialty. This is due to 
the fact that the surgeons know the exact procedure
to perform for this patient, and can give an estimate of the required OR time.
More generally, a learning algorithm could provide an estimation
for the parameters $\mu_i$ and $\sigma_i$ such that 
$P_i\sim\mathcal{L}\mathcal{N}(\mu_i,\sigma_i)$,
\emph{based on observed features of the patient} $i$.
The approach we follow to produce parameters $\mu_i,\sigma_i$
yielding probability distributions 
that are \emph{less dispersed} but still consitent with 
the marginals, is explained in details in Appendix~\ref{sec:gen-dist}.
}

\subsection{Results} \label{sec:results}

\begin{figure}[t]
\centering
		\includegraphics[trim=5 5 5 10, clip,
		scale=0.75, width=8cm]{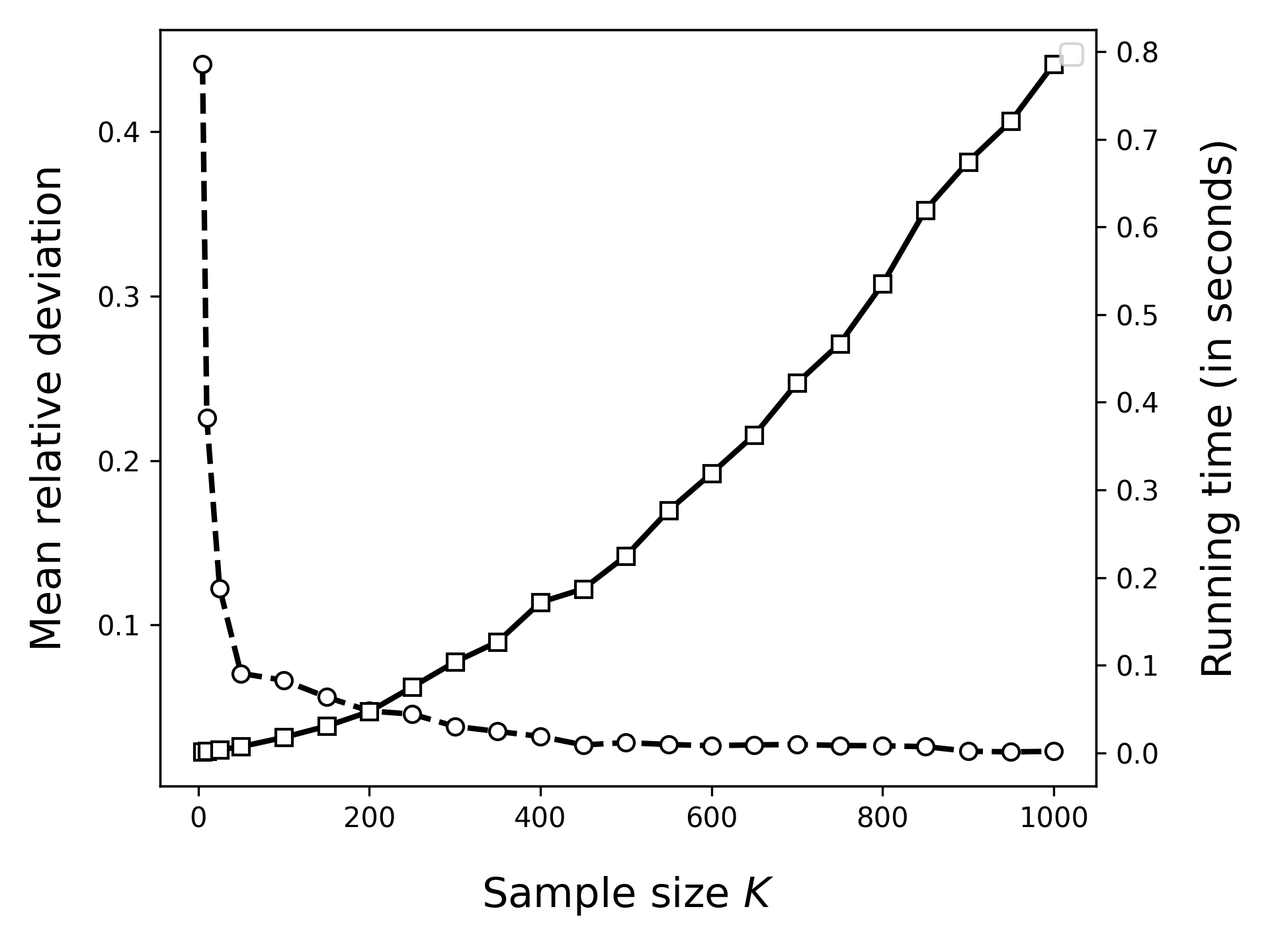}\vspace{-1em}
\caption{\small Impact of sample size on the accuracy of the SAA (dashed line, left axis) and on the CPU time for solving~\eqref{LPstage2} (solid line, right axis).}
\label{effect_of_K}
\end{figure}

\paragraph{Choice of sample size $K$.}
The quality of the sample average approximation (SAA) depends on $K$. While larger values of $K$ yield a better approximation,
they also cause longer computing times. In our case study, we simply selected the value of $K$ to achieve a good tradeoff between accuracy and computing time.
Figure~\ref{effect_of_K} shows the effect of the sample size on the SAA when we solve the LP for the second stage problem~\eqref{LPstage2}. We generated 10 different instances, each consisting of 6 patients from the \texttt{MED} specialty (the number of $6$ patients is chosen to approximately match the regular working time of 8 hours in the block).  We measure the accuracy of the SAA by comparing the
optimal value of Problem~\eqref{LPstage2} to the value obtained by setting $K$ to a very large value ($K=10\, 000$).
The plot shows both the computing time for solving the LP~\eqref{LPstage2} and the mean relative deviation across the 10 instances (measured in terms of objective function value) to the solution for $K=10\, 000$, for $K$ ranging from 10 to 1000. The figure seems to suggest that the accuracy decreases very slowly after $K=450$, hence our choice of $K$ for the computational study.

\paragraph{Experimental setup.}

In our case study, we generate random instances 
corresponding to 50 different seeds for the random number generator, for every combination of number of elective patients $n\in\{70,100,140,200\}$, emergency rate $\lambda\in\{0,1,2,3\}$,
flowtime unit $u\in\{\texttt{day},\texttt{week}\}$
and the six cost structures for $(\co,\cw,\ci)$ discussed in Section~\ref{sec:assumptions}, resulting in $192$ different instances for each seed. 

We evaluate the different approaches with the help of Monte Carlo simulations. For this purpose, we generate for every instance a validation set of $K=450$ scenarios. All MIPs and LPs have been solved with Gurobi~\citep{gurobi} on a server with an AMD EPYC 7302 processor with 16 cores at 1.5GHz.
For each instance,
we have run the greedy policy of Section~\ref{sec:onlinephase} with 
4 different values of $\alpha$ and 3 values of $\Delta$,
as indicated in Table~\ref{tab:varyingparams}, resulting
in 12 simulations per instance.

\begin{table}[t!]
	\centering	
 \caption{\small Varying parameters in our computational study}
 
\medskip 
\begin{tabular}{ll}
\toprule
Parameters & Values\\
\midrule
instance seed & $\{1,2,\ldots,50\}$\\[0.2em]
number $n$ of patients  & $\{70, 100, 140, 200\}$  \\[0.2em]
emergency rate $\lambda$   & $\{0,1,2,3  \}$ \\[0.2em]
flowtime unit $u$ & $\{\texttt{day}, \texttt{week} \}$\\[0.2em]
cost structure & $\{cs_1,cs_{2},cs_{3},cs_{4},cs_5,cs_6\}$\\[0.2em]
cancellation parameter $\Delta$& $\{60,120, 1000 \}$\\[0.2em]
emergency insertion parameter $\alpha$& $\{0.6, 0.7, 0.8, 0.9 \}$\\[0.2em]
\bottomrule
\end{tabular}

\label{tab:varyingparams}
\end{table} 

\paragraph{Effect of the simulation parameters $\alpha$ and $\Delta$.}

\begin{figure}[t!]
	\centering
	\includegraphics[width=0.49\linewidth]{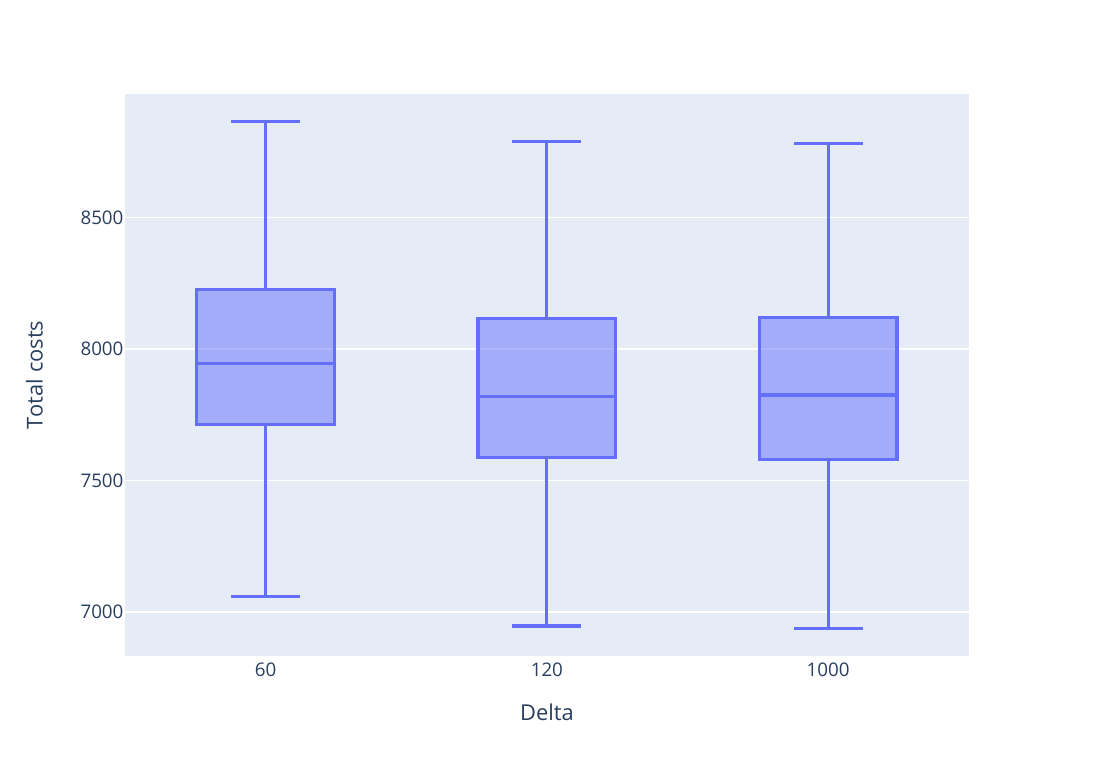}
	\includegraphics[width=0.49\linewidth]{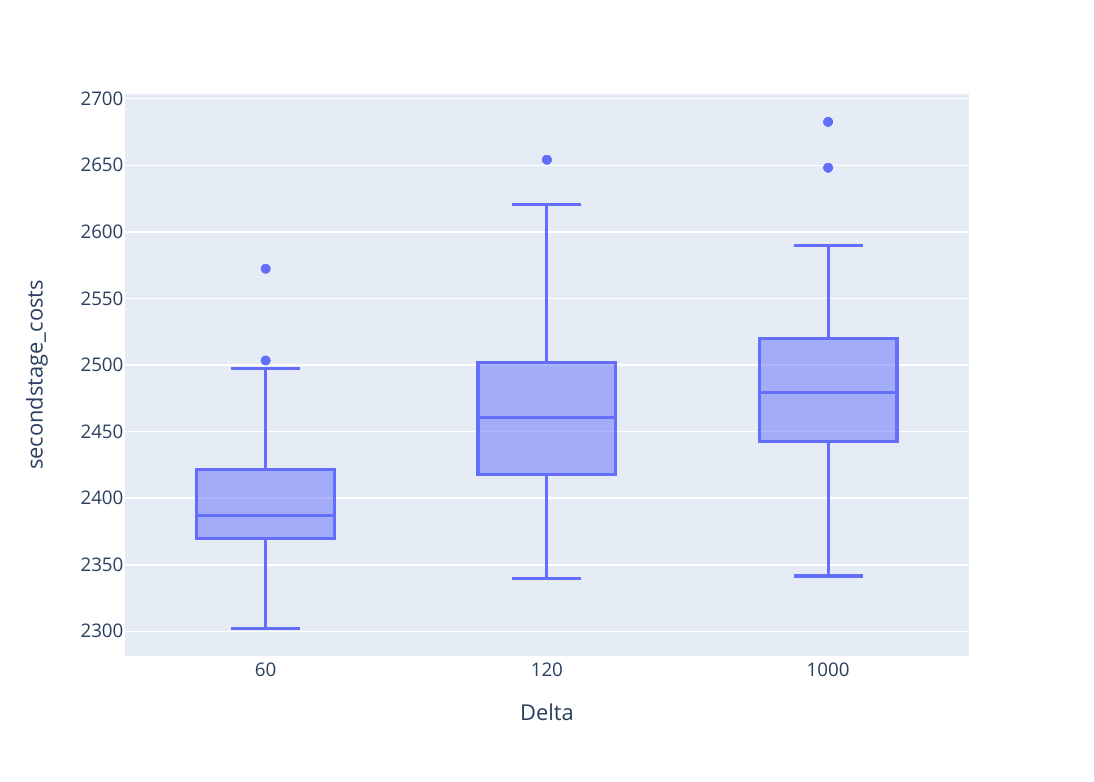}
	\caption{\small Results of simulations for different settings of $\Delta$, for the cost structure $cs_6\ (\co=\cw=\ci=1)$, 
	with $n=200$ patients and emergency rate $\lambda=3$
	in the setting $u=\texttt{day}$. 
	The assignment of patients was computed with the 
	surrogate model-based two-stage stochastic programming approach of Section~\ref{sec:offlinephase}, considering a maximum of $n_e=10$ emergencies per day. The left figure shows boxplots of the total costs, while the right figure shows only the second-stage costs (idling, waiting, overtime).
	\label{fig_delta}}
	\end{figure}

Figure~\ref{fig_delta} shows that the use of our cancellation policy with $\Delta=60$ substantially improves the second-stage costs, but also causes more patient migrations, which results in
larger total costs. Note that the high setting of $\Delta=1000$
basically disables the possibility of cancelling patients during
the online phase.
However in practice, it might be an obligation to use some kind of cancellation rule, as the medical staff cannot do unlimited
overhours. Therefore, for the remainder of this computational study we only consider the results for $\Delta=120$, which
seems to be an acceptable tradeoff. For this value of $\Delta$,
cancellations are very rare, which is due to our setting
of $\cm$ and $c_{ib'}$.

Our experiments also suggest that the parameter $\alpha$ has a negligible impact on the total costs of the solution, but
the reader can check with the demonstrator (see Section~\ref{sec:demonstrator}) that smaller values of $\alpha$ 
yield more insertions of emergencies at non-final positions.
We set $\alpha$ to $0.7$ for the rest of this computational study.

\paragraph{Effect of flowtime unit $u$.}

\begin{figure}[t!]
\centering
	\includegraphics[width=0.8\linewidth]{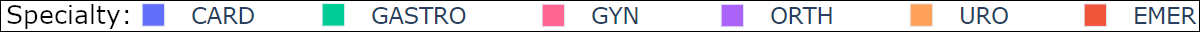}
	
	\bigskip

	\begin{minipage}{0.495\textwidth}
		\includegraphics[clip, trim=0cm 0cm 2.4cm 0cm,width=\textwidth]{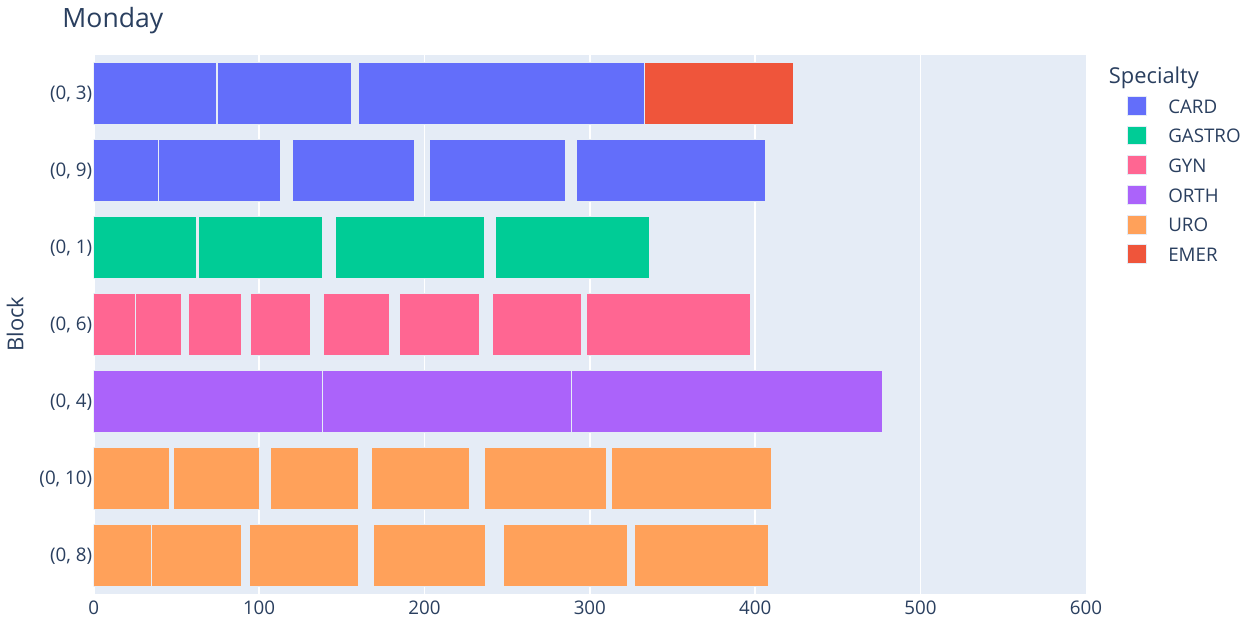}

	\end{minipage}
		\begin{minipage}{0.495\textwidth}
		\includegraphics[clip, trim=0cm 0cm 2.4cm 0cm,width=\textwidth]{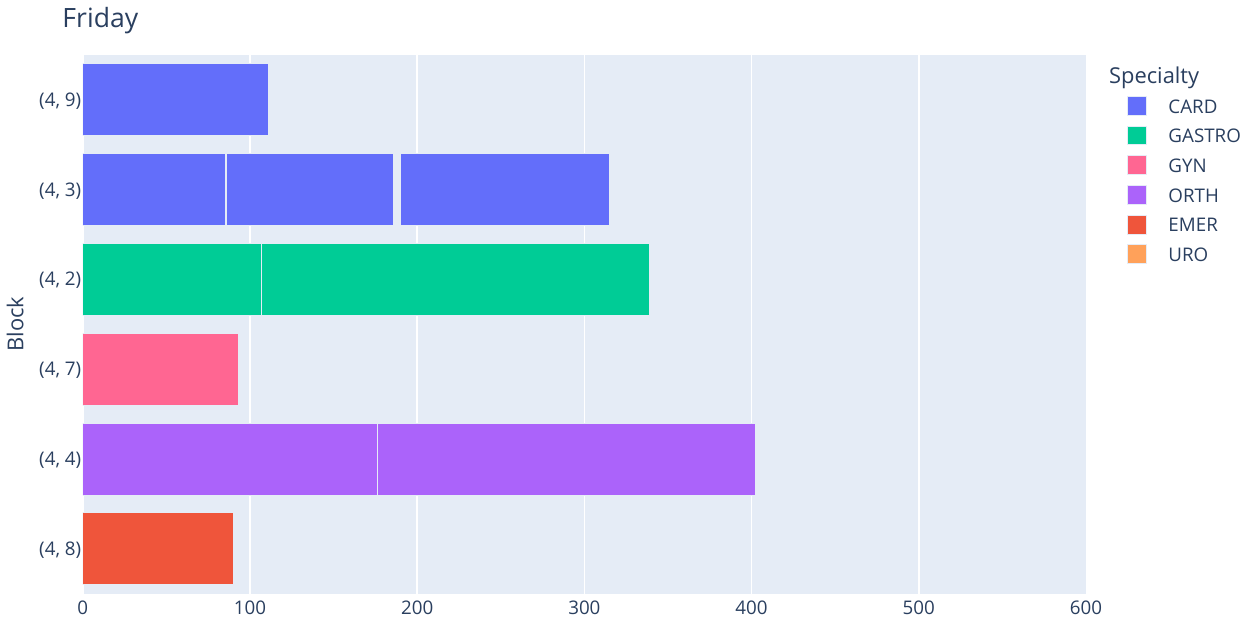}

	\end{minipage}
	
	\medskip
	{\small $u=\texttt{day}$}
	
	\vspace{1em}
	\begin{minipage}{0.495\textwidth}
		\includegraphics[clip, trim=0cm 0cm 2.4cm 0cm,width=\textwidth]{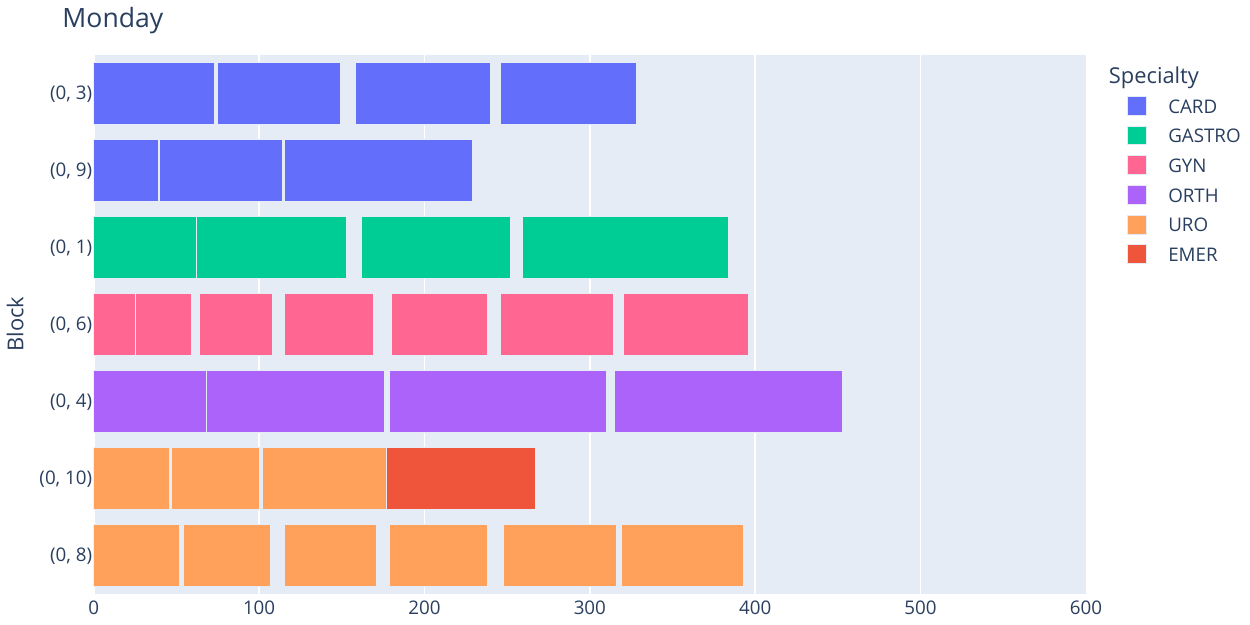}

	\end{minipage}
		\begin{minipage}{0.495\textwidth}
		\includegraphics[clip, trim=0cm 0cm 2.4cm 0cm,width=\textwidth]{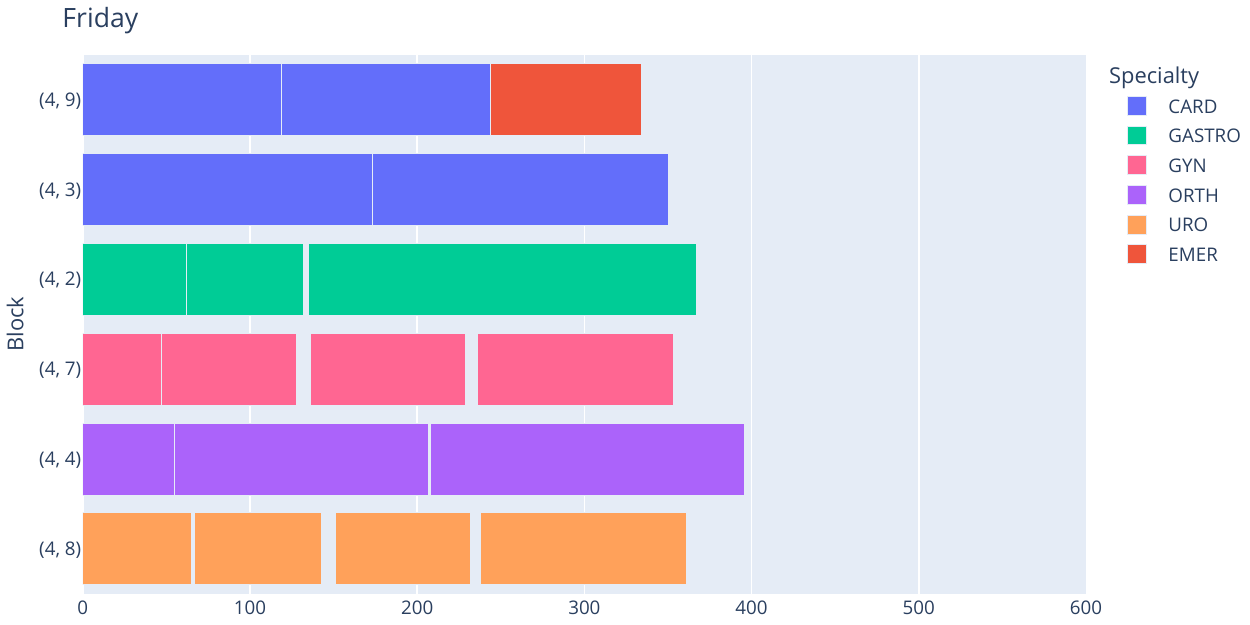}

	\end{minipage}
	
	\medskip
	{\small $u=\texttt{week}$}
\caption{\small Computed assignment (displayed with mean durations) on Monday and Friday only, for an instance with $n=100$ patients, under the setting $u=\texttt{day}$ (first row) and $u=\texttt{week}$ (second row).
\label{fig:gantt}}
\end{figure}

Figure~\ref{fig:gantt} shows the 
assignment and tentative starting times computed with the surrogate model-based two-stage stochastic
method (with a maximum of $n_e=5$ considered emergencies), for 
an instance with $n=100$ patients, rate of $\lambda=1$
emergency per day
and the cost structure $cs_6$ ($\co=\cw=\ci=1$).

Two different settings for the flowtime unit
parameter $u$ are compared. The picture displays the blocks
of Monday and Friday only for the sake of visibility,
which allows us to compare how the load is balanced over the week. 
The Gantt charts show the schedule in the case where
mean scenario (where the duration of each surgery is equal
to its expected value); Idle time occurs when the tentative
starting time is larger than the completion time of the previous operation. Note that since $\lambda=1$, the graph also indicates one emergency per day (with mean duration $m_e=90$), added on the least loaded machine. (Recall that the number of emergencies per day must not be equal to $\lambda=1$ as $|E_d|$ is a Poisson random
variable, but we only display the mean scenario.)

The figure shows that the algorithm produces an assignment which is very dense on Monday and almost empty on Friday with the setting $u=\texttt{day}$, as in that case the scheduling costs incentivize to schedule the emergencies as early as possible. In contrast, with $u=\texttt{week}$ we obtain an assignment of surgeries which is much more balanced over the week.

\paragraph{Comparison of assignment algorithms.}

\begin{figure}[ht!]
\centering
\begin{tabular}{cc}
 \multicolumn{2}{c}{\includegraphics[width=0.8\linewidth]{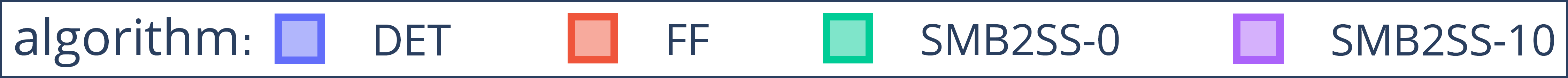}}\\
 \hspace*{-6mm}\includegraphics[clip, trim=0cm 0cm 2.8cm 0cm,width=0.52\linewidth,
 height=5.4cm]{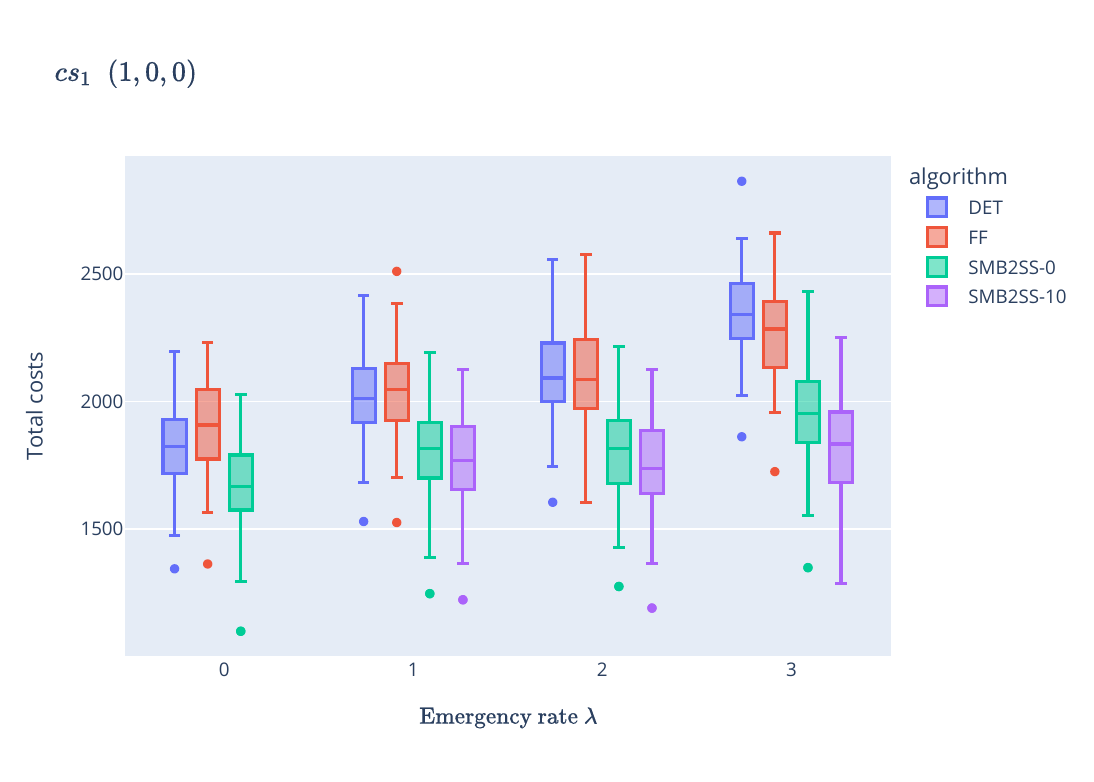}
  &
 \hspace*{-6mm}\includegraphics[clip, trim=0cm 0cm 3.1cm 0cm,width=0.52\linewidth,
 height=5.4cm]{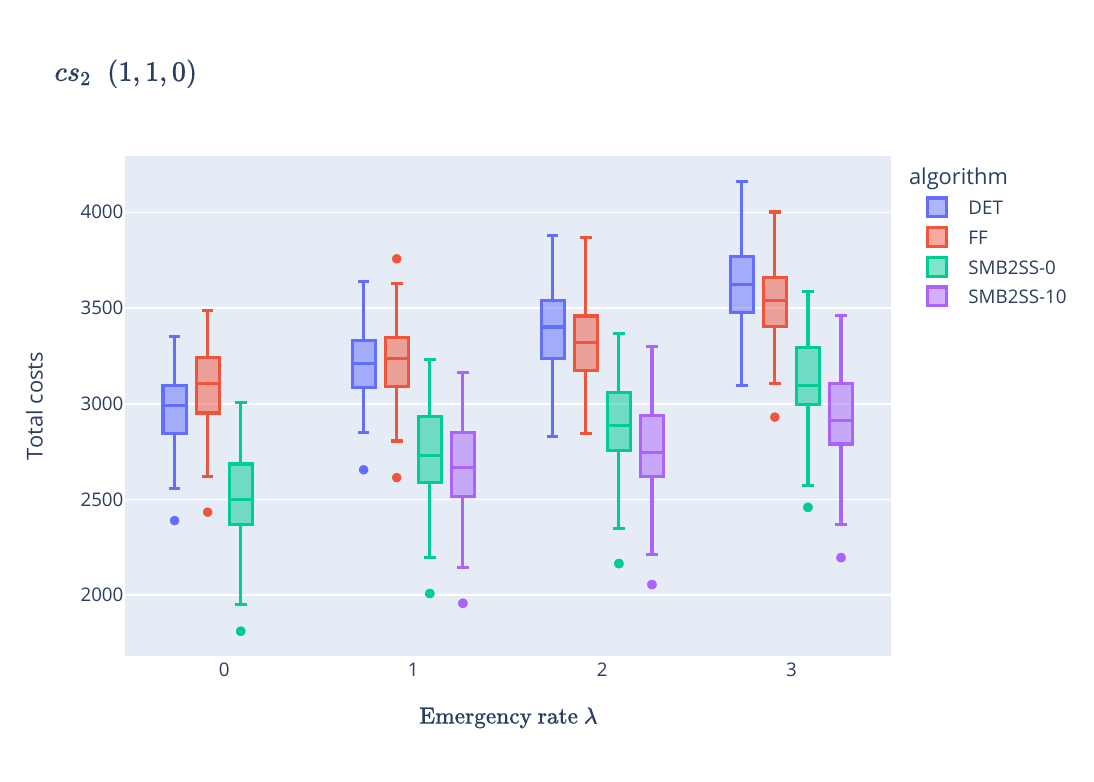}\\[-1.2em]
 \hspace*{-6mm}\includegraphics[clip, trim=0cm 0cm 2.8cm 0cm,width=0.52\linewidth,
 height=5.4cm]{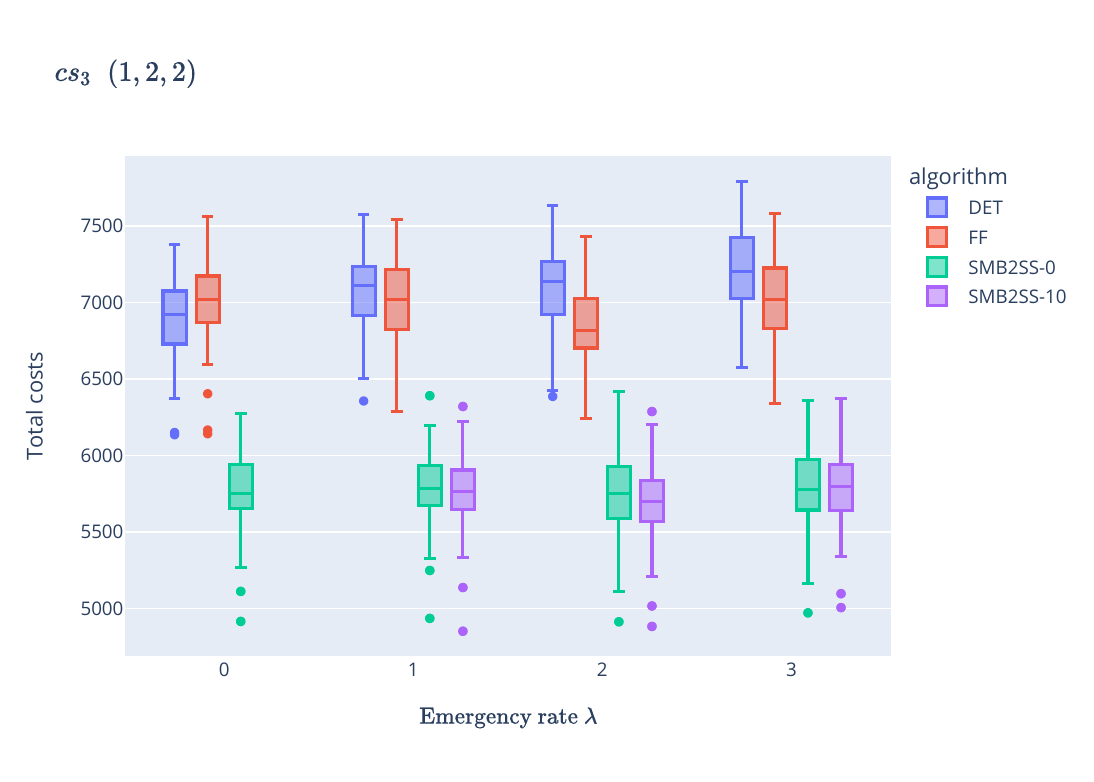}
  &
 \hspace*{-6mm}\includegraphics[clip, trim=0cm 0cm 3.1cm 0cm,width=0.52\linewidth,
 height=5.4cm]{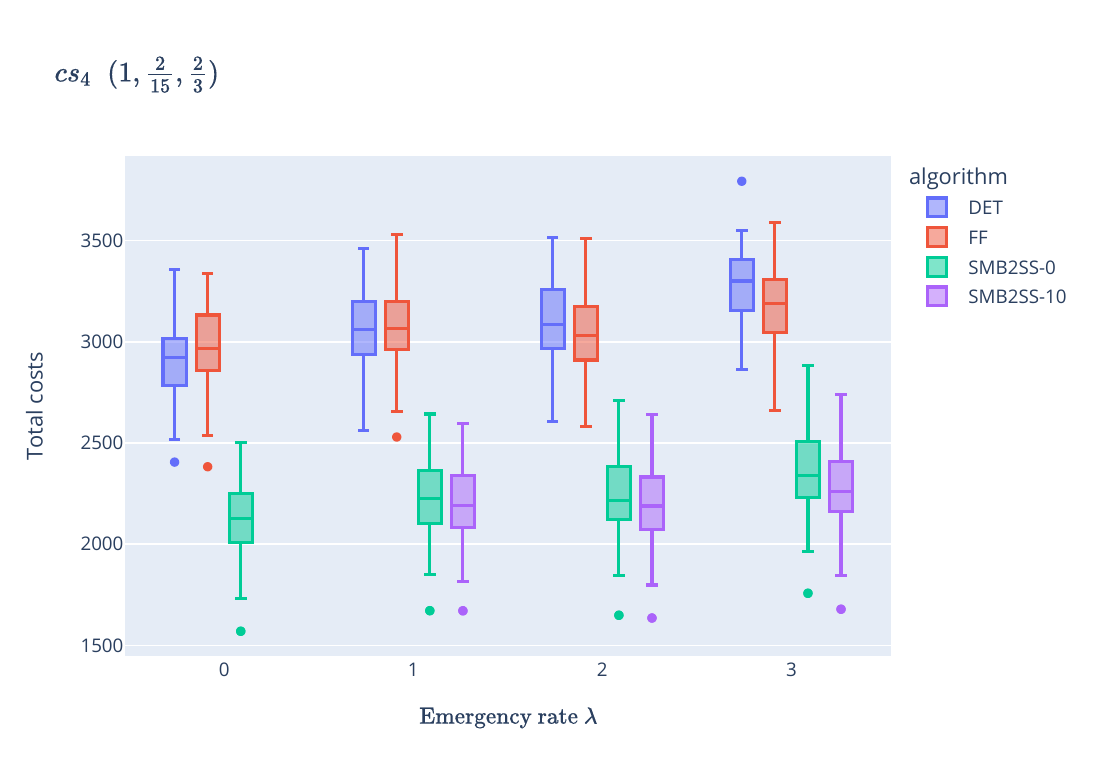}\\[-1.2em]
  \hspace*{-6mm}\includegraphics[clip, trim=0cm 0cm 2.8cm 0cm,width=0.52\linewidth,
 height=5.4cm]{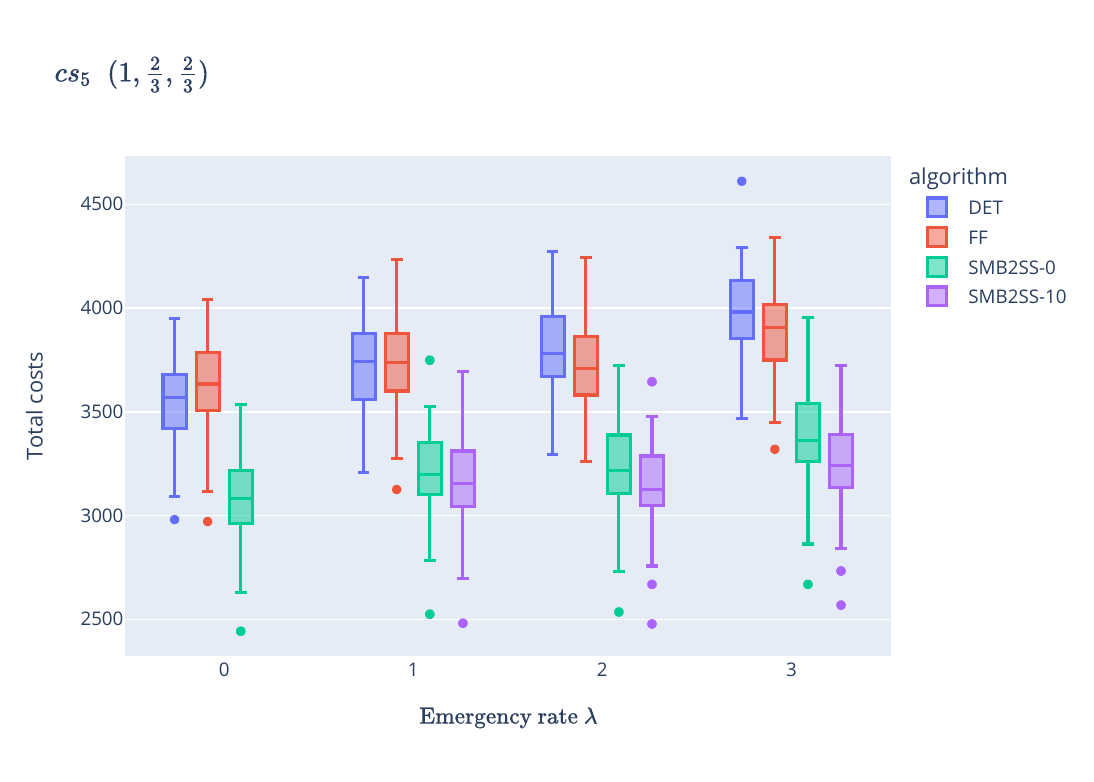}
  &
 \hspace*{-6mm}\includegraphics[clip, trim=0cm 0cm 3.1cm 0cm,width=0.52\linewidth,
 height=5.4cm]{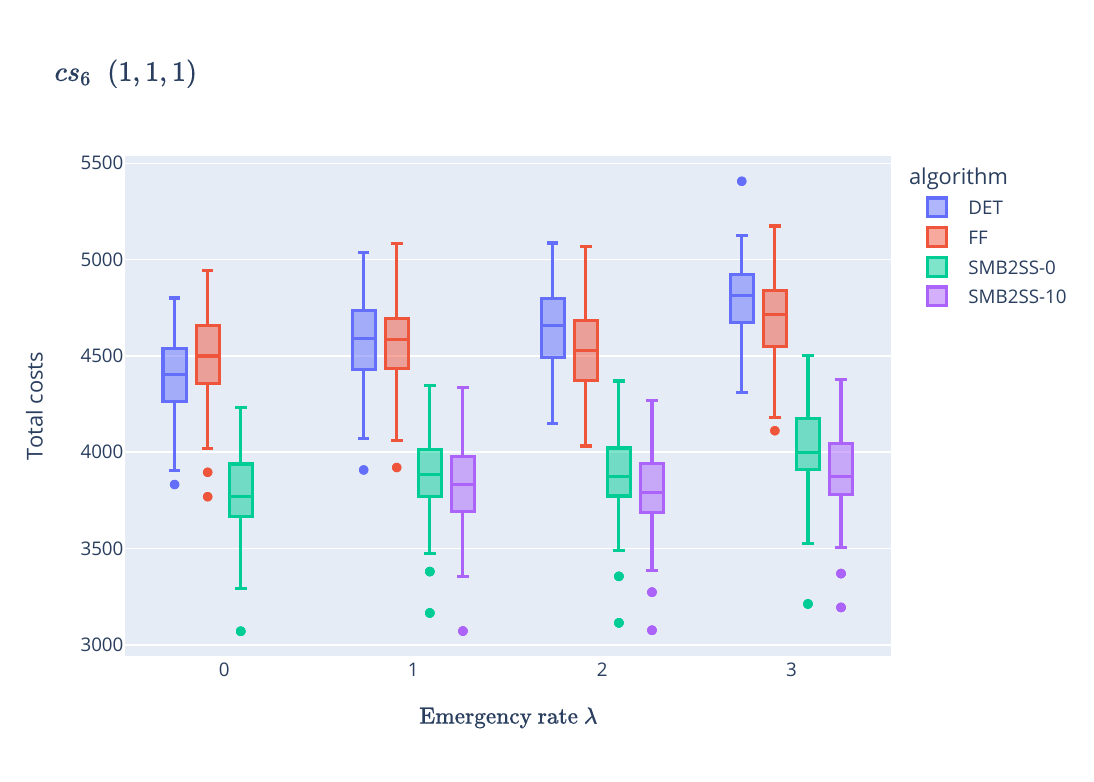}
\end{tabular}

\caption{\small Boxplots comparing the distribution of 
the total costs obtained by Monte-Carlo simulation,
for the assignment of patients computed with $5$ different algorithms, for different values of the emergency
rate $\lambda$ (x-axis) and 6 different cost structures
(one picture for each cost structure).
In all these experiments, the number of patients
is $n=140$ and the flowtime unit is $u=\texttt{day}$.
\label{fig:costs_structure}}
\end{figure}

We first compare our approach to fast heuristic approaches.
Figure~\ref{fig:costs_structure}
shows the distribution of the total costs over all instances with $n=140$
patients and flowtime unit $u=\texttt{day}$,
for each of the six cost structures 
and values of the emergency rate.
Four approaches are compared: the deterministic approach (DET),
the first-fit method (FF), and
the surrogate model-based two-stage stochastic programming approach with parameter $n_e \in \{ 0, \changed{10}\}$ (SMB2SS-0 and SMB2SS-\changed{10}).\enlargethispage{1em}

As can be seen on the figure, our 2-stage stochastic programming approach yields significantly better
results than the other approaches, for all values of $\lambda$
and all cost-structures. Also, taking $n_e=\changed{10}$ emergencies into accounts improves the total costs when $\lambda$ is large,
as expected. 
The improvement depends a lot on the cost structure, though. In the cost structure $cs_1$, where only the overtime is considered, SMB2SS-0 achieves 
$7.9\%$ improvement over DET \changed{on average}
when there are no emergencies ($\lambda=0$) and \changed{SMB2SS-10} has \changed{$22.4\%$} improvement over DET when $\lambda=3$. 
In the cost structure \changed{$cs_4$}, these improvements are of \changed{$27.0\%$} and \changed{$30.8\%$}, respectively.

\bigskip
We next compare the solutions obtained with our approach (SMB2SS)
to the solutions obtained by the Integrated Sampling Average Approximation (SAA) or with the Benders-decomposition approach. As already mentioned, these approaches require
a lot of computation time. Therefore, we restrict ourselves
to analyzing the case without emergencies ($\lambda=0$),
so the problem decomposes in 6 independent subproblems (one
for each specialty), and we set a time limit of 5 minutes for
solving each subproblem. 
For the SMB2SS approach, the time limit was set to $20s$ for each specialty.
To make the large MIPs easier to solve, we also consider running SAA and Benders with only $K=50$
scenarios.

The average results on 50 randomly generated instances 
with different values of $n$ 
are displayed in Table~\ref{benders_cs2} and
Table~\ref{benders_cs3} for the cost structures $cs_2$ and
$cs_3$, respectively.
The results for Benders with $K=50$ are not displayed, because
they are of worse quality than with $K=450$, for a similar 
runtime. Unlike SAA, one advantage of the Benders-decomposition
approach is indeed that it is little affected by the number of scenarios, as this only impacts the time required to solve the
appointment problems, which are linear programs (without integer variables).

As we can see,
the SMB2SS approach requires several orders of magnitude less runtime than SAA or Benders. It was able to solve all instances
with a gap $<0.2\%$ within 8s for $n=70$ and 130s for $n=200$, 
while SAA and Benders still have large gaps after 30 minutes of computation. 
The SAA approach with $K=450$ was not even
able to find solutions within the allocated time limit
as soon as $n\geq 100$. We also observe that SAA beats our 2SS-approach on small instances, but our approach yields the
best results for larger instances.

Also, note the difference between the columns called ``cost'' and
``sim.''. The former indicates the cost in the objective value
of the considered MIP, while the latter are the results
of the Monte-Carlo simulation with our greedy online algorithm
on the validation set of $K=450$ scenarios.  For SAA and Benders, the difference between these two columns mostly comes from the fact that different set of scenarios are used (a training set for solving the MIP, and the validation set to run the simulation). In the SMB2SS approach, we see that the MIP value
is very close to the value obtained by simulation, which
indicates that our surrogate model is very accurate.

\paragraph{Near-optimality of SVF-rule.} 
We also ran additional experiments to confirm that the shortest-variance first (SVF) rule is a good choice for ordering the patients. For a fixed block $b$, we generate $N=100$ different random instances, with a number of patients $n_b$ such that the expected load of the block lies in the interval $[300,550]$, which is a reasonable assumption, given that the regular working time of a block is set to $T=480$ in our study.
Then, we solve the LP~\eqref{LPstage2}
for each of the $(n_b!)$ ordering of the $n_b$ patients in the block. Table~\ref{tab:svf-opt} shows that for the considered cost structure ($cs_5,\ (\co=1, \cw=\ci=\frac{2}{3})$), the SVF-rule is often optimal, and if it is not, then the relative gap between the SVF-ordering of the patients and the optimal ordering is very small (always less than 1\%).

\begin{table}[ht!]
	\centering	
 \caption {{\small Comparison of the SVF-ordering and optimal ordering of the patients, per surgical specialty\label{tab:svf-opt}}}
 \smallskip
\begin{tabular}{crr}
\toprule
 Specialty & Frequency SVF-Optimal & Average Relative Gap\\
 \midrule
\texttt{CARD}  & $75\%$ & $ 1.8  \cdot 10^{-3}$ \\
\texttt{GASTRO}  & $84\%$ & $ 1.1  \cdot 10^{-3}$ \\
\texttt{GYN}  & $50\%$ & $ 4.2  \cdot 10^{-3}$ \\
\texttt{MED}  & $56\%$ & $ 2.2  \cdot 10^{-3}$ \\
\texttt{ORTH}  & $88\%$ & $ 1.0 \cdot 10^{-3}$  \\
\texttt{URO}  & $44\%$ & $ 3.5\cdot 10^{-3}$  \\
\bottomrule
\end{tabular}
\end{table}

\paragraph{Scalability of our approach.}
We performed additional experiments to measure the performance of our approach (SMB2SS) for larger instances. In the input, we multiply the number of patients, the emergency rate and the number of blocks available each day and in each specialty by a scaling factor $ \rho \in \{1, 2, 3, 4, 5 \}$. 
As Table~\ref{tab:scale} shows, our approach can solve very large instances in less than 90 s. For these experiments, we set the relative-gap-limit parameter of the MIP-solver to a tolerance of $0.5\%$.

\begin{table}[ht!]
\centering
\caption{\small Runtime (in seconds)
to solve the MIP~\eqref{2MIPs} with parameter $n_e=8\rho$ to compute the assignment, and the LP~\eqref{LPstage2} for all blocks $b\in B$ to set the tentative times,
depending on the size of scaled instances \label{tab:scale}}
\medskip
{\small
\begin{tabular}{cccccc}
\toprule
 $ \rho$ & $1$ & $2$ & $3$ & $4$ & $5$\\
 $ n$ & $200$ & $400$ & $600$ & $800$ & $1000$\\
 $ \lambda$ & $3$ & $6$ & $9$ & $12$ & $15$\\
 $ |B|$ & $32$ & $64$ & $96$ & $128$ & $160$\\
 \midrule
Total-Time  & $4.4749$ & $ 13.3102$ & $ 18.7838$ & $ 74.4704$\ & $ 84.3221$\\
\midrule
Solving the MIP~\eqref{2MIPs}  & $0.9383$ & $ 5.9060 $ & $ 7.5897$ & $ 59.2352$ & $ 65.2486$\\
\midrule
Solving the second stage LPs~\eqref{LPstage2}  & $3.5366$ & $7.4041 $ & $ 11.1940$ & $ 15.2351$ & $ 19.0734 $\\
\bottomrule
\end{tabular}
}
\end{table}

\section{Conclusion and outlook} \label{sec:conclusion}

We have proposed a practical approach to solve the
\EESPlong (\EESP). It consists of a 
surrogate model-based two-stage stochastic programming approach to solve the offline phase (\APP) and a simple greedy policy
to handle emergencies in the online phase of the problem (\OSP).
Our approach yields schedules
with total expected costs up to \changed{30\%} cheaper than with standard heuristics. Moreover, it runs very fast (less than 1 minute)
and is therefore well suited for an implentation in 
a real-world situation. In contrast, approaches relying 
on optimizing over a training set of scenarios (sampling average approximation) fail to obtain good solutions in a reasonable
amount of time, even for relatively small instances.

\changed{
Another possible direction of research for future work
is to take into account the availability of
downstream resources, i.e., beds in 
post-anesthesia care unit (PACU), intensive care unit (ICU), and ward beds. It is known that the patient assignment problem
becomes much harder when we consider
these resources~\cite{naderi2021increased}. Since the approach
presented in this paper produces a MIP of relatively small size
and which is easy to solve, 
we believe it has the potential to be used as part of a more complex model incorporating these resources.
}

The concept of a solution for the \APP considered in this paper 
is to assign patient to blocks in advance; apart from 
rare rescheduling decisions, with our policy
the patients are operated in the block
corresponding to the computed assignment.
In the field of stochastic parallel machine
scheduling, this is known as a fixed-assignment
policy. It is well known that adaptive policies,
in which the assignment of jobs to machines (in the context of \EESP, the assignment of patients to blocks)
\emph{depends on the realization of random processing times},
can give much better costs. 
In particular, it happens that two blocks of the same specialty run in parallel, e.g.\ \texttt{CARD} or 
\texttt{URO} on Monday, see~Table~\ref{tab:mss}; We believe that in situations where one block causes overtime while the other one is under-utilized, migrating one job from the most-loaded block to the least-loaded one could {significantly} improve the schedule. \changed{In practice however, same-day patient migrations tend
to be avoided as it is a source of stress in the operation theater.
The use of semi-adaptive policies could be an interesting
solution concept to explore in this domain;
see e.g.~\citet{DBLP:conf/esa/SagnolW21}.}

\begin{small}
\setlength{\bibsep}{0.7ex}

\bibliographystyle{apalike}
\bibliography{mopta.bib}{}

\end{small}

\begin{landscape}
	\begin{table}
	\caption{\small Comparison of four assignment algorithms
	for the cost structure $cs_2\ (\co=1, \cw=1, \ci=0)$
	in terms of runtime (CPU time in seconds), MIP gap after a time limit of 5 minutes for each specialty, MIP cost and simulated cost on the validation set.
	\label{benders_cs2}}
	\medskip
	\centering
	\footnotesize
	\medskip

\begin{tabular}{ccc|rrrr|rrrr|rrrr|rrrr}
\\
\multicolumn{19}{c}{$n=70$}\\
\hline\\[-1em]
\multirow{2}{*}{specialty}   &  \multirow{2}{*}{\#ORs}   &  \multirow{2}{*}{\#patients}
& \multicolumn{4}{c|}{SAA ($K=50$)} & \multicolumn{4}{c|}{SAA ($K=450$)} & \multicolumn{4}{c|}{Benders ($K=450$)} & \multicolumn{4}{c}{SMB2SS} \\
& & & CPU & gap & cost & sim.
& CPU & gap & cost & sim.
& CPU & gap & cost & sim.
& CPU & gap & cost & sim.\\ \hline
\texttt{CARD} & 5 & 10 
 &  66.4  & 0.29\% & 42.9  & \bf{44.6} 
 &  298.7  & 11.85\% & 44.9  & 45.6 
 &  299.9  & 14.42\% & 47.3  & 48.8 
 &  0.2  & 0.00\% & 46.0  & 51.0 \\ 
\texttt{GASTRO} & 6 & 13 
 &  190.7  & 9.86\% & 50.5  & \bf{55.9} 
 &  309.7  & 38.27\% & 76.4  & 76.8 
 &  301.8  & 58.21\% & 101.6  & 101.4 
 &  0.8  & 0.00\% & 59.7  & 62.0 \\ 
\texttt{GYN} & 8 & 19 
 &  300.5  & 9.76\% & 83.1  & \bf{86.7} 
 &  318.3  & 49.61\% & 155.5  & 150.2 
 &  302.3  & 73.52\% & 264.5  & 256.1 
 &  2.7  & 0.00\% & 90.7  & 92.9 \\ 
\texttt{MED} & 1 & 3 
 &  0.1  & 0.00\% & 20.2  & 22.3 
 &  0.6  & 0.00\% & 22.2  & 22.3 
 &  0.2  & 0.00\% & 22.2  & \bf{22.3} 
 &  0.0  & 0.00\% & 24.4  & 22.6 \\ 
\texttt{ORTH} & 6 & 12 
 &  209.9  & 2.53\% & 62.5  & \bf{64.5} 
 &  341.7  & 28.93\% & 81.2  & 83.8 
 &  301.9  & 47.45\% & 106.1  & 110.8 
 &  0.9  & 0.00\% & 62.3  & 67.9 \\ 
\texttt{URO} & 6 & 13 
 &  156.0  & 0.13\% & 35.7  & \bf{38.3} 
 &  785.5  & 15.84\% & 40.2  & 40.5 
 &  301.3  & 62.78\% & 92.0  & 98.0 
 &  0.3  & 0.00\% & 40.7  & 39.7 \\ 
\hline
\texttt{TOTAL} & 32 & 70 
 &  925.8  & 5.87\% & 294.9  & \bf{312.3} 
 &  2056.9  & 36.87\% & 420.3  & 419.2 
 &  1510.1  & 58.97\% & 633.6  & 637.3 
 &  7.3  & 0.00\% & 323.7  & 336.3 \\ 
\hline
\\
\multicolumn{19}{c}{$n=100$}\\
\hline
\texttt{CARD} & 5 & 14 
 &  195.1  & 2.49\% & 56.8  & \bf{59.9} 
 &  341.1  & 35.70\% & 76.0  & 77.1 
 &  301.8  & 68.50\% & 135.9  & 138.2 
 &  0.6  & 0.00\% & 59.1  & 68.4 \\ 
\texttt{GASTRO} & 6 & 18 
 &  301.8  & 50.65\% & 258.4  & \bf{275.0} 
 &  312.4  & N/A & N/A  & N/A 
 &  302.9  & 70.24\% & 326.8  & 333.4 
 &  3.1  & 0.00\% & 266.7  & 275.1 \\ 
\texttt{GYN} & 8 & 28 
 &  303.1  & 31.29\% & 170.9  & 183.6 
 &  323.4  & N/A & N/A  & N/A 
 &  304.1  & 75.61\% & 446.4  & 459.2 
 &  3.8  & 0.00\% & 159.3  & \bf{179.3} \\ 
\texttt{MED} & 1 & 5 
 &  0.2  & 0.00\% & 94.6  & 96.8 
 &  2.0  & 0.00\% & 97.8  & 95.9 
 &  1.0  & 0.04\% & 98.2  & \bf{95.8} 
 &  0.0  & 0.00\% & 100.1  & 99.5 \\ 
\texttt{ORTH} & 6 & 17 
 &  301.6  & 36.38\% & 150.1  & 166.6 
 &  1257.9  & N/A & N/A  & N/A 
 &  303.0  & 63.52\% & 234.4  & 237.1 
 &  2.1  & 0.00\% & 163.0  & \bf{165.9} \\ 
\texttt{URO} & 6 & 18 
 &  301.5  & 7.85\% & 73.2  & \bf{75.0} 
 &  313.4  & 32.27\% & 101.4  & 101.3 
 &  301.8  & 67.64\% & 199.4  & 199.0 
 &  0.4  & 0.00\% & 77.4  & 82.6 \\ 
\hline
\texttt{TOTAL} & 32 & 100 
 &  1406.7  & 33.26\% & 804.0  & \bf{857.0} 
 &  2558.4  & N/A & N/A  & N/A 
 &  1518.5  & 66.06\% & 1441.1  & 1462.6 
 &  13.7  & 0.00\% & 825.6  & 870.8 \\ 
\hline
\\
\multicolumn{19}{c}{$n=140$}\\
\hline
\texttt{CARD} & 5 & 20 
 &  301.7  & 51.75\% & 180.0  & \bf{186.1} 
 &  312.3  & 69.25\% & 279.6  & 278.1 
 &  302.7  & 74.55\% & 253.4  & 260.4 
 &  2.0  & 0.00\% & 180.9  & 193.5 \\ 
\texttt{GASTRO} & 6 & 25 
 &  302.2  & 49.74\% & 901.1  & 914.1 
 &  1262.8  & N/A & N/A  & N/A 
 &  304.7  & 55.29\% & 850.7  & 843.3 
 &  7.4  & 0.00\% & 816.9  & \bf{839.9} \\ 
\texttt{GYN} & 8 & 39 
 &  341.1  & 65.65\% & 487.7  & 533.8 
 &  329.8  & N/A & N/A  & N/A 
 &  306.5  & 78.54\% & 692.2  & 720.3 
 &  17.5  & 0.44\% & 352.1  & \bf{413.5} \\ 
\texttt{MED} & 1 & 7 
 &  0.5  & 0.00\% & 159.5  & 171.6 
 &  15.5  & 0.00\% & 172.9  & 169.9 
 &  5.5  & 0.16\% & 171.6  & \bf{169.9} 
 &  0.0  & 0.00\% & 171.1  & 177.2 \\ 
\texttt{ORTH} & 6 & 24 
 &  302.1  & 44.81\% & 803.5  & 806.9 
 &  351.1  & N/A & N/A  & N/A 
 &  340.2  & 45.74\% & 728.4  & 735.5 
 &  19.9  & 0.31\% & 688.7  & \bf{715.8} \\ 
\texttt{URO} & 6 & 25 
 &  302.2  & 28.62\% & 101.6  & \bf{101.7} 
 &  353.1  & 66.11\% & 226.7  & 218.6 
 &  302.9  & 77.41\% & 295.2  & 287.4 
 &  2.1  & 0.00\% & 94.3  & 110.1 \\ 
\hline
\texttt{TOTAL} & 32 & 140 
 &  1554.3  & 47.63\% & 2633.4  & 2714.3 
 &  2673.0  & N/A & N/A  & N/A 
 &  1568.1  & 58.41\% & 2991.5  & 3016.7 
 &  53.9  & 0.17\% & 2303.9  & \bf{2450.0} \\ 
\hline
\\
\multicolumn{19}{c}{$n=200$}\\
\hline
\texttt{CARD} & 5 & 28 
 &  302.0  & 42.81\% & 718.9  & 731.5 
 &  349.3  & N/A & N/A  & N/A 
 &  304.4  & 54.38\% & 736.4  & 737.4 
 &  4.8  & 0.00\% & 694.1  & \bf{707.2} \\ 
\texttt{GASTRO} & 6 & 36 
 &  303.0  & 30.59\% & 1815.6  & 1859.4 
 &  320.7  & N/A & N/A  & N/A 
 &  305.8  & 28.24\% & 1700.1  & 1717.2 
 &  8.8  & 0.01\% & 1673.2  & \bf{1697.3} \\ 
\texttt{GYN} & 8 & 56 
 &  305.3  & N/A & N/A  & N/A 
 &  816.0  & N/A & N/A  & N/A 
 &  343.1  & 60.39\% & 1427.0  & 1440.2 
 &  20.2  & 0.15\% & 1249.9  & \bf{1315.2} \\ 
\texttt{MED} & 1 & 10 
 &  2.4  & 0.00\% & 330.8  & 337.6 
 &  62.8  & 0.00\% & 335.8  & 333.9 
 &  69.9  & 0.43\% & 334.7  & \bf{333.9} 
 &  0.0  & 0.00\% & 335.1  & 346.0 \\ 
\texttt{ORTH} & 6 & 34 
 &  302.7  & 29.82\% & 1804.9  & 1808.9 
 &  319.1  & N/A & N/A  & N/A 
 &  305.6  & 21.76\% & 1556.6  & 1575.2 
 &  19.4  & 0.09\% & 1511.7  & \bf{1548.8} \\ 
\texttt{URO} & 6 & 36 
 &  777.2  & 54.67\% & 431.3  & 470.6 
 &  791.9  & N/A & N/A  & N/A 
 &  305.2  & 69.83\% & 543.9  & 565.4 
 &  8.7  & 0.01\% & 390.9  & \bf{431.4} \\ 
\hline
\texttt{TOTAL} & 32 & 200 
 &  2001.0  & N/A & N/A  & N/A 
 &  2758.0  & N/A & N/A  & N/A 
 &  1640.9  & 38.67\% & 6298.6  & 6369.4 
 &  67.4  & 0.06\% & 5854.9  & \bf{6046.0} \\ 
\hline
\end{tabular}
	\end{table}
	\end{landscape}     
	
	\begin{landscape}
	\begin{table}
	\caption{\small Comparison of four assignment algorithms
	for the cost structure $cs_3\ (\co=1, \cw=2, \ci=2)$
	in terms of runtime (CPU in seconds), MIP gap after a time limit of 5 minutes for each specialty, MIP cost and simulated cost on the validation set.
	\label{benders_cs3}}
	\medskip
	\centering
	\footnotesize
	\medskip

\begin{tabular}{ccc|rrrr|rrrr|rrrr|rrrr}
\\
\multicolumn{19}{c}{$n=70$}\\
\hline\\[-1em]
\multirow{2}{*}{specialty}   &  \multirow{2}{*}{\#ORs}   &  \multirow{2}{*}{\#patients}
& \multicolumn{4}{c|}{SAA ($K=50$)} & \multicolumn{4}{c|}{SAA ($K=450$)} & \multicolumn{4}{c|}{Benders ($K=450$)} & \multicolumn{4}{c}{SMB2SS} \\
& & & CPU & gap & cost & sim.
& CPU & gap & cost & sim.
& CPU & gap & cost & sim.
& CPU & gap & cost & sim.\\ \hline
\texttt{CARD} & 5 & 10 
 &  300.8  & 45.51\% & 186.8  & \bf{199.2} 
 &  306.1  & 78.32\% & 209.0  & 210.9 
 &  301.4  & 81.48\% & 276.1  & 281.7 
 &  0.2  & 0.00\% & 273.4  & 215.2 \\ 
\texttt{GASTRO} & 6 & 13 
 &  301.3  & 86.84\% & 321.8  & \bf{340.5} 
 &  309.6  & 89.08\% & 380.6  & 370.5 
 &  302.2  & 91.66\% & 459.0  & 451.6 
 &  0.8  & 0.00\% & 409.0  & 345.7 \\ 
\texttt{GYN} & 8 & 19 
 &  334.3  & 81.19\% & 408.2  & \bf{434.5} 
 &  789.2  & N/A & N/A  & N/A 
 &  302.5  & 90.48\% & 719.4  & 712.0 
 &  1.5  & 0.00\% & 559.2  & 447.7 \\ 
\texttt{MED} & 1 & 3 
 &  0.1  & 0.00\% & 79.4  & 81.2 
 &  1.2  & 0.00\% & 81.2  & \bf{80.3} 
 &  0.3  & -0.12\% & 81.9  & 80.7 
 &  0.0  & 0.00\% & 111.5  & 81.8 \\ 
\texttt{ORTH} & 6 & 12 
 &  301.2  & 77.88\% & 266.1  & \bf{277.9} 
 &  817.3  & 82.02\% & 315.3  & 320.7 
 &  302.1  & 84.93\% & 361.0  & 370.3 
 &  0.9  & 0.00\% & 317.8  & 282.1 \\ 
\texttt{URO} & 6 & 13 
 &  301.2  & 80.03\% & 195.7  & \bf{206.1} 
 &  309.8  & 85.32\% & 255.1  & 256.5 
 &  301.8  & 90.36\% & 340.6  & 350.0 
 &  0.2  & 0.00\% & 281.7  & 225.3 \\ 
\hline
\texttt{TOTAL} & 32 & 70 
 &  1541.1  & 72.93\% & 1458.0  & \bf{1539.3} 
 &  2540.0  & N/A & N/A  & N/A 
 &  1513.0  & 85.50\% & 2237.9  & 2246.2 
 &  6.0  & 0.00\% & 1952.6  & 1597.6 \\ 
\hline
\\
\multicolumn{19}{c}{$n=100$}\\
\hline
\texttt{CARD} & 5 & 14 
 &  301.2  & 84.14\% & 340.6  & \bf{359.9} 
 &  309.2  & 86.46\% & 379.1  & 383.9 
 &  302.0  & 90.13\% & 430.8  & 438.2 
 &  0.5  & 0.00\% & 409.4  & 365.6 \\ 
\texttt{GASTRO} & 6 & 18 
 &  335.2  & 84.43\% & 733.1  & 766.6 
 &  787.8  & N/A & N/A  & N/A 
 &  303.3  & 88.28\% & 792.0  & 808.0 
 &  7.2  & 0.02\% & 776.8  & \bf{758.3} \\ 
\texttt{GYN} & 8 & 28 
 &  303.1  & 85.72\% & 814.1  & 859.3 
 &  795.8  & N/A & N/A  & N/A 
 &  303.9  & 90.17\% & 1097.9  & 1124.6 
 &  3.9  & 0.00\% & 887.9  & \bf{829.6} \\ 
\texttt{MED} & 1 & 5 
 &  0.3  & 0.00\% & 203.1  & 216.6 
 &  43.5  & 0.00\% & 211.7  & \bf{211.5} 
 &  1.3  & -0.07\% & 213.6  & 213.7 
 &  0.0  & 0.00\% & 267.5  & 232.5 \\ 
\texttt{ORTH} & 6 & 17 
 &  301.5  & 83.19\% & 529.3  & 563.5 
 &  311.4  & N/A & N/A  & N/A 
 &  303.1  & 85.97\% & 603.5  & 596.9 
 &  2.9  & 0.00\% & 540.6  & \bf{547.9} \\ 
\texttt{URO} & 6 & 18 
 &  301.7  & 80.64\% & 365.3  & \bf{392.3} 
 &  312.1  & N/A & N/A  & N/A 
 &  302.1  & 87.10\% & 504.3  & 504.2 
 &  0.5  & 0.00\% & 454.4  & 407.2 \\ 
\hline
\texttt{TOTAL} & 32 & 100 
 &  1546.3  & 78.26\% & 2985.4  & 3158.3 
 &  2577.7  & N/A & N/A  & N/A 
 &  1519.3  & 83.44\% & 3642.2  & 3685.7 
 &  18.5  & 0.01\% & 3336.6  & \bf{3141.1} \\ 
\hline
\\
\multicolumn{19}{c}{$n=140$}\\
\hline
\texttt{CARD} & 5 & 20 
 &  333.4  & 86.98\% & 633.0  & 668.7 
 &  785.2  & 90.35\% & 767.5  & 769.6 
 &  303.0  & 90.81\% & 693.8  & 699.3 
 &  5.6  & 0.01\% & 661.6  & \bf{656.5} \\ 
\texttt{GASTRO} & 6 & 25 
 &  302.2  & 65.66\% & 1424.4  & \bf{1464.5} 
 &  787.7  & N/A & N/A  & N/A 
 &  304.8  & 70.87\% & 1505.3  & 1496.4 
 &  16.1  & 0.03\% & 1447.4  & 1513.9 \\ 
\texttt{GYN} & 8 & 39 
 &  304.4  & 88.12\% & 1374.6  & 1465.1 
 &  329.5  & N/A & N/A  & N/A 
 &  342.5  & 90.77\% & 1603.1  & 1634.2 
 &  16.9  & 0.08\% & 1378.9  & \bf{1376.0} \\ 
\texttt{MED} & 1 & 7 
 &  1.4  & 0.00\% & 304.5  & 317.4 
 &  22.8  & 0.00\% & 310.0  & \bf{310.6} 
 &  7.3  & -0.03\% & 312.0  & 311.8 
 &  0.0  & 0.00\% & 369.3  & 358.3 \\ 
\texttt{ORTH} & 6 & 24 
 &  302.0  & 61.31\% & 1245.8  & 1250.3 
 &  345.9  & N/A & N/A  & N/A 
 &  304.5  & 63.54\% & 1208.8  & 1221.4 
 &  19.7  & 0.13\% & 1166.8  & \bf{1208.6} \\ 
\texttt{URO} & 6 & 25 
 &  302.1  & 88.42\% & 624.8  & 624.8 
 &  314.7  & N/A & N/A  & N/A 
 &  303.3  & 91.26\% & 762.8  & 752.6 
 &  2.4  & 0.00\% & 634.0  & \bf{612.5} \\ 
\hline
\texttt{TOTAL} & 32 & 140 
 &  1549.8  & 71.41\% & 5607.2  & 5790.9 
 &  2642.7  & N/A & N/A  & N/A 
 &  1570.6  & 75.57\% & 6085.8  & 6115.5 
 &  65.5  & 0.06\% & 5658.0  & \bf{5725.7} \\ 
\hline
\\
\multicolumn{19}{c}{$n=200$}\\
\hline
\texttt{CARD} & 5 & 28 
 &  302.0  & 64.18\% & 1269.5  & 1311.2 
 &  313.9  & N/A & N/A  & N/A 
 &  304.6  & 71.18\% & 1336.4  & 1344.3 
 &  12.5  & 0.02\% & 1272.9  & \bf{1297.2} \\ 
\texttt{GASTRO} & 6 & 36 
 &  303.1  & 41.14\% & 2474.9  & 2542.0 
 &  793.0  & N/A & N/A  & N/A 
 &  340.0  & 43.22\% & 2528.2  & 2558.5 
 &  18.5  & 0.03\% & 2472.4  & \bf{2540.6} \\ 
\texttt{GYN} & 8 & 56 
 &  339.0  & N/A & N/A  & N/A 
 &  341.7  & N/A & N/A  & N/A 
 &  311.9  & 76.32\% & 2702.7  & 2742.5 
 &  20.1  & 0.07\% & 2472.9  & \bf{2551.2} \\ 
\texttt{MED} & 1 & 10 
 &  3.5  & 0.00\% & 508.5  & 516.4 
 &  570.8  & N/A & N/A  & N/A 
 &  63.2  & 0.29\% & 514.0  & \bf{514.3} 
 &  0.0  & 0.00\% & 572.9  & 588.9 \\ 
\texttt{ORTH} & 6 & 34 
 &  778.4  & 37.17\% & 2250.2  & 2278.4 
 &  319.5  & N/A & N/A  & N/A 
 &  305.8  & 36.70\% & 2171.3  & 2203.0 
 &  19.8  & 0.07\% & 2105.0  & \bf{2177.0} \\ 
\texttt{URO} & 6 & 36 
 &  779.8  & 83.44\% & 1102.4  & 1164.1 
 &  831.9  & N/A & N/A  & N/A 
 &  305.9  & 86.42\% & 1209.7  & 1203.1 
 &  50.2  & 0.02\% & 1089.1  & \bf{1108.2} \\ 
\hline
\texttt{TOTAL} & 32 & 200 
 &  2511.3  & N/A & N/A  & N/A 
 &  3752.7  & N/A & N/A  & N/A 
 &  1638.3  & 56.51\% & 10462.4  & 10565.7 
 &  126.8  & 0.04\% & 9985.1  & \bf{10263.1} \\ 
\hline
\end{tabular}
	\end{table}
	\end{landscape} 

\newpage
\appendix 

\spacing{1}

\begin{center}
	{\Large \textbf{Appendix}}
\end{center}
\bigskip

\section{Benders decomposition approach}\label{sec:benders}

\changed{
As for the case of the SAA approach (see Section~\ref{sec:other_approaches}), the Benders decomposition technique
typically takes a long time, even on small instances.
We thus restrict to the case without emergencies, so
that there is no coupling between the different surgical specialties.
}

\changed{
The problem to be solved is Problem~\eqref{P1},
but we restrict our attention
to the blocks $b\in B_s$ of specialty $s$ and we assume that
there are no emergencies (so we can set $z=0$). This problem can be
rewritten as:
\begin{align}\label{benders-Master0}
	\min_{\bm{x}\in X}&\quad \sum_{i\in I_s} \sum_{b \in B_s'} x_{ib}\, c_{ib} + \sum_{b \in B_s} \hat{\zeta}_{b}\\ 
	\text{s.t.} 
	&\quad \sum_{b\in B_s'} x_{ib} =1  && \forall \; i \in I_s\nonumber \\
	&\quad  x_{ib} \in \{ 0, 1\} && \forall i \in I_s, \; \forall \; b \in B_s',\nonumber \\
	&\quad \hat{\zeta}_{b} \geq \phi_b(\bm{x},\bm{0}), && \forall b \in B_s \nonumber
\end{align}
}

\changed{
It will be useful to work with second-stage costs which
are linear with respect to the first stage variables $\vec{x}$, without constant term.
We thus define 
$\psi_b(\bm{x}):=\phi_b(\bm{x},\vec{0})
+\frac{{\ci} }{K} \sum_{k=1}^{K} p^{\omega_{k}}_i$,
and make the change of variable 
${\zeta}_b=\hat{\zeta}_b+\frac{{\ci} }{K} \sum_{k=1}^{K} p^{\omega_{k}}_i$
to rewrite~\eqref{benders-Master0} as:}
\begin{align}\label{benders-Master}
\min_{\bm{x}\in X}&\quad \sum_{i\in I_s} \sum_{b \in B_s'} x_{ib}\, \Big( c_{ib} - \frac{{\ci} }{K} \sum_{k=1}^{K} p^{\omega_{k}}_i \Big)  + \sum_{b \in B_s} \zeta_{b}\\ 
\text{s.t.} 
&\quad \sum_{b\in B_s'} x_{ib} =1  && \forall \; i \in I_s\nonumber \\
&\quad  x_{ib} \in \{ 0, 1\} && \forall i \in I_s, \; \forall \; b \in B_s',\nonumber \\
&\quad \zeta_{b} \geq \psi_b(\bm{x}), && \forall b \in B_s. \nonumber
\end{align}
Note that the term $-\frac{\ci}{K}\sum_{k=1}^K \sum_{i\in I_s} x_{ib}\, p_i^{\omega_k}$
has been put in the first-stage part of the objective of~\eqref{benders-Master}. Thus, we have
\begin{align}
\psi_b(\bm{x}) = \min&\quad \frac{1}{ K} \label{LPstage2_benders}
\rlap{$\displaystyle{\sum_{k=1}^K} \Big( \sum_{i=1}^{n_b} {\cw} ( s_{i,k} - t_i)  + {\ci} \, L_k + {\co}\, O_k \Big)$,}\\
&\quad s_{i,k} \geq t_i,  && i=1, \dotsc,  n_b,\quad k=1,\ldots,K \nonumber\\
&\quad s_{i,k} \geq s_{i-1,k} + p^{\omega_{k}}_{i-1}, &&i=2, \dotsc, n_b,\quad k=1,\ldots,K \nonumber\\
& \quad L_k = s_{n_b,k} + p^{\omega_{k}}_{n_b}, &&k=1,\ldots,K\nonumber\\
&\quad O_k \geq L_k - T, &&k=1,\ldots,K\nonumber\\
&\quad s_{i,k},\ t_i,\ O_k \geq 0, &&i=1, ..., n_b,\quad k=1,\ldots,K, \nonumber
\end{align}
where the patients assigned to block $b$ for the assignment vector $\bm{x}$ have been numbered as
$1,\ldots,n_b$, by order of nondecreasing variance.

The idea of the Benders decomposition is to replace the constraints $\zeta_b\geq \psi_b(\bm{x})$ in Problem~\eqref{benders-Master} by a collection of linear inequalities, which can be obtained by solving the dual problem of~\eqref{LPstage2_benders} for a candidate assignment $\bm{x}$.
The dual reads as follows:
\begin{align}
\psi_b(\bm{x}) = \max\ &\quad  \label{dual_LP}
\sum_{k=1}^K \sum_{i=1}^{n_b}  \lambda_{ik}   p^{\omega_{k}}_i -  \sum_{k=1}^K  \alpha_k  T,\\
s.t.\ &\quad \sum_{k=1}^K  \mu_{ik} \geq  c^{\mathbf{w}},  && i=1, \dotsc,  n_b,\quad k=1,\ldots,K \nonumber\\
&\quad \frac{c^{\mathbf{o}}}{K}\geq \alpha_k, &&k=1,\ldots,K\nonumber\\
&\quad \frac{c^{\mathbf{i}}}{K} + \alpha_k \geq \lambda_{n_b-1,k}, &&k=1,\ldots,K\nonumber\\
&\quad \lambda_{1k} + \frac{c^{\mathbf{w}}}{K} \geq \mu_{1k} , &&k=1,\ldots,K \nonumber\\
&\quad \lambda_{ik} - \lambda_{i-1,k} + \frac{c^{\mathbf{w}}}{K} \geq \mu_{ik} , &&i=2, \dotsc, n_b,\quad k=1,\ldots,K \nonumber\\
&\quad \lambda_{ik},\ \mu_{ik},\ \alpha_k \geq 0, &&i=1, ..., n_b,\quad k=1,\ldots,K. \nonumber
 \end{align}
Therefore, given a candidate solution $(\bm{x},\bm{\zeta})$ of~\eqref{benders-Master}, we can check whether an inequality of the form\linebreak $\zeta_b \geq \psi_b(\bm{x})$ is violated by solving the above dual problem. Should the inequality be violated, we can add a cutting plane of the form 
\[\zeta_{b} \geq \psi_b(\vec{x})= \sum_{k=1}^K \sum_{i\in I_s}  x_{ib}\, \bar{\lambda}_{ik}\,   p^{\omega_{k}}_i -  \bar{\alpha}\,  T,\]
where $\bar{\lambda}$ and $\bar{\alpha}$ are obtained from the optimal dual variables $(\bm{\lambda}^*,\bm{\mu}^*,\bm{\alpha}^*)$ of
Problem~\eqref{dual_LP}:
\begin{equation}\label{eq:lambdabar}
\bar{\alpha} = \sum_{k=1}^K \alpha_k^*\quad \text{and}\quad
\bar{\lambda}_{ik}=\left\{
\begin{array}{ll}
\lambda_{jk}^* &  \text{if $i$ is the $j$th patient in block $b$  (in SVF order)}\\
0              & \text{otherwise.}
\end{array}
\right.
\end{equation}
Given a collection of sets $\mathcal{C}_b=\{(\bar{\vec{\lambda}_1},\bar{\alpha}_1),\ldots,(\bar{\vec{\lambda}_m},\bar{\alpha}_m)\}$
of cutting planes for the block $b$, $\forall b\in B_s$, we therefore consider the following restricted master problem:
\begin{align}\label{RMP}
\min_{\bm{x}\in X}&\quad \sum_{i\in I_s} \sum_{b \in B_s'} x_{ib}\, \Big( c_{ib} - \frac{c^{\mathbf{i}} }{K} \sum_{k=1}^{K} p^{\omega_{k}}_i \Big)  + \sum_{b \in B_s} \zeta_{b}\\ 
\text{s.t.} 
&\quad \sum_{b\in B_s'} x_{ib} =1  && \forall \; i \in I_s\nonumber \\
&\quad  x_{ib} \in \{ 0, 1\} && \forall i \in I_s, \; \forall \; b \in B_s',\nonumber \\
&\quad \zeta_b \geq \sum_{k=1}^K \sum_{i\in I_s}  x_{ib}\, \bar{\lambda}_{ik}\,   p^{\omega_{k}}_i -  \bar{\alpha}\,  T, && \forall b \in B_s, \forall (\vec{\bar{\lambda}},\bar{\alpha})\in\mathcal{C}_b \nonumber
\end{align}
It remains to discuss how the set $\mathcal{C}_b$ 
can be initialized.
We know that 
\[\psi_b(\vec{x}) = \phi_b(\vec{x}) + \frac{\ci}{K} \displaystyle{\sum_{k=1}^K }\sum_{i\in I_s} x_{ib}\, p_i^{\omega_k} 
\geq
\frac{\co}{K} \displaystyle{\sum_{k=1}^K }
\max \Big( \sum_{i\in I_s}x_{ib}\, p_i^{\omega_k} - T,0 \Big)
+ \frac{\ci}{K} \displaystyle{\sum_{k=1}^K }\sum_{i\in I_s} x_{ib}\, p_i^{\omega_k},
\]
where the inequality comes from the fact that the idling costs
and waiting costs are nonnegative. Thus, by
\changed{separately considering} the two terms in the $\max()$, 
we obtain the linear bounds
\[
 \psi_b(\vec{x}) \geq 
 \sum_{k=1}^K \sum_{i\in I_s} x_{ib} p_i^{\omega_k} \cdot \frac{\ci}{K}
 \qquad\text{and}\qquad
  \psi_b(\vec{x}) \geq 
 \sum_{k=1}^K \sum_{i\in I_s} x_{ib} p_i^{\omega_k} \cdot \frac{\co+\ci}{K} - \co \cdot T.
\]
We can thus  initialize each $\mathcal{C}_b$ by setting $\mathcal{C}_b=\{(\vec{\bar{\lambda}}_0,0),(\vec{\bar{\lambda}}_1,\co)\}$,
where $\vec{\bar{\lambda}}_0$ and
$\vec{\bar{\lambda}}_1$ are the 
constant vectors given by:  $(\vec{\bar{\lambda}}_0)_{ik} = \frac{\ci}{K}$
and 
$(\vec{\bar{\lambda}}_1)_{ik} = \frac{\co+\ci}{K},\ \forall i\in I_s, \forall k\in\{1,\ldots,K\}$.
Then, we solve iteratively the restricted master problem~\eqref{RMP} to obtain a candidate solution $(\vec{x},\vec{\zeta})$ and we check its feasibility by solving the subproblems~\eqref{dual_LP} for each block $b\in B$, adding cutting planes whenever the constraint $\zeta_b\geq \psi_b(\vec{x})$ is violated.
Our approach is summarized in Algorithm~\ref{alg:benders}.

\begin{algorithm}[t]
Initialize the set of cuts $\mathcal{C}_b=\{(\vec{\bar{\lambda}}_0,0),(\vec{\bar{\lambda}}_1,\co )\}$, $\forall b\in B_s$ \\
$UB \gets \infty$\\
$LB \gets 0$\\
 \While {$\frac{UB-LB}{UB} > \epsilon$}
 {
compute an optimal solution $(\bm{x}, \vec{\zeta})$ to~\eqref{RMP}\\
$LB \gets \displaystyle{\sum_{i\in I_s}} \sum_{b \in B_s'} x_{ib}\, \Big( c_{ib} - \frac{c^{\mathbf{i}} }{K} \sum_{k=1}^{K} p^{\omega_{k}}_i \Big)  + \sum_{b \in B_s} \zeta_{b}$\\
 \For{$ b\in B_s$}
 {
 	solve the LP~\eqref{dual_LP} to get $\psi_b(\vec{x})$, $\bm{\lambda}^*$ and $\bm{\alpha}^*$\\
 \If{$\psi_b(\vec{x}) > \zeta_b $}
 {
  	compute $(\vec{\bar{\lambda}},\bar{\alpha})$ using~\eqref{eq:lambdabar}\\
 	add the cutting plane $(\vec{\bar{\lambda}},\bar{\alpha})$ to $\mathcal{C}_b$
 	}
 }
 $UB \gets \min\left(UB,\ \displaystyle{ \sum_{i\in I_s}} \sum_{b \in B_s'} x_{ib}\, \Big( c_{ib} - \frac{c^{\mathbf{i}} }{K} \sum_{k=1}^{K} p^{\omega_{k}}_i \Big)  + \sum_{b \in B_s} \psi_{b}(\vec{x})\right)$\\
  }
 \caption{\small Benders Decomposition Approach for Solving Problem~\eqref{BMIPs} \label{alg:benders}}
\end{algorithm}

\section{Demonstrator} \label{sec:demonstrator}
As mentioned in the introduction, a demonstrator
has been put online:  
\url{https://wsgi.math.tu-berlin.de/esp_demonstrator/}.
We shortly explain its features in this section.
\changed{\emph{Note:} It is sometimes necessary to empty the cache of the browser and to refresh the page (with Ctrl+F5 or Shift+F5) so that the content loads properly.
}

The user can select many different parameters to set up the instance to be solved and the used algorithm; see Figure~\ref{fig:input}.
For example, there are dropdown menus to select the instance seed value, the number of elective surgeries \(n \in \{70,\, 100\,, 140,\, 200\}\), the cost parameters and the emergency rate \(\lambda \in \{0,\,1,\,2,\,3,\,4\}\).
It is also possible to select the value of the
parameters \(\alpha\) and $\Delta$ used by the online
policy; see Section~\ref{sec:onlinephase}. The default value of  $\Delta$ is set to a very high value (1000), so no patient is ever cancelled.

\bigskip
After clicking on \emph{Optimize Assignment}, 
the computed solution is displayed in the right panel
of the window; see Figure~\ref{fig:input}.
The visualization includes a Gantt chart
showing the schedule on every day of the week, and 
a histogram showing how many patients are scheduled, postponed, cancelled or rescheduled in each specialty.
Note that the output is dynamic, and depends on the 
\emph{current time} \changed{$\tau$} which the user can select
by moving a cursor on the slider below the figures.
The expected outcome is displayed, with
\begin{itemize}
 \item The true realization of the specified scenario up to time \changed{$\tau$.}
 \item Expected values conditioned to observations until~\changed{$\tau$} for the time after \changed{$\tau$}.
\end{itemize}
By moving the cursor, we can see how the prediction evolves, and when a surgery gets cancelled.
The time after \changed{$\tau$} is grayed out on the Gantt Chart.
In particular, if the cursor is placed at the beginning (Monday, 8AM) with a high value of the parameter $\Delta$ that controls rescheduling decisions, we can visualize the assignment and tentative starting times as they were planned after
solving the offline problem (\APP).

Finally, at the bottom of the app, there is a Monte-Carlo Simulation; see Figure~\ref{fig:Monte-Carlo}. Instead of simulating and visualizing a single scenario, a high number of scenarios is considered. The graph shows the distribution of the total costs, and how this cost is distributed between the different parts (waiting, idling, overtime, etc.) for each simulated scenario (ordered from left to right by total cost).

\section{Generating Probability distributions for surgery durations.} \label{sec:gen-dist}
\changed{
The coefficient of variation of patients $i\in I_s$ is defined as $CV_s=\frac{\sqrt{v_s}}{m_s}$,
which is a measure of the dispersion of
durations within specialty $s$. In this section we explain how
we generate individual parameters ($\mu_i,\sigma_i$) 
to create random instances, in a way that the coefficient
of variation of $P_i$ is smaller than the marginal
coefficient of variation $CV_{s(i)}$, but such that 
the overall coefficient of variations of all patients of
specialty $s$ remains equal to $CV_s$.}

We assume that upon observation of patient $i$, a learning algorithm \changed{(which may rely on the estimate given by a surgeon)} produces an estimation of the parameters
$\mu_i,\sigma_i$ of the random duration $P_i\sim\mathcal{L}\mathcal{N}(\mu_i,\sigma_i)$   of patient $i$, such that the coefficient of variation has been reduced (the distribution of a surgery duration conditioned to additional information on a patient $i\in I_s$ is less dispersed
that the overall distribution of durations in specialty $s$).
More precisely, we assume that the coefficient of variation is decreased by roughly
one half\footnote{\changed{The factor $\frac{1}{2}$ is arbitrary,
but seems to be consistent with the companion dataset of~\cite{sagnol2018robust}.}}, i.e., 
\[CV_i := \frac{\sqrt{\mathbb{V}[P_i]}}{\mathbb{E}[P_i]} = \frac{1}{2}\cdot \Delta_i \cdot CV_{s(i)},\]
where $\Delta_i$ is a random parameter drawn from a normal distribution centered about 1, in practice we sample $\Delta_i$ from $\mathcal{N}(1,0.15^2)$.
Since the coefficient of variation of a lognormal random variable $\mathcal{LN}(\mu_i,\,\sigma_i)$ is equal to $\sqrt{e^{\sigma_i^2}-1}$, a simple calculation yields
$$
\sigma_i :=  \log^{1/2}\left( \big(CV_i\big)^2 + 1\right)=
\sqrt{\log\left(\frac{\Delta_i^2\, v_s}{4m_s^2} + 1\right)},\quad
\forall i\in I_s.
$$
Now, we generate a random variable $M_i$
such that the distribution of $P_i$ \emph{conditioned to the event $M_i=\mu_i$}
is lognormal with parameters $\mu_i$ and $\sigma_i$, i.e., we set $P_i|M_i \sim \mathcal{L}\mathcal{N}(M_i,\sigma_i)$.
By the law of total expectation and the law of total variance, we want to set the law of $M_i$ in such a way that 
$$
m_s = \mathbb{E}[\mathbb{E}[P_i|M_i]]
\qquad\text{and}\qquad
v_s = \mathbb{E}[\text{Var}[P_i|M_i]] + \text{Var}[\mathbb{E}[P_i|M_i]].
$$
Again, using standard formula about lognormal random variables, one can check that generating
$M_i$ as $M_i\sim\mathcal{N}(\mu_i',\sigma_i')$, with
\[
\mu_i':=\log \frac{m_s}{\sqrt{1+\frac{v_s}{m_s^2}}}
\quad\text{and}\quad
\sigma_i':=\sqrt{\log\left( \frac{1+\frac{v_s}{m_s^2}}{1+\frac{\Delta_i^2}{4}\frac{v_s}{m_s^2}}\right)}
\]
satisfies the above requirements for the marginal distribution of $P_i$.

We proceed similarly for the emergency patients: We assume that a posterior distribution 
$P_e\sim\mathcal{L}\mathcal{N}(\mu_e,\sigma_e)$
is available on day $d$ for the surgery duration of each patient $e\in E_d$. The parameters $\sigma_e$ and $\mu_e$ are generated as above, to match the marginal mean $m_e=90$ and marginal variance $v_e=70^2$.

\begin{landscape}

	\begin{figure}[H]
		\centering
	
		\centering
		\hspace*{-0.85cm}\includegraphics[width=25cm]{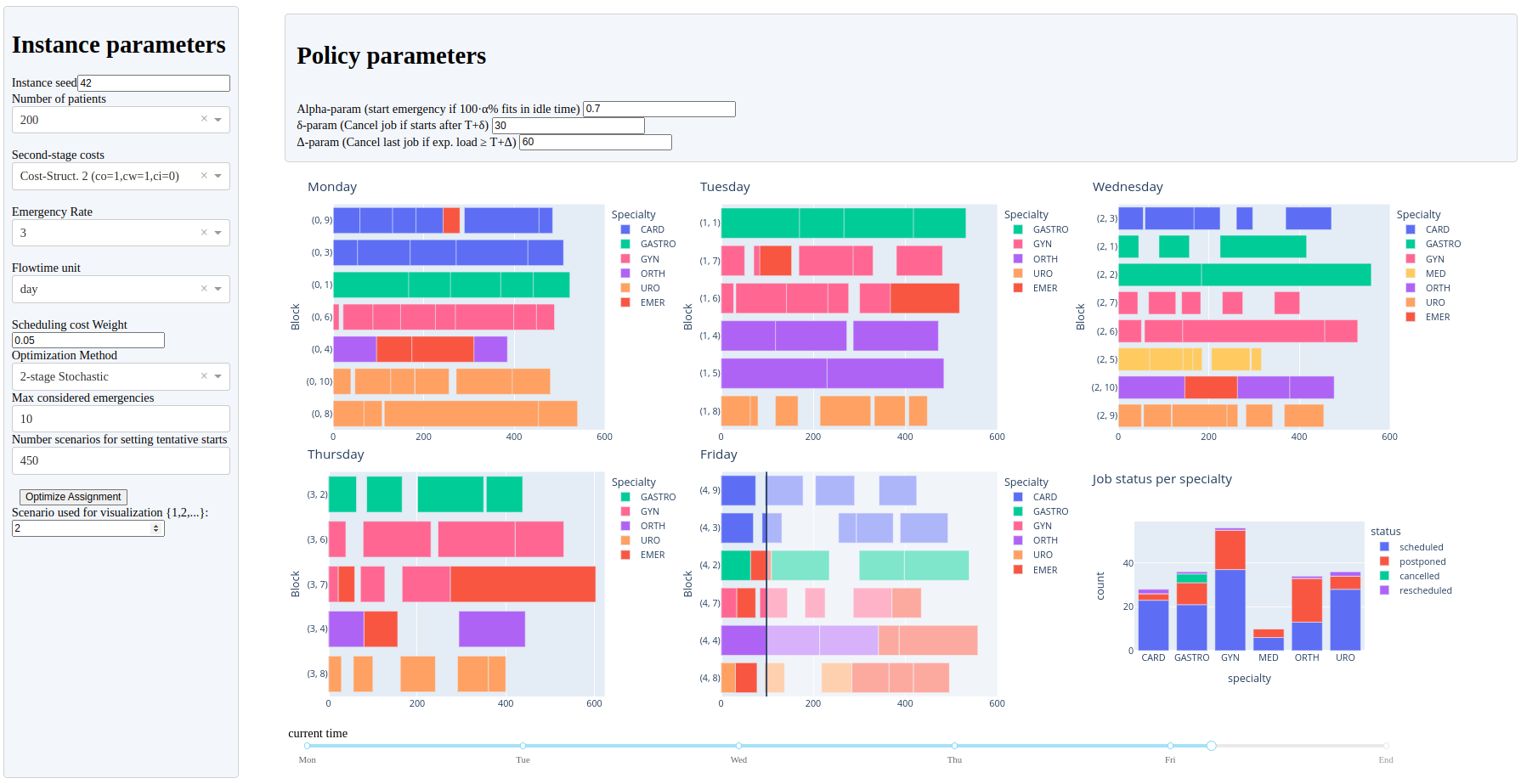}
		\caption{\small Screenshot of the app\label{fig:input}}
	\end{figure}

\begin{figure}[H]
\centering
	\includegraphics[width=25cm]{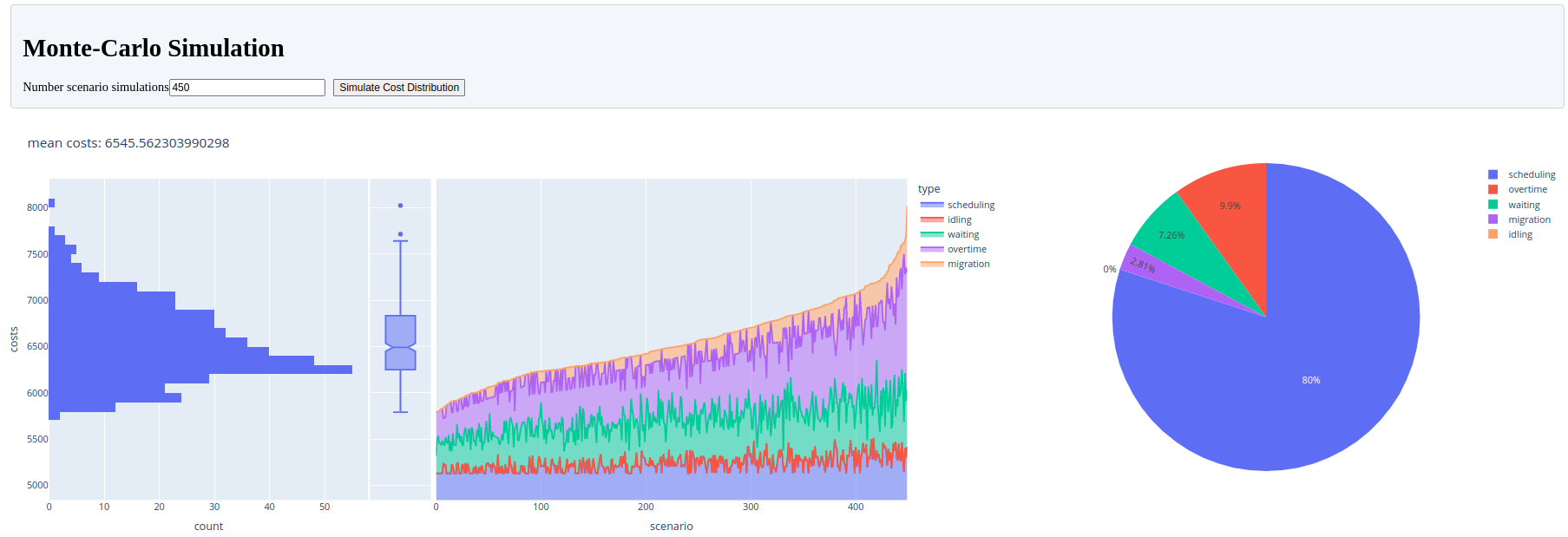}
	\caption{\small Screenshot of the visualization of the Monte-Carlo simulation in the app \label{fig:Monte-Carlo}}
\end{figure}

\end{landscape}

\end{document}